%% file: main_a.tex
\documentclass{agtart_a}
\pdfoutput=1

\usepackage{wasysym}
\usepackage{graphicx}

\title{Legendrian knots and monopoles}

\author{Tomasz Mrowka}
\givenname{Tomasz S}
\surname{Mrowka}
\address{MIT\\77 Massachusetts Avenue\\
Cambridge MA 02139\\USA}
\email{mrowka@math.mit.edu}
\author{Yann Rollin}
\givenname{Yann}
\surname{Rollin}
\address{Imperial College\\
Huxley Building\\
180 Queen's Gate\\
London SW7 2AZ\\UK}
\email{rollin@imperial.ic.ac.uk}

\received{8 November 2005}
\revised{}
\accepted{10 July 2005}
\proposed{}
\seconded{}
\published{24 February 2006}
\publishedonline{24 February 2006}

\volumenumber{6}
\issuenumber{}
\publicationyear{2006}
\papernumber{1}
\startpage{1}
\endpage{69}

\doi{}
\MR{}
\Zbl{}

\subject{primary}{msc2000}{57R17}
\subject{primary}{msc2000}{57M25}
\subject{primary}{msc2000}{57M27}
\subject{primary}{msc2000}{57R57}
\keyword{contact structures}
\keyword{Legendrian knots}
\keyword{Bennequin inequality}
\keyword{excision}
\keyword{monopoles}

\def\overM{\wwbar M}
\def\overZ{\wwbar Z}
\def\overN{\wwbar N}

\makeop{sw}
\makeop{SW}
\makeop{tb}
\makeop{lk}
\makeop{Conn}

\allowdisplaybreaks
\DeclareMathOperator{\id}{id}       %identity
\DeclareMathOperator{\Map}{Map}
\DeclareMathOperator{\vol}{vol}     %Volume
\DeclareMathOperator{\trace}{trace} %trace
\DeclareMathOperator{\End}{End}     %endomorphisms
\DeclareMathOperator{\coker}{coker} %cokernel
\DeclareMathOperator{\Ric}{Ric}     %Ricci curvature
\def\CC{{\mathbb{C}}}
\def\ZZ{{\mathbb{Z}}}
\def\RR{{\mathbb{R}}}
\def\EE{{\mathbb{E}}}
\def\NN{{\mathbb{N}}}
\def\cC{{\mathcal{C}}}
\def\cD{{\mathcal{D}}}
\def\cG{{\mathcal{G}}}
\def\cZ{{\mathcal{Z}}}
\def\cG{{\mathcal{G}}}
\def\cL{{\mathcal{L}}}
\def\cM{{\mathcal{M}}}
\def\cN{{\mathcal{N}}}
\def\cQ{{\mathcal{Q}}}
\def\cS{{\mathcal{S}}}
\def\cW{{\mathcal{W}}}
\def\cU{{\mathcal{U}}}
\def\gG{{\mathfrak{G}}} %\def\gH{{\mathfrak{H}}}
\AtBeginDocument{\renewcommand{\epsilon}{\varepsilon}}
\newcommand{\ip}[1]{\langle #1 \rangle}

\newcommand{\gs}{\mathfrak{s}}
\newcommand{\gtt}{\mathfrak{t}}

\newcommand{\gsu}{\mathfrak{su}}
\newcommand{\spinc}{\ensuremath{\mathrm{Spin}^c}}
\newcommand{\Dirac}{\mathrm{D}}

\newcommand{\del}{\partial}
\newcommand{\delb}{\bar\partial}
\newcommand{\dt}{{\partial_t}}

\newcommand{\BBox}{{ \APLbox}}
\newcommand{\IZ}{\mathbb{Z}}

\numberwithin{equation}{section}
\renewcommand{\theequation}{\thesection.\arabic{equation}}

\renewcommand{\theenumi}{\textup{(\roman{enumi})}}

\makeatletter
\def\cnewtheorem#1[#2]#3{\newtheorem{#1}{#3}[subsection]
\expandafter\let\csname c@#1\endcsname\c@subsubsection}

\def\dnewtheorem#1[#2]#3{\newtheorem{#1}{#3}
\expandafter\let\csname c@#1\endcsname\c@mytheo}

\cnewtheorem{prop}[subsubsection]{Proposition}
\cnewtheorem{cor}[subsubsection]{Corollary}
\cnewtheorem{theo}[subsubsection]{Theorem}
\newtheorem{mytheo}{Theorem}

\dnewtheorem{mycor}[mytheo]{Corollary}
\cnewtheorem{lemma}[subsubsection]{Lemma}
\makeautorefname{mytheo}{Theorem}
\makeautorefname{mycor}{Corollary}

\theoremstyle{definition}
\cnewtheorem{rmk}[subsubsection]{Remark}
\newtheorem*{rmk*}{Remark} 
\newtheorem*{rmks*}{Remarks} 
\cnewtheorem{rmks}[subsubsection]{Remarks}
\newtheorem*{result}{Result} 
\cnewtheorem{dfn}[subsubsection]{Definition}

\makeatother  %  move after \newtheorem block

\newenvironment{remark}{\begin{rmk}}{\end{rmk} }
\newenvironment{remark*}{\begin{rmk*}} { \end{rmk*} }
\newenvironment{remarks*}{\begin{rmks*}$\phantom{99}$ \begin{itemize}}
	       {\end{itemize} \end{rmks*} }
\newenvironment{remarks}{\begin{rmks}$\phantom{99}$  \begin{itemize}}
	       {\end{itemize} \end{rmks} }

\begin{document}

\begin{asciiabstract}
We prove a generalization of Bennequin's inequality for Legendrian
knots in a 3-dimensional contact manifold (Y,\xi), under the
assumption that Y is the boundary of a 4-dimensional manifold M and
the version of Seiberg-Witten invariants introduced by Kronheimer and
Mrowka [Invent. Math.  130 (1997) 209-255] is nonvanishing.  The
proof requires an excision result for Seiberg-Witten moduli spaces;
then the Bennequin inequality becomes a special case of the adjunction
inequality for surfaces lying inside M.
\end{asciiabstract}

\begin{webabstract}
We prove a generalization of Bennequin's inequality for Legendrian
knots in a 3-dimensional contact manifold $(Y,\xi)$, under the
assumption that $Y$ is the boundary of a 4-dimensional manifold $M$ and
the version of Seiberg-Witten invariants introduced by Kronheimer and
Mrowka [Invent. Math.  130 (1997) 209--255] is nonvanishing.  The
proof requires an excision result for Seiberg-Witten moduli spaces;
then the Bennequin inequality becomes a special case of the adjunction
inequality for surfaces lying inside $M$.
\end{webabstract}

\begin{abstract}
We prove a generalization of Bennequin's inequality for Legendrian
knots in a 3--dimensional contact manifold $(Y,\xi)$, under the
assumption that $Y$ is the boundary of a $4$--dimensional manifold $M$
and the version of Seiberg--Witten invariants introduced by Kronheimer
and Mrowka in \cite{KM} is nonvanishing.  The proof requires an
excision result for Seiberg--Witten moduli spaces; then the Bennequin
inequality becomes a special case of the adjunction inequality for
surfaces lying inside $M$.
\end{abstract}

\begin{mathmlabstract}
We prove a generalization of Bennequin's inequality for
Legendrian knots in a 3--dimensional contact manifold
<math xmlns='http://www.w3.org/1998/Math/MathML'>
   <mo>(</mo> <mi>Y</mi><mo>,</mo><mi>&xi;</mi><mo>)</mo>
</math>,        
under the assumption that
<math xmlns='http://www.w3.org/1998/Math/MathML'>
   <mi>Y</mi>              
</math>          
is the boundary of a 4--dimensional manifold 
<math xmlns='http://www.w3.org/1998/Math/MathML'>
    <mi>M</mi>
</math>
and the version of Seiberg--Witten invariants
introduced by Kronheimer and Mrowka in is nonvanishing. The proof
requires an excision result for Seiberg--Witten moduli spaces; then
the Bennequin inequality becomes a special case of the adjunction
inequality for surfaces lying inside
<math xmlns='http://www.w3.org/1998/Math/MathML'>
 <mi>M</mi>
</math>.
\end{mathmlabstract}

\maketitle

\section{Introduction}
This paper is a sequel to the article \cite{KM} by Kronheimer and the
first author,
where the Seiberg--Witten invariants were generalized to invariants of
connected oriented smooth 4--manifolds carrying a contact structure
on their boundary.   An oriented contact structure $\xi$  (or more
generally an oriented 2--plane field)  induces a \emph{canonical}
\spinc--structure $\gs_\xi$ on $Y$. In \cite{KM}, the Seiberg--Witten
invariants were defined for $4$--manifolds with boundary endowed
with a contact structure . The domain of these invariants is the set
$\spinc(\overM,\xi)$ of isomorphism classes of pairs $(\gs,h)$ where
$\gs$ is a \spinc--structure on $\overM$ and $h$ is an isomorphism
between $\gs|_{Y}$ and $\gs_\xi$.  The Seiberg--Witten invariant is a map
$$\sw\co \spinc(\overM,\xi)\to \ZZ$$
well defined up to an overall sign. The main result of this paper is
an excision property for these invariants.  As a corollary we derive
that for overtwisted contact structures these invariants are trivial
and combining this with extension of the Taubes nonvanishing theorem
\cite{Taub} to this context we derive a pseudoholomorphic curve
free proof of Eliashberg's theorem \cite{Etight} that weakly
symplectic fillable contact structures are tight.

\subsection{Some recollections and conventions}
All 4--manifolds will be connected, oriented and smooth unless
otherwise noted.  All 3--manifolds will be oriented and smooth
but not  necessarily connected.

A contact
structure on a compact $3$--dimensional manifold $Y^3$ is
a totally nonintegrable $2$--plane field $\xi\subset TY$ which we
assume to be orientable. Thus there is a nonvanishing  differential
$1$--form $\eta$, such that $\xi=\ker\eta$ and $\eta\wedge d\eta \neq
0$ and so  $Y$ gets an orientation provided by the volume form
$\eta\wedge d\eta>0$.

A knot $K$ in $Y$ is called \emph{Legendrian} if its
 tangent space is contained in the contact plane field $\xi$.
Suppose that $Y$ is the oriented boundary of an oriented
connected smooth 4--manifold $M$ and that $K$ is the
boundary of a smooth connected orientable surface $\Sigma$.
Then $\Sigma$ determines, up to sign, a homology class
\[
\sigma \in H_{2}(M,K;\IZ)
\]
which maps to a generator of $H_{1}(K;\IZ)$ under the coboundary map.

The Thurston--Bennequin invariant and the rotation
number of a knot in $S^{3}$ generalize
to invariants of the pair $K$ and $\sigma$.  The rotation number
is only defined up to sign until an orientation of $K$ is chosen.  Both
invariants arise because a Legendrian knot has a {\it  canonical
framing} obtained by choosing a
vector field $V$  transverse to the contact distribution along $K$.
To generalize the Thurston--Bennequin invariant
choose an arbitrary orientation for $K$ give  $\Sigma$ the compatible
orientation.  Push $K$ off slightly in
the
$V$--direction, thus obtaining a disjoint oriented knot $K'$.
Then push $\Sigma$ off itself to get a surface $\Sigma'$ so that its
boundary
coincides with $K'$. Since the $\Sigma$ and $\Sigma'$ are disjoint
along their boundary they have a well defined intersection number. The
Thurston--Bennequin
invariant relative to $\Sigma$, $\tb(K,\sigma)$ is defined to be this
self-intersection
number.  If $\Sigma$ is contained in $Y$ (so that $K$ is  nullhomologous
in $Y$)
then
$$\tb(K,\sigma) =\lk(K,K'):=\tb(K).$$
Notice that  $\tb(K)$ does not depend on the choice of the initial
orientation for $K$ or $\Sigma$.

The generalization of the rotation
number is obtained as follows. After choosing an orientation, the contact
distribution can be endowed with an almost
complex structure $J_\xi$, which is unique up to homotopy. Therefore
$\xi\to Y$ has the structure of a complex line bundle, hence it has a well
defined Chern class. The isomorphism $h$ induces an isomorphism
of the determinant line $L_{\gs}=\det(W^{+}_{\gs})$ of the bundle
of positive spinors for the \spinc--structure $\gs$ with $\xi$ on the
boundary If we also fix an  orientation for $K$ then we get a preferred
nonvanishing tangent vector field $v$ and so, by the Legendrian property,
a nonvanishing section of $L_{\gs}$ Then the rotation number of $K$
relative to $\Sigma$ is by definition
$$ r(K,\sigma,\gs,h) := \langle c_1(L_{\gs}, v), \sigma\rangle$$
where $c_1(L_{\gs},v)$ is the  relative Chern class with respect to the
trivialization of $ \xi$ along $K$  induced by $v$.
Notice that the rotation number depends \emph{a priori}  on the homology
class of
$\Sigma$ and on the orientation of $K$.
This definition coincides
with the usual rotation number defined on $(\RR^3,\xi_{\mathrm{std}})$
as
the winding
number of $v$ in $\xi_{\mathrm{std}}$.

\subsection{Main results}

With these definitions and notation in place we can state our
generalization
of Bennequin's Inequality.
\begin{mytheo}
\label{theo:main}
Let $(Y,\xi)$ be a  $3$--dimensional closed manifold endowed with a
contact structure $\xi$ and let $\overM$ be a compact
$4$--dimensional manifold with boundary $Y$.
Suppose we have a Legendrian
knot $K \subset Y$,  and a connected, orientable
compact surface $\Sigma\subset \overM$
with boundary
$\partial\Sigma = K$.
Then for every
relative \spinc--structure $(\gs,h)$ with
$\sw_{\overM,\xi}(\gs,h)\not=  0$ we have
$$ \chi(\Sigma) + \tb(K,\sigma) + | r(K,\sigma,\gs,h)| \leq 0,
$$
where $\chi$ denotes the Euler characteristic.
\end{mytheo}
Notice that this result was known before in the case of
compact Stein complex surfaces with pseudoconvex boundary
(see Akbulut--Matveyev and Lisca--Mati{\'c} \cite{AM,LM}).

Here are two corollaries of \fullref{theo:main}.
A contact manifold $(Y,\xi)$ is called \emph{weakly symplectically
fillable}
if it is the boundary of a symplectic manifold $(\overM,\omega)$
such that
$\omega$ is positive on $\xi$.	We say that $(Y,\xi)$ is
\emph{weakly symplectically semi-fillable} if it is a component of a
weakly
symplectically fillable contact manifold.
 A contact structure $\xi$ on $Y$ is called \emph{overtwisted} if $Y$
contains
an embedded disk $D^2$ whose boundary $K=\partial D^2$ is tangent to
$\xi$ while $D^2$ is transverse to $\xi$ at the boundary --- such a
disk is called an \emph{overtwisted disk}. Otherwise we
say that the contact structure is \emph{tight}.
If $\xi$ is overtwisted, we see that the Bennequin inequality does
not hold
 since $\tb(\partial D^2) = 0 $ and $\chi(D^2) = 1$. So we have the
 following result.
\begin{mycor}
\label{mycortight}
The Seiberg--Witten invariant $\sw_{\overM,\xi}$ is identically
zero for  a manifold with an
overtwisted contact boundary.
\end{mycor}

Using the fact that $\sw_{\overM,\xi}\not\equiv 0$ for	weakly
symplectically semi-fillable contact structures \cite[Theorem
1.1]{KM},
we have a pseudoholomorphic curve free proof of the
the following theorem of Eliashberg's theorem result without using
curves.
\begin{mycor}
If $(Y^3,\xi)$ is a weakly symplectically semi-fillable contact
manifold, it
is necessarily tight.
\end{mycor}
This result could also be deduced from Taubes nonvanishing theorem
and Eliashberg's recent result on concave fillings \cite{Efill}.

These results highlight an interesting question.  Are there contact
3--manifolds $(Y,\xi)$ which are not weakly symplectically fillable
and yet there is a 4--manifold $M$ which bounds $Y$ so that
$\sw_{\overM,\xi}$ is not identically zero?  On the other hand
there are contact structures which are tight but not weakly symplectically
fillable.  The first examples are due to Etnyre and Honda \cite{EH}
on certain
Seifert fibered space and later infinite families were discovered by
Lisca and Stipsicz \cite{LS1,LS}.

All these results rely on an excision property for Seiberg--Witten
invariants. Recall that a \emph{symplectic cobordism} between contact
manifolds $(Y,\xi)$ and $(Y',\xi')$ is a compact symplectic manifold
$(\overZ,\omega) $  so that,  with the symplectic orientation,
its  boundary is $\del \overZ = -Y\sqcup
Y'$, where $Y$ and $Y'$ have their orientations induced by the contact
structures. $Y$ is called the {\it concave} end of the cobordism
and $Y'$ is called the {\it convex} end. In addition, it is required
that $\omega$ is strictly
positive on
$\xi$ and $\xi'$ with their induced orientations.
By convention the boundary
components will always be given in the order concave, convex.

A symplectic cobordism is said to be \emph{special} if
\begin{itemize}
\item the symplectic form
is given in a collar neighborhood of the concave boundary by
 a symplectization of $(Y,\xi)$;
\item the map induced by the
  inclusion
\begin{equation}
\label{assumption2}
i^*\co H^1(\overZ, Y')\to H^1(Y)
\end{equation}
is the zero map.
\end{itemize}

A symplectic cobordism carries a canonical
\spinc--structures $\gs_{\omega}$. Moreover, there are isomorphisms,
unique up to
homotopy, which identify the restrictions of $\gs_\omega$ to $Y$ and $Y'$
with the canonical \spinc--structures $\gs_\xi$ and  $\gs_{\xi'}$.

Let $(\gs,h)$ be an element in $\spinc(\overM,\xi)$. There is
a canonical way to extend $(\gs,h)$ to a $\spinc$ structure $\gtt$ on
$\overM'=\overM\cup
\overZ$ together with an isomorphism $h'$ between $\gtt|_{Y'}$ and
$\gs_{\xi'}$. We declare $\gtt:=\gs_\omega$ on $\overZ$. The data of $h$
identifies the restriction  $\gs|_{Y}$ with
$\gs_{\xi}$ while $\gs_\omega|_{Y}$ is identified canonically with
$\gs_\xi$. Together these provides a gluing map and defines $\gtt$.
Thus we have defined  a
 canonical map
$$j\co\spinc(\overM,\xi)\to \spinc(\overM',\xi').$$
The main technical result of this paper is the following.
\begin{mytheo}
\label{theoexcision}
Let $\overM$ be a manifold with contact boundary $(Y,\xi)$, and let $\overZ$
be a special symplectic cobordism between
$(Y,\xi)$ and a second contact manifold $(Y',\xi')$. Let
$\overM'=\overM\cup_Y \overZ$
 be the manifold obtained by gluing $\overZ$ on $\overM$
 along $Y$. Then
  \begin{equation}
\label{eqexcision} \sw_{\overM,\xi}\circ
 j = \pm \sw_{\overM',\xi'}\, .
\end{equation}
\end{mytheo}

\begin{remark}
 Recall the invariants themselves are only defined up to an overall
 sign.	We can refine the statement above to the following. For every
 pair of
 \spinc--structures $(\gs_1,h_1)$
 and $(\gs_2,h_2)$ in $\spinc(\overM , \xi)$ we have:
$$ \sw_{\overM',\xi'}\circ j(\gs_1,h_1)\sw_{\overM,\xi}(\gs_2,h_2) =
\sw_{\overM,\xi}(\gs_1,h_1)\sw_{\overM',\xi'}\circ
 j  (\gs_2,h_2)\,.
$$
\end{remark}
 \begin{remark}
\label{rkassumption}
Assumption~\eqref{assumption2} can be reformulated as follows.  If 
$u\in\Map(\overZ,S^1)$ is homotopic to the identity
along $Y'$, then
$u$ must be homotopic to the identity along $Y$. Without the
assumption~\eqref{assumption2}, the gauge transformation $u|_Y$ may
not extend
to $\overM$. In this case the map $j$ is not generally injective anymore.

All the cobordisms of
interest in this paper
will be shown to verify assumption \eqref{assumption2}, that is to
say $1$ and $2$--handle surgeries.
Assumption \eqref{assumption2} may be removed, but
conclusion~\eqref{eqexcision} of  \fullref{theoexcision} has to be
replaced by
   \begin{equation}
\label{eqexcision2} \sw_{\overM',\xi'} (\gs',h')= \pm\!\!\!\!\sum_{(\gs,h)\in
  j^{-1}(\gs',h')}\!\!\!\! \sw_{\overM,\xi}(\gs,h)
   \, .
 \end{equation}
This generalization is proved by refining the gluing
\fullref{theogluing} as explained in \fullref{rkglue}.	Indeed if
two pairs $(s_{1},h_{1})$ and $(s_{2},h_{2})$ have the same image
under $j$ then we can assume, up to isomorphism, that $s_{1}=s_{2}$ and
$h_{1}=uh_{2}$ where $u$ is an automorphism of $\gs_{\xi}$ which extends
to an automorphism of $\gs_{\omega}$ which is the indentity at
infinity.
 \end{remark}

  A result of
Weinstein~\cite{W} shows that a $1$--handle surgery, or a $2$--handle
surgery along a Legendrian knot $K$ with
framing coefficient $-1$ relative to the canonical framing, on the
boundary of $\overM$
leads to a manifold $\overM'$  given
by
$$\overM' = \overM\cup_Y \overZ$$
 where $(\overZ,\omega)$ is a special symplectic cobordism
 between $(Y,\xi)$ and a  contact boundary $(Y',\xi')$ obtained by
 the surgery.

The strategy to prove \fullref{theo:main} is to cap $\Sigma$
doing a Weinstein surgery along $K$. Bennequin's
 inequality is then obtained by applying the adjunction
inequality \cite{KM2}  to the resulting closed surface	$\Sigma'\subset
\overM'$; however the adjunction inequality
holds provided	$\sw_{\overM',\xi'}\circ j(\gs,h)$ does not vanish.
This is true under the assumption that $\sw_{\overM,\xi}(\gs,h)$
does not vanish
thanks to the
 excision \fullref{theoexcision}.

\subsection{Monopoles and contact structures}
We  briefly describe the Seiberg--Witten theory on a 4--manifold $\overM$
with
contact boundary $(Y,\xi)$ developed by Kronheimer and
the first author  in \cite{KM}.

\subsubsection{An almost K\"ahler cone}
\label{subsubkahler}
We pick a contact 1--form $\eta$ with $\ker \eta =\xi$.  The
contact form determines the  Reeb vector field $R$ by the
properties that $\iota_{R}d\eta =0$ and $\eta(R)=1$.
 Then $(0,+\infty) \times Y $ is endowed with a symplectic structure
    called the
    symplectization of $(Y,\eta)$, and the symplectic form is defined
by
the formula
$$\omega =\tfrac 12d(t^2\eta).$$
Choose an almost complex structure $J$ on $(0,\infty) \times Y$
where $\xi$ and $d\eta$  are invariant under $J$, and such that
$$tJ\dt = R.$$
The almost complex structure $J$ is clearly compatible with the symplectic
form $\omega$ in
the sense that $g=\omega(J\cdot,\cdot)$ is a Riemannian metric and
that $\omega$ is $J$ invariant; thus
we have defined an
almost K\"ahler structure $(g,\omega,J)$ on $(0,\infty)\times Y$
called the \emph{almost K\"ahler cone} on $(Y,\eta, J|_\xi)$.

The metric $g$ has an  expression of the form
$$g=dt^2 + t^2\eta^2 +t^2 \gamma,$$
where $\gamma$ is the positive symmetric bilinear form on $Y$ defined
by
$\gamma = \frac 12 d\eta(J \cdot, \cdot)$.

Let $M$ be  the  manifold
$$ M =\overM \cup C_M,$$
where $C_M = (T, \infty)\times Y$ is endowed with the structure of
    almost K\"ahler cone on $(Y,\eta, J|_\xi)$ described above.
Then we arbitrarily extend the Riemannian metric $g$ defined on the
end $C_M$ to the compact set $ \overM$.

The \spinc--structure $\gs_\xi$ induced by the contact
structure
on $Y$ is canonically identified with the restriction to
$\{T\}\times Y$ of the \spinc--structure $\gs_\omega$ induced by
the
symplectic form $\omega$ on $C_M$.  Therefore, an element
$(\gs,h)\in \spinc(\overM,\xi)$ admits a natural extension over $M$. By a
slight abuse of language, we will still denote it $(\gs,h)$, where $h$
is now an isomorphism between $\gs$ and $\gs_\omega$ over $C_M$.

The spinor bundle $W= W^+\oplus W^-$ of $\gs$ is identified over
$C_M$ with  the
canonical Spin--bundle $W_{J}=W^+_J\oplus W_J^-$ of $\gs_\omega$;
we recall that
$$ W^+_J =\Lambda^{0,0} \oplus \Lambda^{0,2},\quad W^-_J = \Lambda^{0,1}.$$

\subsubsection{Monopole equations}
\label{subsubmonopoleeq}
 Let $A$ be a spin connection on $W$ and $\nabla^A$ be its covariant
 derivative. The
 the Dirac operator is then defined by
$\Dirac_A^g = \sum_i e^{i}\cdot \nabla^A_{e_j}$, where $e_j$ is an
 oriented orthonormal local frame on $M$. $e^{i}$ is the
 dual coframe and acts by Clifford
 multiplication.
It is a first order elliptic operator of order $1$ between the space
 of spinor fields $\Dirac_A\co\Gamma(W^\pm)\rightarrow
\Gamma(W^\mp)$.

Put $\Psi=(1,0)\in \Lambda^{0,0}\oplus \Lambda^{0,2}$. Let $B$ be the
spin-connection in
the spinor bundle $W^+_J = \Lambda^{0,0}\oplus \Lambda^{0,2}$
so that $ \nabla^B \Psi$ is a section of $T^*X\otimes \Lambda^{0,2}$.

The Dirac operator $\Dirac^{can}$ on $W_J$ is obtained by replacing
$\nabla^A$ with $\what\nabla$ in the definition. For $\Phi\in
\Gamma(W_J^+)$, we
have explicitly $\Dirac^{can}\Phi= \sqrt 2(\delb\beta + \delb^*\gamma)$
for $\Phi=(\beta,\gamma)\in
\Gamma(\Lambda^{0,0}\oplus\Lambda^{0,2})$. Then
$\Dirac^{can}\Psi = 0$.

The Chern connection $\what \nabla$ in $\Lambda^{0,0}\oplus\Lambda^{0,2}$
induces a connection in $\Lambda^{0,2} =\det(W_J^+)$ and this connection
agrees with the one induced from $B$.
There is also  an  identity
$$\Dirac_B = \Dirac^{can}.$$
Since the spinor bundle $W$  is identified with $W_J$ via the
isomorphism $h$ over $C_M$,
 $B$ may be considered as a spin connection over $W|_{C_M}$ and
 $\Psi|_{C_M}$ as a section of $W^+|_{C_M}$.  We extend them arbitrarily
 over
 $\overM$.

The domain of the Seiberg--Witten equations on $(M,g,\gs)$ is
the space of
configurations
$$\cC = \{(A,\Phi)\in \Conn(W)\times \Gamma(W^+)\},$$
where $\Conn(W)$ is the space of spin connections
on $W$. Recall that $\Conn(W)$ is an affine space modeled on
$\Gamma(i\Lambda^1)$: for two spin connections $A$ and $\wtilde A$, we
have
$$ \wtilde A = A + a\otimes \id|_W,
$$
where $a$ is a purely imaginary $1$--form. We will simply write the
above identity $\wtilde A = A +a$ in the sequel.
 We introduce the curvature form
\begin{equation}
\label{deffa}
 F_A (X,Y) = \tfrac 14 \trace_\CC R_A(X,Y),
\end{equation}
where the full curvature tensor $R_A$ of $A$ is viewed as a section of
$\Lambda^2\otimes \End_\CC(W)$. With this convention, we have
$$ F_A=\tfrac 12 F_{\what A}
$$
where $\what A$ is the unitary connection induced by $A$ on the determinant
line bundle\footnote{Alternatively $F_A$ can be viewed as  the curvature
of
the unitary connection induced by $A$ on the \emph{virtual line
  bundle}~$L^{1/2}$}
$L$, and $F_{\what A}$ is its usual curvature form.
The Seiberg--Witten equations are
\begin{eqnarray}
\label{eqSW1}
F_A^+- \{\Phi\otimes\Phi^*\}_0 & = & F_B^+- \{\Psi\otimes\Psi^*\}_0 +
\varpi\\
\Dirac_A \Phi & = & 0,
\end{eqnarray}
where
$\varpi$ is a self-dual purely imaginary $2$--form on $M$ and
$\{\cdot\}_0$ represents the trace free part of an endomorphism;
the first equation makes sense since  purely imaginary self-dual
 $2$--forms are identified with the traceless endomorphisms of $W^+$ via
 Clifford multiplication.
The configuration $(B,\Psi)$ is clearly a solution of  the equations
 over $C_M$.

 The space of solutions $\cZ_{\varpi}$ is acted on by a gauge group
 $\cG=\Map(M,S^1)$ and the action is defined on $\cC$ by
\begin{equation}
 u\cdot(A,\Phi) = (A - u^{-1}{du}, u \Phi)\quad \text{for all } u\in\cG .
\end{equation}
The moduli space $\cM_{\varpi}(M,g,\gs)=\cZ_{\varpi}/\cG$, for suitable
generic
$\varpi$, is a compact smooth manifold of
dimension
$$ d= \ip{e(W^+,\Psi),[M,C_ M]},
$$
where $e(W^+,\Psi)$ is the relative Euler class of $W^+$.
If $d\neq 0$, then the Seiberg--Witten invariant is always  $0$; if
$d=0$, the Seiberg--Witten invariant is the
number of points  of $\cM_{\varpi}(M,\gs)$  counted with signs. Following
\cite{KM} there is a trivial bundle over the configuration space (the
determinant line bundle of the appropriate deformation operator)  which
is identified with the orientation bundle of moduli space. Furthermore
a trivialization of this determinant line bundle for one
relative \spinc--structure determines a trivialization
for all others in a canonical manner and determines  a {\it consistent
orientation}.  Thus the
set of consistent orientations is a two-element set.  In particular
unlike the closed case the sign of the invariant
cannot be pined down by a homology orientation rather the ratio of
the signs
of the	values of the invariant for different relative \spinc--structures
is well defined.  It turns out
that this number  depends only on $\overM$, on the contact
structure at the boundary, and on the choice of $(\gs,h)\in
\spinc(\overM,\xi)$; this explains the notation $\sw_{\overM,\xi}$.

\begin{remark*} The
details concerning the regularity of $\cC$ and $\cG$ have been omitted
at the moment; morally, $(A,\Phi)$ is
supposed to behave like $(B,\Psi)$ near infinity and gauge
transformations should be close to identity as well. The relevant
Sobolev
spaces will be introduced later on.
\end{remark*}

 \rk{Acknowledgments} 
This paper obviously owes a large intellectual debt to Peter
 Kronheimer. The second author thanks the MIT where most of this work
 was done. He also thanks the Institute for Advanced Studies where he
 was hosted on several occasions.  The first author was partially
 supported by NSF grants DMS-0111298, DMS0206485, and FRG-0244663. The
 second author was partially supported by NSF grant DMS-0305130.

\section{Excision}
\label{secexi}
The goal of this section is to prove \fullref{theoexcision}. The
strategy is to show, in a more general setting, that the Seiberg--Witten
moduli spaces associated
to $\overM$ and $\overM'$ are diffeomorphic for a suitable
choice of metrics and perturbations.  This is achieved thanks to a gluing
technique in \fullref{theogluing}.

\subsection{Families of AFAK manifolds}
\label{sec:paste}
We set up the analytical framework for	our gluing problem. The
first
step is to construct suitable families of asymptotically flat almost
K\"ahler metrics.

Let $\overM$ be a manifold with a contact boundary $(Y,\xi)$. We
choose a particular contact form $\eta$ and glue an almost K\"ahler
cone on the boundary as in \fullref{subsubkahler}; furthermore the
Riemannian metric is extended arbitrarily to $\overM$.
Hence we have obtained
   $(M,g,\omega,J)$ where $g$ is Riemannian metric and $M$  splits as
\begin{equation}
\label{symplcone}
M=\overM \cup_Y
 C_M
\end{equation}
with an almost K\"ahler structure defined outside a compact set.
More  generally, we recall the definition~\cite[Condition 3.1]{KM} of an
\emph{asymptotically flat
 almost K\"ahler}  manifold (AFAK manifold in short).
\begin{dfn}[AFAK manifolds]
\label{dfnafak}
A  manifold $M$ with an almost K\"ahler structure
$(M,\omega,J)$ defined outside a compact
set $K\subset M$, a Riemannian metric $g$ extending the metric on the
end, and a proper function
$\sigma\co M\rightarrow (0,\infty)$ is called an asymptotically flat
    almost K\"ahler manifold   if it satisfies the
    following conditions:
\begin{enumerate}
\item There is a constant $\kappa>0$, such that the injectivity radius
  satisfies $\kappa \; inj(x) > \sigma(x)$ for all  $x\in
M$.
\item For each $x\in M$, let $e_x$ be the map $e_x\co v \mapsto
  \exp_x(\sigma(x) v/\kappa )$ and $\gamma_{x}$ be the metric on the
unit
  ball in $T_x M $ defined as $e_x^* g/\sigma(x)^2$. Then
  these metrics have bounded geometry in the sense that all covariant
  derivatives of the curvature are bounded by some constants
  independent of~$x$.
\item For each $x\in M\setminus K$, let $o_{x}$ similarly be the
  symplectic form  $e_z^*\omega/\sigma(z)^2$ on the unit
  ball. Then $o_{x}$ similarly approximates the
  translation-invariant form, along
  with all its	derivatives.
\item For all $\epsilon >0$, the function $e^{-\epsilon\sigma}$
  is integrable on $M$.
\item The symplectic form $\omega$ extends as a closed form on
  $K$.
\end{enumerate}
\end{dfn}

{\bf Convention}\qua
When  the end of $M$ has a structure of symplectic cone as
 in~\eqref{symplcone} we choose
 $\sigma(t,y)=t$ on the end $C_M$ and extend it arbitrarily
     with the condition $0<\sigma(x)\leq T$ on	$\overM$.

\begin{dfn}[AFAK ends]
\label{dfncob}
An asymptotically flat almost K\"ahler end is a manifold $Z$ which
admits a decomposition
$C_Z \cup_Y N$,
 where
 $N$ is a not necessarily compact $4$--dimensional manifold, with
 a contact
 boundary $Y$, endowed with a fixed contact form $\eta$, and $C_Z =
 (0,T]\times Y$ for some $T>0$.

In addition, $Z$ is endowed with an almost K\"ahler structure
$(\omega_Z,J_Z)$
and a proper function $\sigma_Z\co Z\rightarrow (0,\infty)$ satisfying:
\begin{itemize}
\item on $(0,T]\times Y \subset Z$, we have $\sigma_Z(t,y)=t$;
\item the almost K\"ahler structure on $C_Z$ is the the one of an almost
  K\"ahler cone on $(Y,\eta)$ as defined in \fullref{subsubkahler};
\item if we substitute	$M$ by $Z$, $K$ by $\emptyset$ and $\sigma $
  by $\sigma_Z$, in \fullref{dfnafak}, the the properties
  $(i)-(iv)$ are verified;
\item the map induced by the inclusion $Y=\del N \subset N$
\begin{equation}
\label{assumption}
H^1_c(N)\longrightarrow H^1(Y)
\end{equation}
where $H^*_c(N)$ is the compactly supported DeRham cohomology,
is identically $0$.
\end{itemize}
\end{dfn}
\begin{remarks}
\label{rqsdef}
\item The last condition is a technical condition, reminiscent
  of the assumption~\eqref{assumption2}
  (see \fullref{rkassumption}). More precisely, we will construct
  some particular AFAK ends out of special symplectic cobordisms in this
  section, and the
  property~\eqref{assumption} will be a consequence of the
  assumption~\eqref{assumption2}. This property will be used for the
  compactness results of \fullref{subsubcomp}. However this
  assumption is not essential as explained in \fullref{remarkglue}.
\item
If we scale the metric by $\lambda^2$, and accordingly
scale $T ,\sigma$ into $\lambda T, \lambda \sigma$, then the constant
$\kappa$ and the constants controlling $\gamma_x$ and $o_x$ remain
unchanged.
In particular, we may always assume $T=1$ in
Definitions~\ref{dfnafak} and \ref{dfncob}.
\item In \fullref{dfncob}, the contact form $\eta$ can be
  replaced by any other $1$--form representing the same contact
  structure. This is shown in the next lemma which allows us to modify
  the symplectic form $\omega_Z$ near the sharp end of the cone $C_Z$.
\end{remarks}

\begin{lemma}
Let $Y$ be a compact $3$--manifold endowed with a contact structure
$\xi$. Let $\eta_1$ and $\eta_2$ be two $1$--forms such that
$\ker\eta_j=\xi$.
Then, for every $\epsilon>0$, there exist $\alpha\in(0,\epsilon)$ and
a symplectic form $\omega$ on $(0,+\infty)\times Y$ such that
\begin{itemize}
\item $\omega = \frac 12 d(t^2\eta_1)$ on  $(0,\alpha)\times Y$,
\item $\omega = \frac 12 d(t^2\eta_2)$ on $(\epsilon ,+\infty)\times
  Y$.
\end{itemize}
\end{lemma}
\begin{proof}
We write $\eta_2 = e^{\mu}\eta_1$ where $\mu$ is a real function on
$Y$. We consider the exact $2$--form
$$ \omega = \frac 12 d(t^2e^{\mu f(t)}\eta_1),
$$
where $f(t)$ is a smooth increasing function equals to $1$ for $t\geq
2\epsilon$ and to $0$ for $t\leq \epsilon$. If we show that we can
choose $f$ in such a way that $\omega$ is definite, then, $\omega$ is
a symplectic form and is a solution for the lemma.

A direct computation shows that
$$ \omega^2 = \frac 12 t^3(2+t\mu f')dt\wedge\eta\wedge d\eta.
$$
Let $c$ be the minimum of the function $\mu$ on $Y$. If $c\geq 0$, we
just
require that $f$ is an increasing function of $t$. If $c<0$, a
sufficient additional condition for having $\omega^2>0$ is
\begin{equation}
\label{eqcond000}
 f'(t) < -\frac 2{tc}.
\end{equation}
The function $f_0 = -\frac 1c \ln (\frac {2et}{\epsilon} )$ verifies the
above identity. Then we define
\begin{itemize}
\item $ f_1(t)=0$ for $t\in(0,\frac\epsilon{2e}]$,
\item $f_1(t)=f_0(t)$ for $t\in [\frac \epsilon {2e},\frac \epsilon 2]$,
\item $f_1(t) = 1$ for $t\in [\frac \epsilon 2,+\infty)$.
\end{itemize}
By definition $f_1$ is a piecewise $C^1$ function which verifies the
condition~\eqref{eqcond000} on each interval. We can regularize $f_1$
into a	positive smooth increasing function $f$ by making a perturbation
on the
interval $[\frac \epsilon {4e}, \epsilon]$ in such a way that we have
$f'\leq
f'_1$ on each interval where $f'_1$ is defined. Therefore, the
condition $f'<	-\frac 2c t$ is preserved hence $\omega^2>0$.
\end{proof}

\subsubsection{Gluing an AFAK end on $\overM$}
\label{secglueafak}
Thanks to \fullref{rqsdef}, we may assume from now on that $T=1$ and that
the contact form $\eta$ is the same in the definition of the almost
K\"ahler  cone
$C_M=(T,\infty)\times Y$ and in the definition of the AFAK end
(\fullref{dfncob}).
We identify an annulus in $C_M\subset M$ with an annulus in $C_Z
\subset Z$ using the dilation map
\begin{align*}
  M  \supset C_M\simeq	(1,+\infty)\times Y \supset   (1,
\tau)\times Y  & \stackrel{\nu^\tau}\longmapsto
   (1/\tau, 1 )\times Y    \subset (0,1) \times Y \simeq C_Z \subset
Z\\
  (t,y)&\longmapsto (t/\tau,y)
\end{align*}
and  define the manifold $M_\tau$ as the union of
$M\cap \{  \sigma_M <\tau\}$ and $ Z\cap \{ \sigma_Z> 1/\tau\}$
and with the identify along the annuli given by the dilation
$\nu^\tau$.  The operation of connected sum along $Y$ we just defined
is represented in the figure below. The gray regions represent the
annuli, the arrows suggest that they are identified by a dilation,
and the dashed regions are the parts of $M$ and $Z$ that are taken
off for the construction of $M_\tau$.

 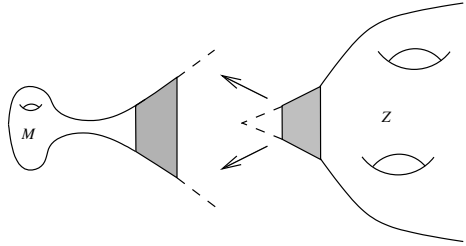
\begin{figure}[ht!]
 \begin{center}
 \input{\figdir/MZ.pstex_t}
 \caption{Construction of AFAK $M_\tau$}
 \end{center}
 \end{figure}

Now $\nu^\tau_* \omega_{C_M} = \tau^2 \omega_Z$ and
$\sigma_Z\circ\nu^\tau= \sigma_M /\tau$, hence, if we scale
$\omega_Z$ by $\tau^2$ and $\sigma_Z$ by $\tau$, all the structure
will match on the annuli. In conclusion $M_\tau$ carries an almost
K\"ahler  structure
$(\omega_\tau,J_\tau)$
 defined outside the compact set $\overM\subset M_\tau$
 and  functions $\sigma_\tau$.

\begin{remark*}
Every compact set of $K\subset M$ is also, by definition, a compact
set
of $M_\tau$ provided $\tau$ is large enough. Similarly, the
structures
$g_\tau, \sigma_\tau, J_\tau, \omega_\tau$ are equal on  every compact
set when $\tau$ is large enough.
\end{remark*}

The following lemma is satisfied by construction.
\begin{lemma}
\label{lemmaboundedgeometry}
The manifolds with almost K\"ahler structure defined outside a compact
set and a proper function
$(M_\tau,g_\tau,J_\tau,\sigma_\tau)$ satisfy \fullref{dfnafak}
uniformly, in the sense that the constant $\kappa, \epsilon$,  the
bounds
on $\gamma_x$, $o_x$ and on the integral of $e^{-\epsilon\sigma_\tau}$
can be chosen independently of $\tau$.
\end{lemma}

\begin{remarks*}
\item A simple consequence of the lemma is the following: for all
  $\epsilon>0$, there exists $T_k$ large enough such that, for every
  $\tau$, the pull-back of the almost K\"ahler structure on a unit
  balls   in
$M_\tau \cap\{\sigma_\tau\geq T_k \}$,	 via the exponential map
  $v\to\exp_x v$, on
  the tangent space $T_xM_\tau $
,  is $\epsilon$--close in $C^k$--norm to the
  euclidean structure $(g_\tau ,\omega_\tau, J_\tau)|_{T_xM_\tau}$.

\item
For all the remainder of \fullref{secexi}, $M$ and $M_\tau$ will
denote the
the manifolds that we just constructed at \fullref{secglueafak}
together with their additional structures (identification of $C_M$
with an almost K\"ahler cone,
proper function $\sigma$, almost K\"ahler structure and
Riemannian
metric.
\end{remarks*}

\subsubsection{\spinc--structures on AFAK manifolds}
\label{secextspinc}
Similarly to the case of a manifold with a contact boundary, we define
the space $\spinc(X,\omega)$ for an AFAK manifold $X$ as the set of
equivalence
classes of pairs $(\gs,h)$, where $\gs$ is a \spinc--structure on $X$
and $h$ is an isomorphism, defined outside a compact set $K_1\subset X$,
between
$\gs_\omega|_{X\setminus K_1}$ and $\gs |_{X\setminus K_1}$. As we saw in
\fullref{subsubkahler}, there is a well defined identification
of $\spinc(\overM,\xi)$ with $ \spinc(M,\omega)$; more generally, there
is a natural identification $ j \co  \spinc(\overM,\xi)\to
\spinc(M_\tau,\omega_\tau)$ when $M_\tau$ is obtained by adding an
 AFAK end $Z$ to the end of $\overM$. Notice that the set of
 equivalence classes $\spinc (M_\tau,\omega_\tau)$ do not depend on
 $\tau$.
However the realization of the \spinc--structure by a spinor bundle is
sensitive to the choice of $\tau$: let $(\gs,h)\in \spinc( M,\omega
)$, and suppose that $h$ is defined on $M\cap\{\sigma > 1\}$ for
simplicity. Alternatively, $h$ can be though of as an isomorphism between
the spinor bundle $W$ of $\gs$	and the spinor bundle
$W_{J}$ of $\gs_\omega$.
We define a family of spinor bundles $W_\tau$ on
 $M_\tau$ by:
$$\left \{
\begin{array}{ll}
W_\tau:= W &  \quad  \mbox{over } M_\tau\cap \{\sigma_\tau <\tau
\}\subset
M \\
W_\tau:=W_{J_\tau} & \quad \mbox{over }  M_\tau\cap \{\sigma_\tau > 1
\}\subset
 Z,
\end{array}
\right .
$$
and the transition map from $W$ to $W_{J_\tau}$ is given by $h$ over the
annulus
$\{1< \sigma_\tau <\tau\}\cap M_\tau \subset C_M$ (where  $W_{J_\tau}=
W_{J}$). At the level of equivalence classes of \spinc--structures, this
procedure  defines an identification of $\spinc(M,\omega)$ with
$\spinc(M_\tau,\omega_\tau)$.

We stress  the fact that  $W_\tau$ is  identified with
the canonical
spinor bundle $W_{J_\tau}$ for it is, by construction, equal to it on
$M_\tau \cap\{\sigma_\tau >1\}$.
Moreover, the spinor bundles $W_\tau$
restricted to any compact set  $K\subset M$ are all identified
provided $\tau$ is
large enough.

In \fullref{subsubmonopoleeq} we defined a configuration
$(B,\Psi)$ for	the spinor bundle $W\to M$. Similarly, we can
define a configuration $(B_\tau,\Psi_\tau)$ for the spinor
bundle $W_\tau\to
M_\tau$: outside $\overM$, the bundle $W_\tau$ is identified with
the canonical spinor bundle $W_{J_\tau}$. The almost K\"ahler
structure induces a Chern connection $\what \nabla_\tau$ on
$L_{J_\tau}=\det
W_{J_\tau}^+=\Lambda^{0,2}$. Hence we deduce a spin connection
$B_\tau$ from the Levi--Civita connection of $g_\tau$ and
$\what\nabla_\tau$.
Let  $\Psi_\tau =(1,0)$ be the spinor of
$\Lambda^{0,0}\oplus\Lambda^{0,2}=W^+_{J_\tau}$.
Then we choose a common extension (independent of $\tau$) of the
configuration
$(B_\tau,\Psi_\tau)$ to the compact $\overM \subset M_\tau$.

In conclusion we have constructed a family of spinor bundles
$W_\tau\to M_\tau$ which are identified with  $W_{J_\tau}$ outside $\overM$
and which are  realizations of the same element $j(\gs,h)\in
\spinc (M_\tau,\omega_\tau)$.

On a
 compact set $K\subset M$, $W$ is canonically identified with
 $W_\tau$ for every
 $\tau$ large enough and the configurations $(B_\tau,\Psi_\tau)$
agree; hence we
will write  $W$ and
 $(B, \Psi)$ instead of $ W_\tau$, $(B_\tau,
\Psi_\tau)$ for simplicity
 of notation.

\subsection{A family of Seiberg--Witten moduli spaces}
\label{sub:approx}
We introduce now Seiberg--Witten equations for the
AFAK manifolds $M$ and $M_\tau$. We show that
in some sense the moduli space	of Seiberg--Witten equations on
$M$  is
a limit of the moduli spaces on $M_\tau$ as $\tau\to +\infty$.

\subsubsection{A family of Seiberg--Witten equations}
\label{secSWeq}
Starting from an element $(\gs,h)\in \spinc(\overM,\xi)$, we
consider the \spinc--structure induced on $M$ and
$j(\gs,h)$ on $M_\tau$.
The Seiberg--Witten equations were introduced in
\fullref{subsubmonopoleeq} on $M$, which has an end modeled on an almost
K\"ahler cone $C_M$. The equations are given on $M_\tau$ in the same
way by
\begin{eqnarray}
\label{eqSW2}
F_A^+- \{\Phi\otimes\Phi^*\}_0 & = & F_B^+- \{\Psi\otimes\Psi^*\}_0 +
\varpi_\tau\\
\Dirac_A \Phi & = & 0 ,
\end{eqnarray}
where $\Phi$ is a section of $W_\tau^+$, $A$ is a spin
connection on $W_\tau$ and $\varpi_\tau$
is a perturbation in $\Gamma(i\Lambda^+M_\tau)$.

As in the case of an almost K\"ahler conical end,
 the almost K\"ahler
structure defined outside $\overM$
induces a Chern connection $\what\nabla$ on $W_\tau$,  with
a  corresponding canonical Dirac operator
$ \Dirac^{can}= \sqrt 2(\delb \oplus \delb^*)$ and we have
 $\Dirac^{g_\tau}_{B}=\Dirac^{can}$, for $B$ the spin connection
deduced from the Levi--Civita connection on $M_\tau$  and  the Chern
connection on $L_\tau\simeq
L_{J_\tau}$.  Therefore $(B,\Psi)$ solves the Seiberg--Witten equations
restricted
 to $M_\tau\setminus \overM$  with $\varpi_\tau=0$.
 $(B,\Psi)$ is called the \emph{canonical solution}.
Notice that the Dirac operator, the projection $+$ and
$(B,\Psi)=(B_\tau,\Psi_\tau)$ depend on $\tau$.

A suitable
family of cut-off functions is now  needed. Let $\chi(t)$ be a smooth
 decreasing function  such that
$$
\left\{
\begin{array}{rll}
\chi(t) = 0 & \mbox{if} & t\geq 1 \\
\chi(t) = 1 & \mbox{if} & t\leq 0
\end{array}\right.
$$
Put
\begin{equation}
\label{eq:cutoff}
\fbox{$\chi_\tau=\chi\left(\frac {t -\tau}{N_0} + 1\right )$}
\end{equation}
where $N_0$ is any number with $N_0\geq 1$; $N_0$ will be fixed later
on
to make the derivatives of $\chi_\tau$ as small as required in our
constructions.	We define a cut-off function on $M_\tau$ by the
formula $\chi_\tau(\sigma_\tau)$. By a slight abuse of notation, the
latter function will be denoted $\chi_\tau$ as well.

For a given perturbation of Seiberg--Witten equation $\varpi$ on $M$,
 the perturbation of the equations on $M_\tau$ is defined by
\begin{equation}
\label{varpitau}
 \varpi_\tau = \chi_\tau \varpi.
\end{equation}

\subsubsection{Linear theory}
\label{seclintheo}
The Study of Seiberg--Witten equations requires introducing suitable
 Sobolev spaces rather than using naive smooth objects
defined in \fullref{subsubmonopoleeq}. We recall very quickly the
 results of \cite{KM} in this section.

The configuration $(B,\Psi)$ is a solution of Seiberg--Witten equations
on the almost K\"ahler end of $M_\tau$; hence we study solutions
$(A,\Phi)$ with the same asymptotic behavior.  We introduce the
configuration space
\begin{multline*}
 \cC_l(M_\tau) =\{ (A,\Phi)\in \Conn(W_\tau)\times \Gamma(W_\tau^+)\; / \;
  A-B\in L_l^2(g_\tau) ,\\ \mbox{
  and }\; \Phi - \Psi \in
L^2_{l}(g_\tau,B) \},
\end{multline*}
 the gauge group
$$ \cG_l(M_\tau) =\{ u:M_\tau\to \CC \; /  \; |u|= 1, \mbox{ and }\;
1-u \in
L^2_{l+1}(g_\tau) \},
$$
acting	on $\cC_l$ by
$u\cdot(A,\Phi) = (A - u^{-1}d
u, u\Phi )$, and, for some fixed $\epsilon_0>0$, the perturbation space
$$ \cN (M_\tau)=e^{-\epsilon_0\sigma_\tau}C^r(i\gsu(W_\tau^+)),
$$
 equipped with the norm
$$\|\varpi\|_{\cN_\tau}=\|e^{\epsilon_0\sigma_\tau}\varpi\|_{C^r(g_\tau)}.$$
The $L^2_k(g_\tau)$--norm is the usual $L^2$ norm with $k$ derivatives
on $M_\tau$ defined  using the metric $g_\tau$. In order to define
a
similar norm on the spinor fields, a unitary connection $A$ is
needed. We put
$$
\|\phi\|^2_{L^2_{l}(g_\tau,A)} = \int_{M_\tau } \left ( |\phi |^2
+|\nabla_A\phi|^2 +
  \cdots + |\nabla_A^l \phi|^2 \right ) \vol^{g_\tau},$$
and define $L^2_{l}(g_\tau,A)$ as the completion of the space of
smooth sections for this norm. For two different connections, $A$ and
$A'$ with $A-A'\in L^2_l(g_\tau)$, Sobolev multiplication theorems
show
that the norms $L^2_l(A,g_\tau)$ and $L^2_l(A',g_\tau)$ are
commensurate.

\begin{remark*}
For any pair $\tau,\tau'$, it is easy to construct a diffeomorphism
$f\co M_\tau \to M_\tau'$
 covered by  an isomorphism $F$ between $W_\tau$ and
$W_{\tau'}$ which are a dilations near infinity. Therefore
 $F^*\cC_l(M_{\tau'})= \cC_l(M_{\tau})$ and $f^*\cG_l(M_{\tau'})=
 \cG_l(M_{\tau})$.
In this sense, the spaces $\cC_l$ and $\cG_l$ are
in fact independent of $\tau$. However, the fact that the norms depend on
$\tau$ will become crucial for analyzing the compactness
properties of the family of moduli spaces on $M_\tau$.
\end{remark*}
Of course the choice of $l\ge 2$ is actually perfectly arbitrary
thanks to elliptic regularity. However it will be chosen with $l\geq
4$ so that  we have the inclusion $L^2_l\subset C^1$.
The Sobolev multiplication theorem shows that $\cG_l$ is a Hilbert
  Lie group acting smoothly on the Hilbert affine space $\cC_l$.
Furthermore action of $\cG_l$ is \emph{free}:  if we have
$u\cdot(A,\Phi)=(A,\Phi)$, then $du=0$ hence $u$ must be constant. Now
  $u-1 \in L^2_{l+1}$ therefore $u=1$.

Let $\cZ_l^\tau$ be the space of configurations $(A,\Phi)\in\cC_l$ which
verify the Seiberg--Witten equations \eqref{eqSW2} on $M_\tau$ with
perturbation $\varpi_\tau$. Then
$\cZ_l^\tau$ is
invariant under the gauge group action and we define
$$\cM_l(M_\tau)  =\cZ^\tau_l/\cG_l.
$$
We drop the reference to the index $\tau$ at the moment, for
simplicity of notation. All of what we say in the rest of
\fullref{seclintheo} holds for $M$ and $M_\tau$ or indeed any
AFAK manifold.
The linearized action of the gauge group at an arbitrary configuration
$(A,\Phi)\in \cC_l$
is given by
a differential operator
\begin{align*}
\delta_{1,(A,\Phi)} \co   L^2_{l+1}(i\RR) \longrightarrow  &
L^2_l(i\Lambda^1) \oplus
L^2_{l,B}(W^+) \\
v \longmapsto & (-dv,v\Phi)
\end{align*}
and its formal adjoint is given by
$$ \delta_{1,(A,\Phi)}^* (a,\phi) = -  d^* a + i\Im \langle\Phi,
\phi\rangle ;
$$
notice that with our convention, the Hermitian product
$\ip{\cdot,\cdot}$ is anti-complex linear in the first variable.

A tangent vector $(a,\phi)$ is $L^2$--orthogonal to the orbit through
$(A,\Phi)$ if and only if $\delta_1^*(a,\phi)=0$; the orbit space
$\cC_l/\cG_l$ is a smooth Hilbert manifold, and its
tangent space at $(A,\Phi)$ is identified with $\ker
\delta_{1,(A,\Phi)}^*$.

The linearized
Seiberg--Witten equations at $(A,\Phi)$  are given as well by a
differential
operator:
\begin{align*}
\delta_{2,(A,\Phi)} \co    L^2_{l+1}(i\Lambda^1) \oplus
L^2_{l+1,B}(W^+) &\rightarrow  L^2_l(i\gsu W^+) \oplus
L^2_{l,B}(W^-) \\
 (a,\phi) \mapsto  & (d^+a -\{\Phi \otimes \phi ^* +\phi\otimes
\Phi ^*\}_0, \Dirac_A\phi + a \cdot \Phi )
\end{align*}
\textbf{Remark}\qua
There  is a slight inconsistency in the conventions
of~\cite{KM}. The linear theory  studied  there is exactly the one
presented in the current paper. However, this is not the one of the
equations written
in~\cite{KM} where  $F_A$ is replaced by the
curvature
$F_{\what A}$ of the unitary connection $\what A$ induced by $A$ on the
determinant line
bundle $L$. Since $F_A = 2 F_{\what A}$, the corresponding linearized
equations should be $ \delta_{2,(A,\Phi)}(a,\phi) = (2d^+a -\{\Phi
\otimes \phi ^* +\phi\otimes
\Phi ^*\}_0, \Dirac_A\phi + a \cdot \Phi )$.

At at solution of Seiberg--Witten equations $(A,\Phi)$, the operators
verify
$\delta_{2,(A,\Phi)}\circ \delta_{1,(A,\Phi)}=0$, so we have an
elliptic complex:
\begin{multline*}
 0\rightarrow L^2_{l+2}(i\RR) \stackrel {\delta_{1,(A,\Phi)}}
\longrightarrow   L^2_{l+1}(i\Lambda^1) \oplus
L^2_{l+1,B}(W^+) \\
 \stackrel {\delta_{2,(A,\Phi)}}\longrightarrow
L^2_l(i\gsu W^+) \oplus
L^2_{l,B}(W^-) \longrightarrow 0 
\end{multline*}
$H_0=0$ since the action of the gauge group is free,
$H^1= \ker \delta_2/\Im \delta_1$
is the virtual tangent space to the Seiberg--Witten moduli
space at $(A,\Phi)$, and  $H^2= \coker \delta_2$
is the obstruction space.
Equivalently $H_1$ can be viewed as the kernel of the elliptic
operator
\begin{equation}
\label{eqlinop}
\cD_{(A,\Phi)} = \delta_{1,(A,\Phi)}^* \oplus
\delta_{2,(A,\Phi)} .
\end{equation}

{\bf Facts}
\begin{itemize}
\item The moduli space $\cM_l$ is compact.
\item By elliptic regularity, $\cM_l=\cM_{l+1}=\cM_{l+2}=\cdots=\cM_{r+1}
$;
  this is precisely why the choice of $l$ does not matter. So
  the moduli space is simply referred to by $\cM$.
\item By Sard--Smale theory, we may always  assume
that $H^2=0$ after choosing a suitable
\emph{generic} perturbation $\varpi_\tau$. Then $\cM$ is unobstructed;
it is
a smooth manifold of
dimension equal to the virtual dimension
$$ d= \ip{e(W^+,\Psi),[\overM, Y]},
$$
which is nothing else but the index of $\cD_{(A,\Phi)}$.
\end{itemize}

From now on, $\varpi$ will be	a generic
perturbation of Seiberg--Witten equations on $M$; then we  define
$\varpi_\tau$ by~\eqref{varpitau}
on $M_\tau$. We will show in \fullref{corunobstr} that
for $\tau$ large enough, then the moduli space on $M_\tau$ becomes
unobstructed. We may now state the general version of our excision
theorem.
\begin{mytheo}
\label{theogluing}
Let $\overM$ be a manifold with a contact boundary $(Y,\xi)$ and
an element $(\gs,h)\in \spinc(\overM,\xi)$.
 Let $Z$ and $Z'$ be two AFAK ends
compatible with the contact structure of $Y$. Let $M_\tau$ and
$M'_\tau$ be the AFAK manifolds obtained as the connected sum of
$\overM$
with $Z$ or $Z'$ along $Y$, together with the Riemannian metrics
$g_\tau$ and $g'_\tau$ and the particular Seiberg--Witten equations
with perturbation $\varpi_\tau$ constructed in \fullref{secSWeq}.

Then, for $\tau$ large enough, the moduli spaces $\cM(M_\tau)$ and
$\cM(M'_\tau)$ are generic, and there is a
diffeomorphism
$$ \gG\co \cM_{\varpi_\tau}(M_\tau)\to \cM_{\varpi_\tau}(M'_\tau).
$$
Furthermore there is a canonical identification of the set
of consistent orientations for $M_\tau$ and  $M'_\tau$.   Using
this canonical identification the above diffeomorphism becomes
orientation preserving.
\end{mytheo}

\begin{remark}
\label{remarkglue}
  If we remove the assumption~\eqref{assumption} for the cobordisms
  $Z$ and $Z'$, then the extension maps  $ j \co  \spinc(\overM,\xi)\to
\spinc(M_\tau,\omega_\tau)$ and  $ j' \co  \spinc(\overM,\xi)\to
\spinc(M'_\tau,\omega'_\tau)$ are not
  injective in general.

Then we may still prove a  generalization of
  \fullref{theogluing}: assume that $Z$ is just the
  symplectization of $(Y,\xi)$ and that $Z'$ is an AFAK end as before,
  but without assuming the property~\eqref{assumption}.
To discuss the generalization we need to make the notation more
precise. We denote by $\cM_{\varpi_\tau}(M_\tau,j(\gs,h))$ the moduli
space for
  some choice of $(\gs,h)\in\spinc(\overM,\xi)$. Similarly we
  have the moduli space $\cM_{\varpi_\tau}(M'_\tau,j'(\gs,h))$.
Then, we may form the moduli space $\wtilde
  \cM_{\varpi_\tau}(M_\tau,(\gs',h'))$ for some choice of
  $(\gs',h')\in\spinc(M'_\tau,\omega'_\tau)$ defined by
\begin{equation}
\label{eqbigmoduli}
\wtilde \cM_{\varpi_\tau}(M_\tau,(\gs',h')):= \bigsqcup_{
\begin{array}{c}
(\gs,h)\in\spinc(\overM,\xi)  \\
  j'(\gs,h)=(\gs',h')
\end{array}}   \cM_{\varpi_\tau}(M_\tau,j(\gs,h)).
\end{equation}
Them the conclusion of \fullref{theogluing} is the same if we replace
$\cM_{\varpi_\tau}(M'_\tau)$ by
$\cM_{\varpi_\tau}(M'_\tau,$ $(\gs',h'))$ and $\cM_{\varpi_\tau}(M_\tau)$
by $\wtilde
\cM_{\varpi_\tau}(M_\tau,(\gs',h'))$. In particular, there is an
orientation preserving diffeomorphism
$$ \gG\co  \wtilde \cM_{\varpi_\tau}(M_\tau,(\gs',h')) \to
\cM_{\varpi_\tau}(M'_\tau,(\gs',h')).
$$
\end{remark}

The rest   of \fullref{secexi} is devoted to constructing the map
$\gG$ by a gluing technique and to showing that it is a diffeomorphism.

\subsubsection{Compactness}
\label{subsubcomp}
In this section, we refine the result of compactness for one fixed moduli
space $\cM(M_\tau)$, by showing that a sequence of solutions
$(A_\tau,\Phi_\tau)$ of Seiberg--Witten  equations on $M_\tau$ converge
in some sense as $\tau\to\infty$, up to extraction of a subsequence,
and modulo gauge
transformations, to a solution of Seiberg--Witten equations on $M$.

We review the arguments proving the compactness of one particular
moduli space
in \cite{KM} and explain how to apply them to the family  $M_\tau$.

\begin{lemma}
\label{lemmac0}
There exist constants $\kappa_1,\kappa_2$ such that for every $\tau$
and every
solution of Seiberg--Witten equations $(A,\Phi)$ on $M_\tau$ we have
the estimate
$$\|\Phi\|_{C^0}\leq \kappa_1 + \kappa_2 \|\varpi_\tau
\|_{C^0}.$$
\end{lemma}
\begin{proof}
 By
construction of the moduli space, $\Phi-\Psi\in L^2_3$; now the
Sobolev inclusion $C^0\subset L^2_3$ together with
\fullref{lemmaboundedgeometry} tells us that the pointwise norm
$|\Phi-\Psi|_{C^0} \rightarrow 0$ near infinity on $M_\tau$.

Hence, either $|\Phi|\leq 1$, and we are done, either this is not
true, and
$|\Phi|$ must have a local maximum at a point $x\in M_\tau$. We apply
the
maximum principle at $x$:
\begin{align*}
0\leq & \frac 12 \Delta |\Phi|^2
 = \ip {\nabla^*_A\nabla_A^{\null}\Phi,\Phi} - |\nabla_A^{\null}\Phi|^2 \\
\leq & \ip {\nabla^*_A\nabla_A^{\null}\Phi,\Phi}
 = \ip {D^2_A\Phi,\Phi} - \ip{F_A^+\cdot \Phi,\Phi} - \frac s4
|\Phi^2|
\end{align*}
where the last identity follows from the Lichnerowicz formula.
Using Seiberg--Witten equations, we have
$$ 0\leq  -\frac 12 |\Phi|^4 -	 \Bigl \langle (F_{B}^+-
  \{\Psi\otimes\Psi^*\}_0) \cdot \Phi,\Phi \Bigr \rangle -
  \ip{\varpi_\tau\cdot\Phi,\Phi} - \frac s4 |\Phi|^2.
$$
The lemma follows from the fact that the pointwise norm
of  $s$, $F_B$, $\Psi$ is bounded independently of $\tau$ by
\fullref{lemmaboundedgeometry}.
\end{proof}

As we saw, the AFAK structure induces a  Chern connection $\what\nabla$
and and a spin connection $B$ on the bundle $W_\tau$
 restricted
to $M_\tau\setminus \overM$.  In the rest of this paper, we will
constantly use the following notation: for every spin connection
$A=B+a$ on $L_\tau$, we define the \emph{twisted Chern connection} on
$W_\tau|_{M_\tau\setminus \overM}$ by:
$$\what \nabla _A \Phi:=\what \nabla\Phi +  a\otimes \Phi $$
 In particular $\what \nabla_B:=\what\nabla $ with this
notation. Notice that $\nabla_B \neq \what\nabla_B$ unless the almost
 complex structure is integrable.

\begin{prop}
\label{propexpdecay}
There exists a compact $K\subset M$ large enough and $\delta>0$ such
that
for every integer $k$ there is a constant $c_k>0$ so that, for every
$\tau$ large enough  and
every solution of Seiberg--Witten equations $(A,\Phi)$ on $M_\tau$, we
have the pointwise estimate on $M_\tau\setminus K$
\begin{equation}
\label{eqgaugeinv}
|1-|\beta|^2-|\gamma|^2|,\quad |\gamma|,\quad |\what\nabla_A^{\null}
\Phi|,\quad
|\what\nabla_A^ 2 \Phi| , \cdots \cdots, |\what\nabla_A^k \Phi| \quad\leq
c_k e^{-\delta\sigma_\tau},
\end{equation}
where $\Phi = (\beta,\gamma)\in \Lambda^{0,0}\oplus \Lambda^{0,2}$.
\end{prop}
\begin{remark*}
The quantities~\eqref{eqgaugeinv} controlled by the lemma are gauge
invariant.
\end{remark*}
\begin{proof}
The lemma was proved in \cite[Proposition 3.15]{KM} for $\tau$
fixed. It is
readily checked that it extends as stated for the family $M_\tau$.
We recall what the ingredients of the proof are.

A configuration $(A,\Phi)$ on $M_\tau$ has an \emph{energy} which is a
gauge invariant quantity defined by
$$ E_\tau(A,\Phi) = \int_{M_\tau \setminus K_0}
\left ( (|\beta|^2+|\gamma|^2-1)^2 +|\gamma|^2 +|\what\nabla_A\Phi |^2
\right ) \vol^{g_\tau},
$$
where all the norms, connections are taken with respect to the
structures defined on $M_\tau$ and $K_0$ is a compact in $M$ containing
$\overM$.

\begin{lemma}
\label{lemmaboundednrj}
There exist a compact $K_0\subset M$ large enough,
and  some  constants $\kappa_3$ and $\kappa_4$, such that for every
$\tau$ large enough and  every solution of Seiberg--Witten equations
$(A,\Phi)$ on $M_\tau$, we have
$$ E_\tau(A,\Phi) \leq \kappa_3 + \kappa_4 \|\varpi_\tau\|_{\cN_\tau}^2.
$$
\end{lemma}
\begin{proof}
The proof of this lemma is the same proof than for \cite[Lemma
  3.17]{KM}. The fact that
$\kappa_3$ and $\kappa_4$ do not depend on $\tau$
is insured by \fullref{lemmaboundedgeometry}.
\end{proof}

If we look carefully at the proof, using the notation of
\cite[page 232]{KM}, we read the claim that
$$ \int_{K_3} d a\wedge \omega
$$
can be controlled, for $K_3$ a compact domain large enough,
$\omega$ a closed form extending the symplectic form on the whole
manifold as in \fullref{dfnafak}, and $a=A-B$ decaying
exponentially fast. No explanation of this is given and we
provide one now.

Pick an arbitrarily small $\epsilon>0$.
We have
$$ \left |  \int_{K_3} d a\wedge \omega \right |\leq   \frac \epsilon
2 \int_{K_3} |d
a|^2 +	\frac 1{2\epsilon} \int_{K_3} |\omega|^2 \leq \frac \epsilon 2
\int_{M} |d a|^2
 +  \frac 1{2\epsilon} \int_{K_3} |\omega|^2 .
$$
Using the exponential decay of $a$, we have as in the case of a
compact manifold
$$\int_{M} |d a|^2 = 2 \int_{M} |d^+ a|^2.$$
Therefore
$$ \left |  \int_{K_3} d a\wedge \omega \right |\leq   \epsilon
 \int_{M\setminus K_3} |d^+
a|^2 + \epsilon \int_{K^3} |d^+ a|^2
 +  \frac 1{2\epsilon} \int_{K_3} |\omega|^2.
$$
Using the $C^0$ bound on $d^+a$, we deduce a constant $C>0$,
such that
$$ \left |  \int_{K_3} d a\wedge \omega \right |\leq   \epsilon
 \int_{M\setminus K_3} |d^+
a|^2 +	C,
$$
and the proof of \cite[Lemma 3.17]{KM} is now complete.

\begin{lemma}
\label{lemmacontrolbeta}
For every $\epsilon>0$ there exists a compact set $K_1\subset M$ large
enough, such that for every   $\tau$ large enough, every solution
$(A,\Phi)$
of Seiberg--Witten
equations   on $M_\tau$ verifies
$$ |\beta |\geq  1-\epsilon,
\mbox{ on $M_\tau\setminus K_1$, where $\Phi =(\beta,\gamma)\in
\Lambda^{0,0}\oplus \Lambda^{0,2}$.}
$$
\end{lemma}
\begin{proof}
Assume the lemma is not true. Then, there is a sequence of
$\tau_j\rightarrow \infty$,  a sequence of
      Seiberg--Witten solutions $(A_j,\Phi_j)$ on $M_{\tau_j}$
such that we have
$|\beta_j(x_j)|< 1-\epsilon$  (using the decomposition
$\Phi_j=(\beta_j,\gamma_j)$)  at some  points $x_j\in M_{\tau_j}$
verifying
$\sigma_{\tau_j}(x_j) \rightarrow \infty$.

We restrict $(A_j,\Phi_j)$ to the ball of center $x_j$ and radius
$\sigma_{\tau_j}(x_j)/\kappa$. Thus we have a sequence of Seiberg--Witten
solutions on a sequence of balls of increasing radius. By
\fullref{lemmaboundedgeometry}, the Riemannian metric and symplectic
form on the balls converge to the standard structures on $\RR^4$.

After taking suitable gauge transformations, we can extract a
subsequence  converging smoothly on every compact set  to
a solution $(A,\Phi)$ of Seiberg--Witten equations
on $\RR^4$, with its standard metric, symplectic form, and
perturbation $\varpi=0$. In particular $\Phi=(\beta,\gamma)$, with
$|\beta|\leq 1-\epsilon$ at the origin since $|\beta_j(x_j)|<
1-\epsilon$ by assumption. Using \fullref{lemmaboundednrj}, we
deduce that $(A,\Phi)$ has a bounded energy.  By
Lemma~\cite[3.20]{KM},
$(A,\Phi)$ must be gauge equivalent to the standard solution
$(B,\Psi)$ of Seiberg--Witten equations on the Euclidean space;
this is a contradiction
since $\Psi=(1,0)\in\Lambda^{0,0}\oplus \Lambda^{0,2}$, hence we should
have $|\beta|=1$.
\end{proof}

\fullref{propexpdecay} now follows from \fullref{lemmaboundedgeometry},
\fullref{lemmacontrolbeta} and the  local behavior of Seiberg--Witten
solutions~\cite[Proposition 3.22]{KM}.
\end{proof}

\subsubsection{Gauge with uniform exponential decay}
The next corollary is essential to study the compactness property of
$\cM_{\varpi_\tau}(M_\tau)$.
\begin{cor}
\label{gaugeexpdecay}
There exist constants $c_k,\delta>0$ and a compact set $K\subset M$
such that for every $\tau$ large enough, every solution of
Seiberg--Witten equations  on $M_\tau$ admits a particular gauge
representative $(A,\Phi)$
such that we have the pointwise estimates
\begin{multline}
\label{eqexpdecay}
 |A-B|, \quad |\nabla (A-B)|,\cdots \cdots, |\nabla^k(A-B)|, \\
 |\Phi-\Psi|, \quad
 |\what \nabla _A^{\null} \Phi|, \quad | \what \nabla_A^2 \Phi|,\cdots
 \cdots,
 |\what\nabla_A^{k-1} \Phi| \quad \quad\mbox{ and }\quad
|\what\nabla_A^k \Phi| \leq c_k e^{-\delta \sigma_\tau}
\end{multline}
on $M_\tau\setminus K$.
\end{cor}
\begin{remark*}
We could state a similar Corollary by replacing
$\what\nabla^j_A\Phi$ in  \eqref{eqexpdecay} with $\nabla_A^j(\Phi -
\Psi)$. Indeed, the spin connection $\nabla_B$ tends uniformly
to the
Chern connection $\what \nabla$ on unit balls with center going to
infinity
thanks to
\fullref{lemmaboundedgeometry}.
\end{remark*}
\begin{proof}
 Let
$(A,\Phi)\in \cZ_l^\tau$ be  a solution of Seiberg--Witten equations on
 $M_\tau$. We
 follow closely the proof of \cite[Corollary 3.16]{KM}: we begin to
 increase the
regularity of $(A,\Phi)\in\cC_l$. For that purpose, we can apply a
gauge transformation $u_1\in\cG_l$ after which $(A,\Phi)$ satisfies
the Coulomb condition
$$\delta_{1,{(A,\Phi)}}^* (A-B,\Phi-\Psi) = 0
$$
near infinity. The linear theory for Seiberg--Witten equations shows
that in this gauge we have an additional regularity
$(A-B,\Phi-\Psi)\in L^2_{l+1}$. From this point, we can define the map
$u\co M_\tau \setminus K\rightarrow S^1$ by
$$u = \frac {|\beta|}\beta.
$$
Kato inequality and Sobolev multiplication theorem show that $1-u\in
L^2_{l+1}$. In particular, it is homotopic to $1$ outside a very large
compact set. By the assumption~\eqref{assumption}, we deduce that $u$ is
also homotopic to $1$ on $M_\tau\cap\{1\leq \sigma_\leq \tau\}$. In
particular, $u$
can be extended to $\overM$ thus defining a gauge
transformation in $\cG_l(M_\tau)$.

The gauge transformed solution $(a-u^{-1}du, u\beta,u\gamma) $ has
 exponential decay thanks to \fullref{propexpdecay}, using
the fact that
$u\beta= |\beta|$ and the identity
$$ a - u^{-1}du=
(u\what\nabla_A\beta-d|\beta|)/|\beta|.\proved
$$
\end{proof}

\begin{theo}
\label{theocompact}
Let $(A_\tau,\Phi_\tau)_{\tau\in I}$ be a sequence of solutions of
Seiberg--Witten
equations on $M_\tau$ where $I$ is an unbounded subset of $\RR^+$.
Then, we can extract a subsequence $\tau_j\rightarrow \infty$ and
apply gauge transformations $u_j$, such that the gauge transformed
solutions
$u_j\cdot (A_{\tau_j},\Phi_{\tau_j})$ have uniform exponential
decay, in the sense that they verify
\eqref{eqexpdecay}, and
converge smoothly on every
compact set of $M$ toward a solution of Seiberg--Witten equations
$(A,\Phi)$ on $M$.
\end{theo}
\begin{proof}
As suggested in the Theorem, we apply first
\fullref{gaugeexpdecay}, i.e.\ we apply	a gauge transformations
$u_\tau$ to
$(A_\tau,\Phi_\tau)$ such that
\fullref{gaugeexpdecay} is verified by the gauge transformed
solutions $(A'_\tau,\Phi'_\tau)=u_\tau\cdot (A_\tau,\Phi_\tau)$. Now,
$ (A'_\tau,\Phi'_\tau) $ have  uniform exponential decay outside a
certain compact set $K\subset M$; in particular,
$\|(A,\Phi)-(B,\Psi)\|_{L^2_l(g_\tau,B)}$  is uniformly bounded; a
diagonal argument shows that
we can extract a subsequence $\tau_j\to\infty$ such that
$(A_{\tau_j}',
\Phi_{\tau_j}')$ converge strongly on every compact set of $M\setminus
K$  in the $L^2_{l-1}$ sense to a weak limit $(A',\Phi') \in
L^2_l(g,B)$
solution of Seiberg--Witten equation on $M\setminus K$.
Elliptic regularity shows that, in fact,
$(A'_{\tau_j},\Phi'_{\tau_j})$ converge  smoothly to
$(A',\Phi')$ on every compact of
$M\setminus K$.

Consider the  compact manifold with boundary  $K_2 = M\cap \{\sigma
\leq T_2\}$ and choose $T_2$ large enough, so that $K$
is properly
contained in $K_2$.  There exists a sequence of gauge transformations
$v_j\in L^2_{l+1}(K_2)$ such that after passing to a subsequence,
the transformed solutions $v_j\cdot (A_{\tau_j},\Phi_{\tau_j})$ converge
smoothly
in $K_2$ (see \cite{KM2} for a proof of this statement).

Now, the ratio $w_j=u_jv_j^{-1}$ must converge to a gauge
transformation $w$ on $K_2\setminus K$. Choose $j_0$ sufficiently
large so that $|w_j-w_{j_0}|\leq 1/2$ for all $j\geq j_0$; then $j$
large, we may write
$w_j=\exp(2i\pi\theta_j) w_{j_0}$ where $\theta_j$ is a real valued
function
with $|\theta_j|<1$.
The gauge transformation $u_{j_0}$ on
$M_{\tau_{j_0}}$ is homotopic to $1$ near infinity. Therefore, by
assumption \eqref{assumption}, it must be homotopic to $1$ on
$K_2\setminus K$ and we write $u_{j_0}=\exp(2i\pi \mu)$ for some real
function $\mu$ on $K_2\setminus K$.
 Then, we define gauge transformations on
$M_{\tau_j}$ by
$$ f_j = \left \{
\begin{array}{l}
v_jv_{j_0}^{-1} \exp(2i \pi (1-\chi_{T_2})(\theta_j +\mu)) \mbox{ on
$K_2\setminus K$} \\
u _j \mbox{ on $M_\tau\setminus K_2$}.
\end{array} \right .
$$
now the $f_j\cdot (A_{\tau_j},\Phi_{\tau_j})$ converge smoothly on
every compact
set and
have the uniform exponential decay property as required in the
theorem.
\end{proof}

\begin{remark}
Notice	that the assumption~\ref{assumption} was used for the first
time in the proof of this compactness theorem. Withouth  this
assumption, the gauge transformation $u_{j_0}$ could be  non
homotopic to $1$. Then, a sequence of solutions of Seiberg--Witten
equations on $M_\tau$ for some choice of \spinc--structure
$(\wtilde\gs,\wtilde h)\in \spinc(M_\tau,\omega_\tau)$ still converge up
to gauge
transformation on every compact
to a solution of Seiberg--Witten equations on $M$. However, the
relevant \spinc--structure on $M$ is not uniquely defined and the limit
is now an
element of the enlarged moduli space
$$\wtilde
  \cM_{\varpi}(M,(\wtilde \gs,\wtilde h)) = \bigsqcup_{
\begin{array}{c}
(\gs,h)\in\spinc(\overM,\xi)  \\
  j(\gs,h)=(\wtilde \gs,\wtilde h)
\end{array}}   \cM_{\varpi}(M, j(\gs,h)).$$
\end{remark}

\begin{cor}
\label{boundedderivatives}
For every integer $k$, there exists a constant $c_k>0$, such that for
 every $\tau$ large enough
 and every solution of Seiberg--Witten equations $(A,\Phi)$ on
$M_\tau$, we have
the pointwise estimates
$$ |\nabla_A\Phi|,\quad |\nabla_A^2\Phi|, \cdots\cdots,
|\nabla_A^{k}\Phi|\quad  \leq
c_k.
$$
Moreover, we may assume after applying a suitable gauge transformation
to any other solution that we have
$$ |A-B|, \quad |\nabla(A-B)|, \cdots\cdots, |\nabla^{k}(A-B)|\quad
\leq c_k.
$$
\end{cor}
\begin{proof}
Suppose first assertion is false. Then, we have a sequence
$\tau_j\rightarrow
\infty$, solutions $(A_j,\Phi_j)$ of Seiberg--Witten equations	and
 a sequence of points $x_j\in M_{\tau_j}$ with	$
 |\nabla_{A_j}^{r}\Phi|(x_j)\rightarrow
\infty$. After
applying gauge transformations and extracting a subsequence, we can
assume that $(A_j,\Phi_j)$ converge smoothly on every compact
set. Therefore,
the point $x_j$ must go to infinity on the end of $M_{\tau_j}$ in the
sense that $\sigma_{\tau_j}(x_j)\to \infty$.
But this is a contradiction according to
\fullref{gaugeexpdecay}
and the first part of the  corollary is proved.

Suppose the second part is false for $k=0$. Then, we have a sequence
$\tau_j\to\infty$,
and solutions $(A_j,\Phi_j)$ of
Seiberg--Witten on $M_{\tau_j}$ such that, for every
sequence of gauge transformations $u_j$, we have  $\|u_j\cdot A_j-
B\|_{C^k(g_{\tau_j})}\to \infty$.
This is in contradiction with \fullref{theocompact} and
\fullref{gaugeexpdecay}.
\end{proof}

As an immediate, although important, consequence  of
 \fullref{lemmac0}, \fullref{boundedderivatives}
and the
Sobolev embedding $L^2_1\subset L^4$. , we have the following
corollary.
\begin{cor}
\label{coreqnorms}
For every $k$, there exists a constant $a_k\in(0,1)$, such that for
 every $\tau$ large enough,  every
solution of Seiberg--Witten equations  on
 $M_\tau$ admits a gauge representative $(A,\Phi)$ such that
$$ a_k \|\cdot \|_{L^2_{k}(g_\tau,A)}\leq \|\cdot
\|_{L^2_{k}(g_\tau,B)
} \leq \frac 1{a_k}
\|\cdot \|_{L^2_{k}(g_\tau,A)}.
$$
\end{cor}

Another important consequence of uniform exponential decay is that the
constant involved estimates for  local elliptic regularity can be
chosen uniformly, as in the next proposition.
\begin{prop}
\label{prop:contrd1}
There exist  and a compact $K\subset M$, and for every $k\geq 0$ a
constant $c_k>0$,
such that for every $\tau$ large enough and every solution $(A,\Phi)$
 of Seiberg--Witten equations on
$M_\tau$, we have
$$
\|\cD_{(A,\Phi)}(a,\phi)\|_{L^2_k(g_\tau,A)} +
 \|(a,\phi)\|_{L^2(g_\tau,K)}\geq
c_k \|(a,\phi)\|_{L^2_{k+1}(g_\tau,A)},$$
for all $ (a,\phi)\in \Lambda^+(M_\tau)\times W^+(M_\tau)$,
where $L^2(g_\tau,K)$ is the $L^2$--norm restricted to $K$.
\end{prop}
\begin{proof}
We use the Weitzenbock formula derived in \cite[Proposition
  3.8]{KM}: for every section $(a,\phi)$ on $M_\tau$, we have outside
a compact
set $K\subset M$,
\begin{multline*}
\|\cD_{(A,\Phi)}(a,\phi)\|_{L^2(g_\tau, K')}^2 =  \int_{M_\tau\setminus
K'} \Bigg(
|\nabla a|^2 +|\nabla_A\phi|^2 +|\Phi|^2(|\phi|^2+ |a|^2) \\
+  \Ric(a,a)+ \frac s4 |\phi|^2 + \ip{F_A^+\cdot\phi,\phi}
-2 \ip{a\otimes\phi,\nabla_A\Phi} \Bigg )
\vol^{g_\tau}.
\end{multline*}
Using the uniform exponential decay and the fact that the metrics
$g_\tau$
are uniformly asymptotically flat, we deduce that, provided $K$ is
large enough, we have control
$$\|\cD_{(A,\Phi)}(a,\phi)\|_{L^2(g_\tau,M_\tau\setminus K)} \geq
c_1\|(a,\phi)\|_{L^2(g_\tau, M_\tau\setminus K)}
$$
where $L^2(g_\tau, M_\tau \setminus K)$ means the $L^2$--norm restricted
to the set $M_\tau \setminus K$, and where
$c_1>0$ is a constant independent of $\tau$ and $(A,\Phi)$.
Clearly
\begin{equation}
\label{eqcontrell1}
\|\cD_{(A,\Phi)}(a,\phi)\|_{L^2(g_\tau)}  +
 \|(a,\phi)\|_{L^2(g_\tau, K)} \geq
c_2\|(a,\phi)\|_{L^2(g_\tau)},
\end{equation}
for $c_2=\min(1,c_1)$.

Now, consider  two balls   $B_1 \subset B_2$  in $M_\tau$
with small radii $\alpha/2$ and $\alpha$, centered at a point
$x$. Then, by local elliptic regularity,  there exists a constant
$c_3>0$ such
that
$$  \|\cD_{(A,\Phi)}(a,\phi)\|_{L^2_k({g_\tau},A , B_2) } +
\|(a,\phi)\|_{L^2({g_\tau}, B_2) } \geq
c_3  \|(a,\phi)\|_{L^2_{k+1}({g_\tau},A,B_1)},
$$
where $L^2_k({g_\tau},A,B_j )$ is the $L^2_k({g_\tau}, A)$--norm
restricted to the ball $B_j$.

Now, using the fact that the geometry of the $M_\tau$ is uniformly
bounded in
the sense of \fullref{lemmaboundedgeometry},
\fullref{boundedderivatives} and \fullref{coreqnorms}, we
conclude that $c_3$ can be
chosen in such a way that it does not depend neither on $\tau$,
$(A,\Phi)$ nor on the center of
the balls $x$.

Therefore, we can find a constant $c_4>0$, independent of $\tau$ and
$(A,\Phi)$  such that globally
$$  \|\cD_{(A,\Phi)}(a,\phi)\|_{L^2_k({g_\tau},A ) } +
\|(a,\phi)\|_{L^2({g_\tau}) } \geq
c_4 \|(a,\phi)\|_{L^2_{k+1}({g_\tau},A)};
$$
the latter inequality together with the control \eqref{eqcontrell1}
proves the proposition.
\end{proof}

\subsection{Refined slice theorem}
\label{refslice}
 An important  step in developing
 Seiberg--Witten theory, and in particular,
 in order to give a differentiable structure to the moduli space,
consists in
proving a \emph{slice theorem}.

Let $[A,\Phi]$ be the orbit under the gauge group action of a
 configuration $(A,\Phi)\in
 \cC(M_\tau)$.
Then we consider the equivariant smooth map
\begin{equation}
\label{mapslice}
\begin{array}{rlcl}
\nu\co & U   \times \cG  &\longrightarrow & \cC \\
& ((a,\phi), u)& \longmapsto & u\cdot(A+a,\Phi+\phi),
\end{array}
\end{equation}
where $U$ is a neighborhood of $0$ in $\ker
 \delta_{1,(A,\Phi)}^*$.
The usual  slice theorem says that, if $U$ is an open ball small
 enough centered at $0$, then $\nu$ is a diffeomorphism onto a
 $\cG$--invariant open neighborhood of the orbit $[A,\Phi]$.
Therefore  $U\times\cG$ defines an equivariant local chart on $\cC$.
 Such coordinates $\nu$ endow
$\cC/\cG$ with a structure of smooth Hilbert manifold.

In this section, we show that $U$ can be taken
to be a ball of fixed size (independent of both $\tau$ and the
solution to the Seiberg--Witten equations).  Furthermore $\nu(U\times
\cG)$
contains a ball about the solution whose size is also
independent of both $\tau$ and the solution.
This refinement of the slice theorem will  actually be crucial in the
proof that the gluing  map is an embedding.

Before going further, we need to define a metric on the orbit space
$\cC/\cG$: let $(A,\Phi)$, $(\wtilde A, \wtilde\Phi)\in \cC_l$. Put
$$d_{l,[A,\Phi]}([A,\Phi],[\wtilde A, \wtilde \Phi])= \inf_{u\in \cG_l}
\|(A,\Phi)- u\cdot (\wtilde A,\wtilde \Phi)\|_{L^2_l(g_\tau, A)}.
$$
Then we have the following theorem:
\begin{theo}
\label{theoslice}
For every $\alpha_1>0$ small enough, there exists $\alpha_2
\in(0,\alpha_1]$ such that for every $\tau$ large
enough and
every solution of Seiberg--Witten equations $(A,\Phi)$ on $M_\tau$,
the equivariant map $\nu \co  U_{\alpha_1} \times \cG_2  \to \cC_2$
defined by
\eqref{mapslice}, where
$$ U_{\alpha_1}=\{(a,\phi)\in \ker \delta^*_{1,(A,\Phi)}\subset L^2_2
\quad | \quad
\|(a,\phi)\|_{L^2_2(g_\tau,A)} < \alpha_1\}
$$
is a diffeomorphism onto a $\cG$--invariant open neighborhood of
$(A,\Phi)$.
Furthermore  $\nu(U\times \cG_2)$
contains a ball of radius $\alpha_2$; more
precisely, we have $B_{\alpha_2}\subset \nu(U_{\alpha_1}\times \cG_2)$
where
$$B_{\alpha_2}= \{ [\wtilde A,\wtilde \Phi] \quad | \quad  (\wtilde A,\wtilde
\Phi)\in \cC_2 \mbox { with }
d_{2,[A,\Phi]}([A,\Phi],[\wtilde A,\wtilde \Phi]) <\alpha_2\} .
$$
\end{theo}
\begin{proof}
Let $\alpha_3>0$ and
$$G_{\alpha_3} =\{v\in L^2_3(i\RR)\quad |\quad \|v\|_{L^2_3(g_\tau)}
<\alpha_3 \}.$$
We consider the map $\bar \nu$ rather than $\nu$ defined by
\begin{equation} \label{mapslice2}
\begin{array}{rlcl}
\bar \nu\co & U_{\alpha_3}   \times G_{\alpha_3}  &\longrightarrow &
\cC_2 \\
& ((a,\phi), v)& \longmapsto & e^v\cdot(A+a,\Phi+\phi).
\end{array}
\end{equation}
The map $\bar\nu$ has the decomposition
$$
\bar \nu((a,\phi),v) = (A,\Phi) + (a,\phi) + \delta_{1,(A,\Phi)}(v) +
 Q_{\bar\nu} (v,\phi)
$$
where
$$
Q_{\bar \nu}(\phi,v) =\Bigl (0,(e^v-v-1)\Phi + (e^v-1)\phi \Bigr)$$
is the nonlinear part of $\bar \nu$. Clearly
\begin{equation}
\label{d0nu}
d_0\bar \nu ((\dot a,\dot \phi), \dot v) = (\dot a,\dot \phi) +
\delta_{1,(A,\Phi)}(\dot v).
\end{equation}

We study first the linearized problem for finding an inverse to
$\bar \nu$. We define the Laplacian
$$\Delta_{1,(A,\Phi)}^{\null} v =  \delta_{1,(A,\Phi)}^*
\delta^{\null}_{1,(A,\Phi)} (v) = d^*dv + |\Phi|^2 v.
$$
Then, the spectrum of $\Delta_{1,(A,\Phi)}$ is bounded from below
  according to the next lemma.

\begin{lemma}
\label{lemmacontrolP}
For every integer $k\geq 0$, there exists a constant $b_k>0$ such that
for every $\tau$ large enough
and every solution of Seiberg--Witten equations $(A,\Phi)$ on
$M_\tau$, the
operator $\Delta_{1,(A,\Phi)}$ satisfies
\begin{equation}
\label{eqcontrpv}
 b_k \|v\|_{L^2_{k+2}(g_\tau)} \leq
 \|\Delta_{1,(A,\Phi)}(v)\|_{L^2_k(g_\tau)}\leq
 \frac 1{b_k} \|v\|_{L^2_{k+2}(g_\tau)} .
\end{equation}
\end{lemma}

We finish the proof of \fullref{theoslice}
before proving \fullref{lemmacontrolP}. In particular, we have
\begin{equation}
\label{eqcontrd1}
 b_1 \|v\|_{L^2_{3}(g_\tau)} \leq
 \|\Delta_{1,(A,\Phi)}(v)\|_{L^2_1(g_\tau)}
\end{equation}
so $\Delta_{1,(A,\Phi)} \co  L^2_3\to L^2_1$  is Fredholm and
  injective. On the other hand $\Delta_{1,(A,\Phi)}$ is
  self-adjoint, hence its index is $0$, therefore it must be an
  isomorphism.
It is now readily seen that the differential $d_0\bar \nu$
(see~\eqref{d0nu}) has an inverse of the form
$$ d_0\bar \nu^{-1} (\wtilde a, \wtilde \phi) = \Bigl( (\wtilde a, \wtilde
\phi) -
\delta^{\null}_{1,(A,\Phi)}\Delta^{-1}_{1,(A,\Phi)}\delta_{1,(A,\Phi)}^*
      (\wtilde a, \wtilde \phi) ,
      \Delta^{-1}_{1,(A,\Phi)}\delta_{1,(A,\Phi)}^*
      (\wtilde a, \wtilde \phi) \Bigr ).
$$
We define a map
$$F\co	U_{\alpha_3}\times G_{\alpha_3} \to \{(\wtilde
a,\wtilde \phi)	\in L^2_2,\quad  \delta^*_{1,(A,\Phi)}(\wtilde a,\wtilde
\phi)=0\} \times L^2_3(i\RR)$$
 by
\begin{equation}
\label{funcf}
 F((a,\phi),v) =  d_0\bar \nu^{-1}\Bigl(\bar \nu((a,\phi),v)-
 (A,\Phi)\Bigr ).
\end{equation}
Then, a straightforward computation shows that
$$ F\Bigl ((a,\phi),v \Bigr ) = \Bigl ( (a,\phi) , v \Bigr ) +
  Q(\phi,v)$$
where
$$
Q(\phi,v) =   \Bigl ( -
 \delta^{\null}_{1,(A,\Phi)}\Delta^{-1}_{1,(A,\Phi)}
R(\phi,v),  \Delta^{-1}_{1,(A,\Phi)}R(\phi,v)  \Bigr )
$$
with
$$R(\phi,v)= \delta^*_{1,(A,\Phi)} Q_{\bar\nu}(\phi,v)= i\Im \Bigl (
\ip {\Phi, (e^v-1)\phi} +|\Phi|^2(e^v-v-1)
\Bigr ).
$$
We are going to show that $F$ is a local diffeomorphism about $0$
using the fixed point theorem. In order to apply it, we need to show
 that $Q$ is
locally contracting in the sense of \fullref{lemmacontract}.
Before hand, it is required to start with the following technical
 lemma.

\begin{lemma}
\label{lemmatech}
For every $\kappa >0$ there exists $\alpha >0$ such that for every
$\tau$ large enough and every pair of complex valued functions $w$ and
$\wtilde w$ on $M_\tau$ and every open set $D\subset M_\tau$, we have
\begin{multline}
\label{eqcontrexp}
\|w\|_{L^2_3(g_\tau,D)},\|\wtilde w\|_{L^2_3(g_\tau,D)} \leq \alpha
\Rightarrow \\ \|e^w - e^{\wtilde w} -
(w-\wtilde w)\|_{L^2_k (g_\tau,D)}\leq  \kappa \|(w-\wtilde
w)\|_{L^2_k(g_\tau, D) },
\mbox{for $k=0,1,2,3$},
\end{multline}
where $L^2_k(g_\tau,D)$ is the usual $L^2_k(g_\tau)$--norm taken on the
open set $D$.
\end{lemma}
\begin{proof}
For $\tau$ fixed, the property \eqref{eqcontrexp} is readily
 deduced from the Sobolev multiplication theorems, the inclusion
 $L^2_3\hookrightarrow	 C^0$ and the fact that the function $f(w)
 =e^w - w -1$ is analytic, with a zero of multiplicity $2$ at $w=0$.

The geometry of $M_\tau$ is
uniformly controlled at infinity by \fullref{lemmaboundedgeometry},
hence the Sobolev constants involved in the above inclusion or
multiplication theorems can be chosen independently of~$\tau$. This
shows that $\alpha$ can be chosen independently of $\tau$, and the
lemma is proved.
\end{proof}

Before going further, we need to introduce a suitable complete Euclidean
norm on $L^2_2(i\Lambda^1\oplus W^+) \oplus L^2_3(i\RR)$ defined by
$$ \|((a,\phi),v)\|_{L^2_{2,3}(g_\tau,A)}:=
\max\Bigl (\|(a,\phi)\|_{L^2_{2}(g_\tau,A)} ,
\|v\|_{L^2_{3}(g_\tau)}\Bigr).
$$
We will use the shorthand
$$ \|(\phi,v)\|_{L^2_{2,3}(g_\tau,A)}:=
\|((0,\phi) ,v)\|_{L^2_{2,3}(g_\tau,A)}.
$$

\begin{lemma}
\label{lemmacontract}
For all  $\kappa >0$
there exists
$\alpha >0$, such that for every
 $\tau$ large enough, every solution of Seiberg--Witten equations
	$(A,\Phi)$ on $M_\tau$, we have for  every  pair $(\phi,v)$
and $(\wtilde \phi,\wtilde v)$
\begin{multline}
\label{eqcontr}
 \|(\phi,v)\|_{L^2_{2,3}(g_\tau,A)}, \|(\wtilde \phi,\wtilde v)
 \|_{L^2_{2,3}(g_\tau,A)} \leq
\alpha \Rightarrow   \\
\|Q(\wtilde \phi,\wtilde v)- Q ( \phi, v) \|_{L^2_{2,3}(g_\tau,A)}\leq
\kappa
\|(\wtilde \phi,\wtilde v) - (\phi,v)\|_{L^2_2(g_\tau,A)},
\end{multline}
\end{lemma}
\begin{proof}
Let $\kappa >0$,  and two pairs  $(\phi,v)$
and $(\wtilde \phi,\wtilde v)$.
Then
$$ R(\wtilde \phi,\wtilde v)-R(\phi,v) = i\Im\Bigl \langle\Phi,
\Bigl ( \left (e^{\wtilde v}-e^v-(\wtilde v-v) \right ) \Phi
+(e^{\wtilde v}-e^v)\wtilde \phi + (e^v-1)(\wtilde \phi-\phi \Bigr)
\Bigr \rangle.
$$
Using \fullref{boundedderivatives}, we deduce that for some
constant $c_1>0$ independent of
$\tau$ and $(A,\Phi)$  we have
\begin{multline*}
c_1 \| R(\wtilde \phi,\wtilde v)-R(\phi,v)
\|_{L^2_1(g_\tau)}
\leq  \|e^{\wtilde v}-e^v-(\wtilde v-v)\|_{L^2_1(g_\tau)} \\ +
\|(e^{\wtilde v}-e^v)\wtilde \phi \|_{L^2_1(g_\tau, A)} +
\|(e^{\wtilde v}-1)(\wtilde \phi - \phi) \|_{L^2_1(g_\tau,A )}  .
\end{multline*}
The Sobolev multiplication theorems show that for some constant
$c_2>0$:
$$ \|(e^{\wtilde v}-e^v)\wtilde \phi \|_{L^2_1(g_\tau,A)}\leq c_2
\|e^{\wtilde v}-e^v \|_{L^2_2(g_\tau)}\|\wtilde \phi \|_{L^2_2(g_\tau,A)};
$$
$$ \|(e^{\wtilde v}-1)(\wtilde \phi-\phi) \|_{L^2_1(g_\tau,A)}\leq c_2
\|e^{\wtilde v}-1 \|_{L^2_2(g_\tau)}\|\wtilde \phi-\phi
\|_{L^2_2(g_\tau,A)};\leqno{\hbox{and}}
$$
 thanks to \fullref{lemmaboundedgeometry} and
 \fullref{boundedderivatives}, the constant
$c_2$ can be chosen independently of $\tau$.

Eventually, using \fullref{lemmatech}, we see that for	$\alpha>0$
small enough, we have
\begin{multline}
\label{eqcontr2}
 \|(\phi,v)\|_{L^2_{2,3}(g_\tau,A)}, \|(\wtilde \phi,\wtilde v)
 \|_{L^2_{2,3}(g_\tau,A)} \leq
\alpha \Rightarrow   \\
\|R(\wtilde \phi,\wtilde v)- R ( \phi, v)
 \|_{L^2_{1}(g_\tau)}\leq
 \kappa_0
\|(\wtilde \phi,\wtilde v) - (\phi,v)\|_{L^2_2(g_\tau,A)}.
\end{multline}
The estimate~\eqref{eqcontrd1} implies that
$$
\|\Delta_{1,(A,\Phi)}^{-1}\Bigl ( R(\wtilde \phi,\wtilde v)- R ( \phi,
v)\Bigr ) \|_{L^2_{3}(g_\tau)}\leq  \frac {\kappa _0}{b_1}
\|(\wtilde \phi,\wtilde v) - (\phi,v)\|_{L^2_2(g_\tau,A)}.
$$
By \fullref{boundedderivatives}, there is a constant $c_3>0$
independent of $\tau$ and of the solution $(A,\Phi)$ of Seiberg--Witten
equations such that for every $(a,\phi)$
$$\| \delta_{1,(A,\Phi)}(a,\phi)\|_{L^2_2(g_\tau,A)}\leq
c_3\|(a,\phi)\|_{L^2_3(g_\tau,A)}.
$$
In particular
$$
\| \delta_{1,(A,\Phi)}^{\null}\Delta_{1,(A,\Phi)}^{-1}\Bigl ( R(\wtilde
\phi,\wtilde v)- R ( \phi,
v)\Bigr ) \|_{L^2_{2}(g_\tau,A)}\leq  \frac {c_3 \kappa _0}{b_1}
\|(\wtilde \phi,\wtilde v) - (\phi,v)\|_{L^2_2(g_\tau,A)}.
$$
Now, if we take $\kappa_0= \min (\frac {b_1\kappa}{2c_3}, \frac
{b_1\kappa}{2})$ from the beginning,  we have
\begin{multline*}
\|Q(\wtilde \phi,\wtilde v) - Q(\phi,v)\|_{L^2_{2,3}(g_\tau,A)}= \|
\delta_{1,(A,\Phi)}^{\null}\Delta_{1,(A,\Phi)}^{-1}\Bigl ( R(\wtilde
\phi,\wtilde v)- R ( \phi,
v)\Bigr ) \|_{L^2_{2}(g_\tau,A)}\\
 + \| \Delta_{1,(A,\Phi)}^{-1}\Bigl ( R(\wtilde \phi,\wtilde v)- R ( \phi,
v)\Bigr ) \|_{L^2_{3}(g_\tau)}\leq  \kappa
\|(\wtilde \phi,\wtilde v) - (\phi,v)\|_{L^2_2(g_\tau,A)},
\end{multline*}
and the lemma holds.
\end{proof}

Recall that there is an effective  version of the contraction
mapping theorem.
\begin{prop}
\label{propcontractmapping}
Let $S\co \EE\to \EE$ be a smooth function on a Banach space
 $(\EE,\|\cdot\|)$ such
 that $S(0)=0$; assume
 that there exist
constants $\alpha >0$ and $\kappa\in(0,1/2)$ such that
$$ \text{for all } x, y\in\EE,\quad \|x\|,\|y\|\leq \alpha \Rightarrow
\|S(y)-S(x)\|\leq  \kappa \| y-x\|.
$$
Then, for every $y$ in the open ball $B_{\alpha/2}(0)$,  the equation
$y=x+S(x)$ has a unique
solution $x(y)\in B_{\alpha}(0)$; moreover, we have $\|x(y)-y\|\leq \frac
\kappa{1-\kappa} \|y\|$, and the smooth function
$$
\begin{array}{rccl}
F\co & B_{\alpha}(0)   &\longrightarrow & B_{3\alpha/2}(0) \\
& x& \longmapsto & x+S(x),
\end{array}
$$
restricted to $F\co F^{-1}(B_{\alpha/2}(0))\to B_{\alpha/2}(0)$ is a
diffeomorphism. In addition $B_{\alpha/3}(0)\subset
F^{-1}(B_{\alpha/2}(0))$ and $F(B_{\alpha/3}(0))$ is an open
neighborhood of $0$ containing the ball $B_{\alpha/6}(0)$.
\end{prop}
\begin{proof}
    Under these assumptions $x\mapsto y-S(x)$ maps the ball of radius
    $\alpha$
    in $\EE$ to itself and is a contraction mapping there.  These
    estimates for the fixed points follow immediately.
 The last claim is deduced easily  by replacing $\alpha$ with $\alpha/3$
in  the first part of the proposition.
\end{proof}

An important ingredient in the refined slice theorem is to study
the size of the open sets on $\cG$ provided by the  exponential map.
More specifically, we have the following corollary.
\begin{cor}
\label{corexp}
For every $\alpha>0$ small enough and every $\tau$ large enough and
any open set $D\subset M_\tau$, the
exponential defines a
smooth map
$$ \exp \co  B_{\alpha}(0) \to B_{3\alpha/2}(1)
$$
where $B_r(u_0)$ is the space of functions $ u \in L^2_3(\CC)$ on $D$
such that, $\|u-u_0\|_{L^2_3(g_\tau,D)} <
r$. Moreover, $\exp$ is a diffeomorphism onto an open neighborhood
of $1$ which contains $B_{\alpha/2}(1)$.
\end{cor}
\begin{proof}
This is a direct consequence of \fullref{lemmatech} and
\fullref{propcontractmapping}.
\end{proof}

As a second application, we have the following result for $F$.
\begin{cor}
\label{corF}
For every $\alpha_3>0$ small enough, every $\tau$ large enough and
every solution $(A,\Phi)$ of Seiberg--Witten equations on $M_\tau$, the
 map
$$F\co	U_{\alpha_3}\times G_{\alpha_3} \to \{(\wtilde
a,\wtilde \phi)	\in L^2_2,\quad  \delta^*_{1,(A,\Phi)}(\wtilde a,\wtilde
\phi)=0\} \times L^2_3(i\RR)$$
 defined
at~\eqref{funcf} is diffeomorphism onto an open neighborhood of
$0$. Moreover,
we have
$$U_{\alpha_3/2}\times G_{\alpha_3/2}\subset F(U_{\alpha_3}\times
G_{\alpha_3})\subset  U_{3\alpha_3/2}\times G_{3\alpha_3/2}.$$
\end{cor}
 \begin{proof}
Put $x=((a,\phi),v)$,	$y=((\wtilde a,\wtilde \phi),\wtilde v)$ and
$S(x) = Q(\phi,v)$. Then $F(x) = x + S(x)$ and we use  the
 $L^2_{2,3}(g_\tau, A)$--norm. Thanks to \fullref{lemmacontract},
for every $\alpha>0$ small enough, every $\tau$ large enough and every
solution $(A,\Phi)$
of Seiberg--Witten equations on $M_\tau$, the  assumption of
\fullref{propcontractmapping} holds.
Then $\alpha_3 = \alpha/3$ satisfies the conclusions of the corollary.
\end{proof}

We state a serie of preparation lemmas for proving \fullref{theoslice}.
\begin{lemma}
\label{dnu0est}
There exists a constant $c_4>0$ such that for every  $\tau$ large
 enough
and every
  solution $(A,\Phi)$ of Seiberg--Witten equations on $M_\tau$, we have
 for all $((\dot a, \dot\phi), \dot v)$, with $(\dot a,\dot
\phi)\in \ker\delta_{1,(A,\Phi)}^*$ the estimate
$$ \|((\dot a,\dot \phi),\dot v)\|_{L^2_{2,3}(g_\tau,A)}\leq c_4
\|d_0\bar \nu((\dot a, \dot \phi), \dot v)\|_{L^2_{2}(g_\tau,A)}.
$$
\end{lemma}
\begin{proof}
Thanks to Lemmas~\ref{boundedderivatives} and
\ref{lemmaboundedgeometry}, there is a constants $c_5, c_6>0$
independent of $\tau$ and of the solution
$(A,\Phi)$  of Seiberg--Witten equations such that
$$ \|\delta_{1,(A,\Phi)}^{*}d_0\bar \nu((\dot a,\dot\phi),\dot
v)\|_{L^2_{1}(g_\tau)}\leq c_5
\|d_0\bar \nu((\dot a,\dot \phi),\dot  v)\|_{L^2_{2}(g_\tau,A)}
$$
and
$$
c_6 \|\delta_{1,(A,\Phi)}(\dot v)\|_{L^2_{2}(g_\tau,A)}\leq
\|\dot	v\|_{L^2_{3}(g_\tau)}.
$$
On the other hand, using the fact that $d_0\bar \nu((\dot a,\dot
\phi),\dot
v) = (\dot a,\dot \phi)+ \delta_{1,(A,\Phi)}(\dot v)$, we have
$$\delta_{1,(A,\Phi)}^{*}d_0\bar \nu((\dot a,\dot \phi),\dot v) =
\Delta_{1,(A,\Phi)} \dot v.
$$
Using the estimate~\eqref{eqcontrd1}, we obtain
$$ \|\dot v\|_{L^2_3(g_\tau)}\leq \frac {c_5}{b_1} \|d_0\bar \nu((\dot
a,\dot
\phi),\dot  v)\|_{L^2_{2}(g_\tau,A)},
$$
hence
$$  \|\delta_{1,(A,\Phi)}(\dot v)\|_{L^2_{2}(g_\tau,A)} \leq
\frac{c_5}{b_1c_6}\|d_0\bar \nu((\dot a,\dot \phi),\dot
v)\|_{L^2_{2}(g_\tau,A)}.
$$
Therefore,
\begin{multline*}
(1 +  \frac{c_5}{b_1c_6} + \frac{c_5}{b_1})  \|d_0\bar \nu((\dot
  a,\dot \phi) ,\dot
v)\|_{L^2_{2}(g_\tau,A)}
  \geq \\
 \Bigl (\|(\dot a,\dot
\phi)\|_{L^2_2(g_\tau,A)} - \|\delta_{1,(A,\Phi)} (\dot
v)\|_{L^2_2(g_\tau,A )} \Bigr ) +  \|\delta_{1,(A,\Phi)}(\dot
v)\|_{L^2_2(g_\tau,A)} +  \|\dot v\|_{L^2_3(g_\tau)} =	\\
 \|(\dot a,\dot
\phi)\|_{L^2_2(g_\tau,A)} +\|\dot
v\|_{L^2_3(g_\tau)} ,
\end{multline*}
and the lemma is proved.
\end{proof}

\begin{lemma}
\label{lemmainjprep}
For every $\alpha_3>0$, there exists a constant $\alpha'_1>0$ such that
for every $\tau$ large
enough and every solution of
Seiberg--Witten equations $(A,\Phi)$ on $M_\tau$, we have for all
$(a,\phi)$, $(\wtilde a,\wtilde \phi)\in \ker \delta_{1,(A,\Phi)}^*$ and
$u\in \cG_2$
$$
\left . \begin{array}{l}
 \|(a,\phi)\|_{L^2_2(g_\tau)}, \|(a,\phi)\|_{L^2_2(g_\tau)}\leq
\alpha'_1 \mbox { and }\\
  u\cdot (A+a,\Phi+\phi) =  (A+\wtilde a,\Phi+
\wtilde \phi)
\end{array}
\right \} \Rightarrow \|u\|_{L^2_3} < \alpha_3/2
$$
\end{lemma}
\begin{proof}
By assumption  $\|\wtilde a-a\|_{L^2_2(g_\tau)} \leq
2\alpha'_1$ and $u^{-1}du =\wtilde a - a$, hence
\begin{equation}\label{eqestu4}
\left \|\frac {du}u \right \|_{L^2_{2}(g_\tau)}= \| {du}
\|_{L^2_{2}(g_\tau)}\leq 2\alpha'_1.
\end{equation}
First, we derive an $L^2_2$ bound on $du$.
\begin{align}
\nabla(\wtilde a- a) = \nabla \left (\frac{du}u \right ) & = - \frac
{du\otimes du}{u^2} +
\frac {\nabla^2u}u \\
\label{eqestu}
&= - (\wtilde a -a)\otimes (\wtilde a -a) + \frac {\nabla^2u}u
\end{align}
The Sobolev multiplication  theorem $L^2_1 \otimes L^2_1 \subset L^2$
says that for some
constant $c>0$ we have
$$ \| (\wtilde a -a)\otimes (\wtilde a -a)\|_{L^2(g_\tau)} \leq c_7 \|
\wtilde a -a\|^2_{L^2_1(g_\tau)} \leq 4c_7\alpha^2_1.
$$
Thanks to \fullref{lemmaboundedgeometry}, the constant $c_7$ can be
chosen independently of $\tau$.
We deduce a bound from~\eqref{eqestu}
$$ \|\nabla^2u\|_{L ^2(g\tau)} \leq (4c_7\alpha'_1+2)\alpha'_1.
$$
Then, we take care of derivatives of order $2$. We have
\begin{equation}\label{eqestu2}
\nabla^2(\wtilde a-a)= \nabla^2\left ( \frac{du}u \right )
= - \nabla\Bigl( (\wtilde a -a)\otimes (\wtilde a -a) \Bigr)
  -2  \frac {du \otimes \nabla^2u}{u^2} + \frac
       {\nabla^3u}u
\end{equation}
Since
$$ \nabla\Bigl( (\wtilde a -a)\otimes (\wtilde a -a) \Bigr) = \nabla
(\wtilde a
-a)\otimes (\wtilde a -a) + (\wtilde a -a)\otimes \nabla (\wtilde a -a),
$$
we can use again the Sobolev embedding theorem to deduce the estimate
\begin{equation}\label{eqestu3}
 \left \| \nabla\Bigl( (\wtilde a -a)\otimes (\wtilde a -a) \Bigr)
\right \|_{L^2(g_\tau)} \leq 8c_7 (\alpha'_1)^2.
\end{equation}
The next  term of~\eqref{eqestu2} can be written
$$-2  \frac {du \otimes \nabla^2u}{u^2} = -2 \frac {du}u\otimes
\nabla\left ( \frac {du}u\right ) -  \frac 2u \frac {du\otimes du}{u^2}.
$$
Again the Sobolev multiplication theorem gives control
\begin{align*}
\left \| -2  \frac {du \otimes \nabla^2u}{u^2}\right\|_{L^2(g_\tau)}
&\leq 2c_7 \left \|\frac {du}u \right \|_{L^2_1(g_\tau)}  \left \|\nabla
\left(
\frac {du}u\right) \right \|_{L^2_1(g_\tau)} + 2c_7  \left\|\frac {du}u
\right\|_{L^2_1(g_\tau)}^2 \\
&\leq 16c_7(\alpha'_1)^2.
\end{align*}
Returning to the identity~\eqref{eqestu2}, we derive the estimate
$$  \left \| \nabla^3 u \right \|_{L^2(g_\tau)}
\leq 2\alpha'_1 +8c_7(\alpha'_1)^2 +  16c_7(\alpha')_1^2 = 2
\alpha'_1(1+12c_7 \alpha'_1).
$$
In conclusion we have proved the estimate
\begin{equation}
\label{eqestu6}
 \| du \|_{L^2_2(g_\tau)} \leq	2\alpha'_1(3+14c_7 \alpha'_1).
\end{equation}
We show now that we have an $L^2$--estimate on $u-1$ outside a compact
set.
By \fullref{propexpdecay}, there is a compact set
$K\subset M$ such that
for every $\tau$ large enough and every solution of Seiberg--Witten
equations $(A,\Phi)$ on
$M_\tau$, we have $|\Phi|\geq \frac 12 $. We write $(1-u)\Phi =
(u-1)\phi + (\phi -\wtilde \phi)$ and
 take  the
$L^2(g_\tau)$--norm on the noncompact set $M_\tau\setminus K$. We have
 the estimate
\begin{equation}
\label{eqinterm1}
 \frac 12 \|u-1\|_{L^2(g_\tau,M_\tau\setminus K)} \leq \|
(u-1)\phi \|_{L^2(g_\tau,M_\tau\setminus K) } + \|
\phi - \wtilde \phi \|_{L^2(g_\tau,M_\tau\setminus K)}.
\end{equation}
By the Sobolev multiplication theorem $L^2_1\otimes L^2_1\hookrightarrow
L^2$, we deduce that for some constant $c_8>0$ (independent of $\tau$
and the choice of a solution $(A,\Phi)$ of Seiberg--Witten equations on
$M_\tau$ by
\fullref{lemmaboundedgeometry} and \fullref{coreqnorms})
\begin{align*}
\|(1-u)\wtilde \phi \|_{L^2(g_\tau,M_\tau\setminus K) } &
 \leq c_8 \|
1-u \|_{L^2_1(g_\tau,M_\tau\setminus K) } \|
\wtilde \phi \|_{L^2_1(g_\tau,A,M_\tau\setminus K) }\\
&\leq  c_8 \alpha'_1 \|
1-u \|_{L^2_1(g_\tau,M_\tau\setminus K) }.
\end{align*}
So if we chose $\alpha'_1\leq 1/4c_8$, we deduce from~\eqref{eqinterm1}
$$ \frac 14 \|u-1\|_{L^2(g_\tau,M_\tau\setminus K)} \leq
\frac 14 \|du\|_{L^2(g_\tau,M_\tau\setminus K)} +
 \|
\wtilde \phi-\phi \|_{L^2(g_\tau,A,M_\tau\setminus K)} .
$$
Using the estimate~\eqref{eqestu4}, we have
\begin{equation}\label{eqestu5}
 \|u-1\|_{L^2(g_\tau,M_\tau\setminus K)} \leq 10\alpha'_1.
\end{equation}
In order to exploit the estimates~\eqref{eqestu6} and~\eqref{eqestu5},
the following lemma is needed.
\begin{lemma}
\label{lemmapoincare}
Let $K\subset M$ be a compact set. There exists a constant $c>0$ such
that for every $\tau$ large enough and every complex valued function
$w$ on $M_\tau$ we have
$$ c \|w\|_{L^2_1(g_\tau)}\leq \|dw\|_{L^2(g_\tau)} +
\|w\|_{L^2(g_\tau,M_\tau \setminus K)}.
$$
\end{lemma}
\begin{proof}
Suppose this is not true. Then we have a compact set $K$ and a
sequence $\tau_j\to\infty$, together with a sequence of functions
$w_j$ on $M_{\tau_j}$ such that
$$ \|w_j\|_{L^2_1(g_{\tau_j})} = 1,\quad \|dw_j\|_{L^2(g_{\tau_j})} \to
0,\quad  \mbox{and} \quad  \|w_j\|_{L^2(g_{\tau_j}, M_{\tau_j}\setminus
K)} \to 0.
$$
Hence, up to extraction of a subsequence, $w_j$ converge in the strong
$L^2$ sense on every compact toward a weak limit $w\in
L^2_2(g,M)$. Moreover $w$ satisfies $dw=0$ on $M$ and $w=0$ on $M\setminus
K$.
Therefore  $w\equiv 0$.

Thus, we have
$$
 \|dw_j\|_{L^2(g_{\tau_j})}
+  \|w_j\|_{L^2(g_{\tau_j}, M_{\tau_j}\setminus K)} +
 \|w_j\|_{L^2(g_{\tau_j}, K)} =  \|w_j\|_{L^2_1(g_{\tau_j})} .
$$
The first two terms converge to $0$ by assumption. The third one
converge to $0$ since $w_j\to 0$ in the $L^2$ sense on every
compact. This is a contradiction since RHS equals $1$ and the lemma
is proved.
\end{proof}

We return now to the proof of \fullref{lemmainjprep}. Thanks to
\fullref{lemmapoincare}, there is a constant $c_9>0$ such that for
every  $\tau$ large enough
$$ c_9 \|u\|_{L^2_1(g_\tau)}\leq \|du\|_{L^2(g_\tau)} +
\|u\|_{L^2(g_\tau,M_\tau \setminus K)}.
$$
Using the estimates~\eqref{eqestu4} and~\eqref{eqestu5}, we deduce
that
$$  \|u-1 \|_{L^2_1(g_\tau)}\leq 12\alpha'_1/c_9.
$$
Together with~\eqref{eqestu6}, it gives
\begin{equation}
\label{eqinterm2}
\|u-1\|_{L^2_{3}(g_\tau)}\leq  12\alpha'_1/c_9 +
2\alpha'_1(3+14c_7\alpha_1').
\end{equation}
For $\alpha'_1$ small enough, we have
$ 12\alpha'_1/c_9 +
2\alpha'_1(3+14c_7\alpha_1') <\alpha_3/2$
and the lemma is proved.
\end{proof}

We can complete now the proof of \fullref{theoslice}.
Let $\alpha_3>0$ be a constant as in \fullref{corF}. Then any
constant $\alpha_1\in (0,\alpha_3]$ also satisfies the Corollary. We
  will in addition assume that $\alpha_1\leq \alpha'_1$, where $\alpha'_1$
  satisfies \fullref{lemmainjprep}.

In particular, we have $U_{\alpha_1/2}\times G_{\alpha_1/2} \subset
 F(U_{\alpha_1}\times G_{\alpha_1})$. Let $c_4$ be the constant
from \fullref{dnu0est}.
Then for every $\tau$ large enough and every solution of
 Seiberg--Witten equations $(A,\Phi)$ on $M_\tau$, we have
$$ \Bigl \{(a,\phi)\in L^2_2,\quad \|(a,\phi)\|_{L^2_2(g_\tau)}<
c_4\alpha_1/2
 \Bigr \}
\subset d_0\bar \nu (U_{\alpha_1/2}\times G_{\alpha_1/2}).
$$
Since $\bar \nu  = d_0\bar\nu \circ F + (A,\Phi)$, we have
$$ \Bigl \{(\wtilde A,\wtilde \Phi)\in \cC_2,\quad \|(\wtilde
A,\wtilde\Phi)-(A,\Phi)\|_{L^2_2(g_\tau,A)} <
c_4\alpha_1/2 \Bigr \} \subset \bar \nu ( U_{\alpha_1}\times G_{\alpha_1}
).
$$
In conclusion, $\bar \nu \co  U_{\alpha_3}\times G_{\alpha_3} \to \cC_2$
is a diffeomorphism onto an open neighborhood of $(A,\Phi)$ and $\bar \nu(
U_{\alpha_1}\times G_{\alpha_1}) $ contains
a ball of radius $\alpha_2:= c_4\alpha_1/2$ as above.

For a choice of $\alpha_3>0$ small enough, \fullref{corexp} says
that the exponential map
$$ \exp \co  G_{\alpha_3} \to \cG_2,
$$
is a diffeomorphism onto an open neighborhood
$G'_{\alpha_3}=\exp(G_{\alpha_3})$ of $1$ containing a ball
centered at $1$ with $L^2_3(g_\tau)$--radius $\alpha_3/2$.

By equivariance, it follows that $\nu \co  U_{\alpha_3}\times \cG_2 \to
  \cC_2$ is
a local diffeomorphism onto an open neighborhood of the orbit of
  $(A,\Phi)$. Moreover, the ball
  $B_{\alpha_2}$ defined in \fullref{theoslice}, is contained in
  $\nu(U_{\alpha_1}\times\cG_2)$.

We just need to show that $\nu\co U_{\alpha_1}\times\cG_2 \to \cC_2$ is
injective  in order to prove that $\alpha_1$ and
$\alpha_2$ satisfy the theorem.
Suppose that $\nu ((a,\phi),u) = \nu
((\wtilde a,\wtilde \phi),\wtilde u)$, for some $(a,\phi), (\wtilde
a,\wtilde \phi)\in U_{\alpha_1}$ and $u, \wtilde u \in \cG_2$. By
equivariance, we have $\nu ((a,\phi),u\wtilde u^{-1} ) = \nu
((\wtilde a,\wtilde \phi),1)$. Since $\alpha_1\leq\alpha'_1$,
the assumption of \fullref{lemmainjprep} are verified, therefore
we must have  $\|u\wtilde u^{-1}\|_{L^2_3(g_\tau)} <\alpha_3/2$, hence
$u\wtilde u^{-1} \in G'_{\alpha_3}$. By	injectivity of
$\bar\nu \co  U_{\alpha_3} \times G'_{\alpha_3}\to \cC_2$  we conclude
that
 $u\wtilde u^{-1} =1$, therefore $u=\wtilde u$ and
$(a,\phi)= (\wtilde a,\wtilde \phi)$.
\end{proof}

We now	return to the proof of \fullref{lemmacontrolP}.
\begin{proof}
According to \fullref{propexpdecay}, there exists a compact
$K\subset M $ such that $|\Phi|^2\geq 1/2$ outside $K$ for every
 solution of Seiberg--Witten equations $(A,\Phi)$ on $M_\tau$ and for
every $\tau$ large enough);
let $\chi_K$ be the characteristic function of $K$.
Then, for every smooth function $v$, we have
\begin{multline}
%\label{eqcontr2}
\int \ip{\Delta_{1,(A,\Phi)}(v),v} \vol^{g_\tau} = \int \left (|dv|^2 +
|\Phi|^2 |v|^2
  \right )\vol^{g_\tau}
 \\
 \geq \|dv\|_{L^2(g_\tau)}^2 +	 \frac 12 \|(1-\chi_K)
  v\|_{L^2(g_\tau)}^2  .
\end{multline}
Therefore, according to \fullref{lemmapoincare}, there exists a
constant $c>0$ (depending only on $K$)
such that for every $\tau$ large enough, we have
\begin{equation}
\label{eqfirstcontrol}
 \|\Delta_{1,(A,\Phi)}(v)\|_{L^2(g_\tau)}\geq c \|v\|_{L^2_1(g_\tau)}.
\end{equation}
The next step is to obtain a control on higher derivatives.
For this purpose, we use the fact that the  operator $\Delta_1$
is elliptic;
 we consider two balls	 $B_1 \subset B_2$  in $M_\tau$
with small radii $\alpha/2$ and $\alpha$, centered at a point
$x$. Then there exists a constant $c'_k>0$ such
that
$$  \|\Delta_{1,(A,\Phi)}(v)\|_{L^2_k({g_\tau} , B_2) } +
\|v\|_{L^2_1({g_\tau}, B_2) } \geq
c'_k  \|v\|_{L^2_{k+2}({g_\tau},B_1)},
$$
where $L^2_k({g_\tau},B_j )$ is the $L^2_k({g_\tau}, A)$--norm
restricted to the ball $B_j$.

Now, using the fact that the geometry of the $M_\tau$ is uniformly
bounded in
the sense of \fullref{lemmaboundedgeometry} and the bounds on
$(A,\Phi)$ obtained from  \fullref{lemmac0} and
\fullref{boundedderivatives}, we deduce that the
coefficients of $\Delta_{1,(A,\Phi)}$ are bounded independently of
$(A,\Phi)$ and $\tau$. We
conclude that $c'_k$ can be
chosen in such a way that it does not depend neither on the center of
the balls $x$, $\tau$
nor on $(A,\Phi)$.
Therefore, we can find	constants $c''_k>0$, independent of $\tau$ and
$(A,\Phi)$  such that globally
$$  \|\Delta_{1,(A,\Phi)}(v)\|_{L^2_k({g_\tau}) } + \|v\|_{L^2_1({g_\tau})
} \geq
c''_k  \|v\|_{L^2_{k+2}({g_\tau}) }.
$$
Hence
$$  (1+1/c)\|\Delta_{1,(A,\Phi)}(v)\|_{L^2_1({g_\tau}) }\geq   c''_k
\|v\|_{L^2_3({g_\tau}) },
$$
  where $c$ is the constant of
 the control \eqref{eqfirstcontrol}. Then  $b_k= \frac {cc''_k} {c+1}$
   satisfies
$$ b_k \|v\|_{L^2_{k+2}(g_\tau)} \leq
\|\Delta_{1,(A,\Phi)}(v)\|_{L^2_k(g_\tau)}.$$
Finally, we can always take a smaller value for $b_k>0$ such that the
second inequality of the lemma is verified. This is a trivial
consequence of \fullref{boundedderivatives}.
\end{proof}

\subsection{Approximate solutions}
\label{subappsol}
This section deals with the first step for constructing the
gluing map $\gG$ of \fullref{theogluing}.
Recall that the family of AFAK manifolds $M_\tau$ was constructed in
\fullref{sec:paste} by
adding an AFAK end $Z$ to a manifold with  $\overM$ with a
contact boundary $(Y,\xi)$.
Let $Z'$ be another AFAK end compatible with the contact boundary
$(Y,\eta)$ and, similarly we have,
$(M'_\tau,g'_\tau,\omega_\tau',J'_\tau,\sigma'_\tau)$ the family of
AFAK manifolds constructed
out of $\overM$ and $Z'$.

Moreover, an element $(\gs,h)\in \spinc(\overM,\xi)$ gives rise to
a \spinc--structure with an identification with $\gs_{\omega'_\tau}$
outside $\overM$, that is to say an element $j'(\gs,h)\in
\spinc(M'_\tau,\omega'_\tau)$. Let $W'_\tau$ be the spinor bundle
of $j'(\gs,h)$ on $M_\tau'$ constructed similarly to $W_\tau\to M_\tau$
(see \fullref{secextspinc}).

The compact domains $ \{\sigma \leq
\tau\} \subset M$, $ \{\sigma_\tau \leq
\tau\} \subset M_\tau$ and $\{\sigma_\tau' \leq \tau\} \subset
M_\tau'$
are equal by construction and on these sets, all the structures match
(Riemannian metrics, functions
$\sigma$, almost K\"ahler structures,	\spinc--structures, spinor
bundles, canonical solutions $(B,\Psi)$ and $(B',\Psi')$).

We will now explain  how to construct an approximate solution of
Seiberg--Witten
equations on $M'_\tau$ from a solution of Seiberg--Witten equations
on $M_\tau$.

\subsubsection{Construction of the spinor bundle}
On the end $\{1 \leq \sigma_\tau \}\subset M_\tau$, the spinor bundle
$W_\tau$ is by definition identified with
 the spinor bundle $W_{J_\tau}$ of $\gs_{\omega_\tau}$.
Let $(A,\Phi)$ be a solution of Seiberg--Witten equations on $M_\tau$.
The spinor $\Phi$ can be regarded as $(\beta, \gamma)\in W_{J_\tau}^+=
\Lambda^{0,0}\oplus \Lambda^{0,2} $.

By \fullref{propexpdecay}, there is a $T$ large enough,
 independent of $\tau$ and $(A,\Phi)$, such
 that, say, $|\beta| \geq 1/2$ on $M_\tau\setminus K$ where $K=\{T\leq
 \sigma_\tau\}$.
Hence we may define the map
$$h_{(A,\Phi)}\co  M_\tau\setminus K \rightarrow S^1$$
$$ h_{(A,\Phi)} = \frac{|\beta|}\beta . \leqno{\hbox{by}}
$$
By assumption, $(A,\Phi)\in\cC_l$ with $l\geq 4$. Therefore
$1-h_{A,\Phi}\in L^2_{l-1}$; using the Sobolev inclusion $L^2_3\subset
C^0$, we see that $h_{(A,\Phi)}$ tends to $1$ in $C^0$--norm near
infinity. So $h_{(A,\Phi)}$ is homotopic
to $1$ near infinity, and by assumption~\ref{assumption}, it implies
that the restriction of $h_{(A,\Phi)}$	to the annulus $\{T\leq
\sigma_\tau\leq\tau\}$ is
homotopic to $1$.

We define now a spinor bundle $W_{(A,\Phi)}$ on $M_\tau'$ as follows:
$$\left \{
\begin{array}{ll}
W_{(A,\Phi)}:= W &  \quad  \mbox{over } M'_\tau\cap \{\sigma'_\tau
<\tau \}\subset
M \\
W_{(A,\Phi)}:=W_{J'_\tau} & \quad \mbox{over }	M'_\tau\cap \{\sigma'_\tau
>
T \}\subset
 Z'
\end{array}
\right .
$$
and the transition map from $W$ to $W_{J'_\tau}$ is given by
$h'_{(A,\Phi)}=
h_{(A,\Phi)} \circ h$ over the annulus
$\{T< \sigma_\tau <\tau\}$.

The spinor bundle $W_{(A,\Phi)}$ together with its preferred
identification with $W_{J'_\tau}$ on the end of $M'_\tau$ define an
element of $\spinc(M'_\tau,\omega'_\tau)$.
In fact, this is nothing else but $j'(\gs,h)$,
 since $h_{(A,\Phi)}$ is homotopic to $1$ over the patching region.

The construction is compatible with the the gauge group action, in
the sense that for every  $u\in\cG(M_\tau)$, then
$$h_{u\cdot(A,\Phi)} = u^{-1}h_{(A,\Phi)};
$$
therefore $u$ induces an isomorphism $u^\sharp\co W_{(A,\Phi)}\rightarrow
W_{u\cdot(A,\Phi)}$ equal to $1$ on the end once the spinor bundles
are identified with $W_{J'_\tau}$.

\subsubsection{Definition of the approximate solutions}
\label{secpreg}
We suppose from now on that $\tau>T$ as in the previous section, so that
for any
solution of  Seiberg--Witten equations $(A,\Phi)$ on $M_\tau$, we can
construct the spinor bundle $W_{(A,\Phi)}$ on $M'_\tau$.

The bundles $W_{(A,\Phi)}$ and	$W_\tau$ restricted to the
region $K_T=\{\sigma \leq\tau\}\subset M$ are equal by construction. Hence
 $(A,\Phi)$  defines a configuration for the
spinor bundle $W_{(A,\Phi)}|_{K_T}$. We explain now how to extend it to
the entire
 manifold $M'_\tau$.

As before, we express $\Phi$ as a pair $(\beta,\gamma)\in
\Lambda^{0,0}\oplus \Lambda^{0,2}$ using the identification between
$W_\tau$  and $W_{J_\tau}$ on the end on $M_\tau$.
By construction, the bundle $W_{(A,\Phi)}$ is identified via
$h_{(A,\Phi)}$ to
$W_{J_\tau'}$ on the end $\{T < \sigma_\tau'\}$ furthermore
 and we can write modulo this
identification
 $$ h'_{(A,\Phi)}\cdot (A,\Phi) = h_{(A,\Phi)}\cdot
 (B+a,(\beta,\gamma)):=  (B+\what  a, (\what  \beta, \what  \gamma)).$$
The main effect of the	isomorphism  $h_{(A,\Phi)}$, is that
$\what  \beta$ is now a real function and $\what	\beta \geq \frac
12$;  hence the following definition makes sense:
$$(A,\Phi)^{\sharp} = (B' + \chi_\tau \what  a, (\what \beta^{\chi_\tau},
\chi_\tau \what \gamma)).
$$
Here $\chi_{\tau}$ is the function defined in Equation~\eqref{eq:cutoff}.
This extends naturally as a configuration relative to the spinor
bundle $W_{(A,\Phi)}$  on
$M'_\tau$ by setting $(A,\Phi)^{\sharp}:=(B,\Psi)$ on the end
$\{\sigma'_\tau\geq\tau\}\subset M_\tau'$.
We will also use the notation
$(A^\sharp,\Phi^\sharp):=(A,\Phi)^\sharp$. We stress on the fact that
if $(\wtilde A,\wtilde \Phi)$ is another solution of Seiberg--Witten
equations on
$M_\tau$, then $(A,\Phi)^\sharp$ and $(\wtilde A,\wtilde \Phi)^\sharp$
are not
\emph{a priori} defined  on the same bundles. However
this construction is compatible with the gauge group action in the
sense that
$$  u^\sharp\cdot(A,\Phi)^\sharp = (u\cdot(A,\Phi))^\sharp ;
$$
thus we have defined a map, called  the \emph{pregluing map},
$$\fbox{$\sharp\co \cM_{\varpi_\tau}(M_\tau)\rightarrow \cC/\cG(
M'_\tau).$}$$
We will see shortly that in fact this is a smooth map.
\begin{remark}\label{rkglue}
  In case we want to remove assumption~\ref{assumption} for the AFAK
  end $Z'$ and $Z$ is just a symplectic cone, we can construct
  similarly a
  pregluing map
$$\sharp\co \wtilde \cM_{\varpi_\tau}(M,(\gs',h'))\rightarrow
  \cC/\cG( M'_\tau,(\gs',h')),$$
where $(\gs',h')\in\spinc(M'_\tau,\omega_\tau')$ and  $\wtilde
  \cM_{\varpi_\tau}(M_{\tau},(\gs',h'))$ is the enlarged moduli space
  defined
  at~\eqref{eqbigmoduli}. The reader can check that the gluing
  theory applied starting from this pregluing map  leads to the result
  mentioned in \fullref{rkassumption}.
\end{remark}

\subsection{Rough gauge fixing on the target}
\label{secrough}
Given two solutions of Seiberg--Witten equations $(A,\Phi)$ and
$(\wtilde A,\wtilde \Phi)$ on $M_\tau$, the approximate solutions
$(A,\Phi)^\sharp$ and $( \wtilde A, \wtilde \Phi)^\sharp$
 are defined on different bundles $W_{(A,\Phi)}$
  and  $W_{(
  \wtilde A, \wtilde \Phi)}$ on $M'_\tau$. In this section, we explain
  how to
 construct a preferred	isomorphism provided $(A,\Phi)$ and
 $(\wtilde A,\wtilde \Phi)$ are close enough.

\subsubsection{Definition}
Using the identification between $W_\tau$ and $W_{J_\tau}$  on the end $
\{\sigma_\tau \geq T\}\subset M_\tau $,
we may write
$\Phi=(\beta,\gamma)$, $ \wtilde \Phi =( \wtilde \beta,\wtilde \gamma)\in
\Lambda^{0,0}\oplus \Lambda^{0,2}$, $A=B+a$ and $\wtilde A=B+\wtilde a$.

In order to simplify notation, we may assume after applying a gauge
transformation
that $(A,\Phi)$ is already in a gauge with exponential decay as in
\fullref{gaugeexpdecay}.
In particular, we have	$\beta=|\beta|$ hence $h_{(A,\Phi)}=1$,
therefore $W_{(A,\Phi)}$ is
equal to the spinor bundle $W'_\tau$
constructed in \fullref{secextspinc}.

If we assume that $\wtilde \Phi$ is sufficiently close to $\Phi$ in
$C^0$--norm,
then
$h_{( \wtilde A,\wtilde \Phi)}
= |\wtilde \beta|/ \wtilde \beta$ is also close to $1$;
in particular, we can write
$h( \wtilde A, \wtilde \Phi) = \exp(-\wtilde v)$, where $\wtilde v$ is a
purely imaginary
function completely determined	by the requirement  $|\wtilde v|<
 \pi$.	With this notation, we have $\wtilde \beta =
  \exp(\wtilde u+\wtilde v)$, where $\wtilde u$ is a real function such that
$|\wtilde \beta|=\exp \wtilde u$.

Put
 $$k_{(\wtilde A,\wtilde \Phi)}:= \exp ( -\chi_\tau  \wtilde v) ;$$
this is an isomorphism of $W_{J'_\tau}$, restricted to the annulus
$M_\tau\cap\{T<\sigma'_\tau<\tau \}$, equals to
$h_{(\wtilde A,\wtilde \Phi)}$ in a neighborhood
of $\{\sigma'_\tau=T\}$, and  to $1$ along
$\{\sigma'_\tau=\tau\}$. Therefore, we can extend $k_{(\wtilde A,\wtilde
\Phi)} $
to an isomorphism
$$k_{(\wtilde A,\wtilde \Phi)}\co  W_{(A,\Phi)} \to  W_{(\wtilde A,\wtilde
\Phi)} $$
 by setting
\begin{itemize}
\item
$k_{(\wtilde A,\wtilde \Phi)}:=\id|_{W_{J'_\tau}} $  for $\sigma_\tau'\geq
\tau$ and
 \item	$k_{(\wtilde A,\wtilde \Phi)}:=\id|_{W} $ for $\sigma'_\tau\leq T$.
\end{itemize}

\subsubsection{Estimates for the pregluing in rough gauge}
$(A,\Phi)^\sharp$ and
$k_{(\wtilde A,\wtilde \Phi)}^{-1}\cdot(\wtilde A,\wtilde \Phi)^\sharp$
are two configurations
defined  with respect
to the same
spinor bundle $W_{(A,\Phi)}$.
We are going to show that, if $(\wtilde A,\wtilde \Phi)$ is in Coulomb
gauge with
respect to $(A,\Phi)$, then, $k_{(\wtilde A,\wtilde \Phi)}^{-1}\cdot
(\wtilde A,\wtilde \Phi)^\sharp$ is
in some sense very close to be in Coulomb gauge with respect to
$(A,\Phi)^\sharp$.  We will also prove similar estimates for the
linearized Seiberg--Witten equations.

Since the isomorphism $k_{(\wtilde A,\wtilde \Phi)}$ apparently differs
from  the
identity only on the annulus $\{T\leq\sigma'_\tau\leq \tau\}$, we may
focus our study of the
pregluing map  on this region. With our
particular choice of gauge for $(A,\Phi)$, we have $|\beta|=\beta$ hence
$W_{(A,\Phi)}=W'_\tau$; using the the identification $W'_\tau\simeq
  W_{J'_\tau}$ for $\sigma'_\tau \geq T$, we have by definition of
  $\sharp$,
\begin{equation}\label {eqrealgauge}
 (A,\Phi)^\sharp = (B'+\chi_\tau a ,
\exp(\chi_\tau u), \chi_\tau \gamma),
\end{equation}
where $u$ is a real function such that $\beta =\exp u$.
\begin{remark} \label{rqexpdecay}
The fact that $(A,\Phi)$ is in a gauge with exponential decay together
with the identity~\eqref{eqrealgauge} show that for some constant
$c>0$, we have for every $N_0\geq 1$ (cf definition~\eqref{eq:cutoff}
of $\chi_\tau$), every $\tau$ large enough and every solution
$(A,\Phi)$ of Seiberg--Witten equations on $M_\tau$
$$ \chi_\tau \Bigl | (A,\Phi)^\sharp -(A,\Phi)\Bigr |\leq c \chi_\tau
(1-\chi_{(\tau - N_0)}) e^{-\delta \sigma_\tau}.
$$
Furthermore, similar estimates hold for all the covariant derivatives of
$\chi_\tau ((A,\Phi)^\sharp -(A,\Phi))$. Notice that the constants
involved do not
depend on $N_0\geq 1$ either, for the simple reason that the derivatives
of
$\chi_\tau$ are uniformly bounded as  $N_0 $ varies.
\end{remark}

Indeed, the map $\sharp$ provides some \emph{very good}
approximate solutions of Seiberg--Witten equations on $M'_\tau$. This
claim is made precise in the next lemma. Beforehand, it will be
convenient to package the equations on $M'_\tau$  into a single equation:
$$\SW( A', \Phi') = (F_{ A'}^+	-\{ \Phi' \otimes  (\Phi')^*\}_0  -
F_{B'_\tau}^+  +
\{ \Psi'_\tau \otimes (\Psi'_\tau)^* \}_0 - \varpi_\tau,\Dirac_{ A'}
\Phi').$$
\begin{lemma}
\label{lemmagoodapprox}
There exist $\delta >0$ and $T$ large enough such that for every
$N_0\geq 1$,
$k\in\NN$, $\tau\geq T + N_0$ and every
 solutions $(A,\Phi)$ of Seiberg--Witten
equations on $M_\tau$, we have
$$ \SW(A,\Phi)^\sharp = 0 \mbox{ on } \{\sigma'_\tau \leq T\} \subset
M'_\tau $$
and
$$  |\SW(A,\Phi)^\sharp|_{C^k(g_\tau',A^\sharp)} \leq c_k
e^{-\delta \sigma_\tau} \mbox{ on } \{\sigma'_\tau \geq T\} \subset
M'_\tau $$
\end{lemma}
\begin{proof}
It is trivial that $ \SW(A,\Phi)^\sharp =  \SW(A,\Phi)= 0$ on
$\{\sigma_\tau'\leq \tau-N_0\}$ and  $ \SW(A,\Phi)^\sharp$  $=
\SW(B',\Psi')= 0$ on $\{\sigma'_\tau\geq\tau  \}$.

If we put $(A,\Phi)$ in a gauge with exponential
decay, we obtain the estimate of the lemma thanks to
\fullref{rqexpdecay}. It follows that the lemma
 is true for any gauge representative.
\end{proof}

A direct
computation shows that
$$ k_{(\wtilde A,\wtilde \Phi)}^{-1}\cdot (\wtilde A,\wtilde \Phi)^\sharp =
\Bigl(B'+\chi_\tau \wtilde a
-\wtilde vd\chi_\tau,
\exp(\chi_\tau(\wtilde u+\wtilde v)),
\chi_\tau \wtilde \gamma\exp((\chi_\tau -1)\wtilde v)\Bigr),
$$
hence
\begin{equation}
\label{eqdiffrough}
 k_{(\wtilde A,\wtilde \Phi)}^{-1}\cdot (\wtilde A,\wtilde \Phi)^\sharp -
 (A,\Phi)^\sharp =
(\chi_\tau
(\wtilde a-a) +\zeta_1, \chi_\tau (\wtilde \beta-\beta) +\zeta_2, \chi_\tau
(\wtilde \gamma-\gamma) +\zeta_3),
\end{equation}
where
\begin{align}
\label{eq:zeta}
\zeta_1 =& -\wtilde vd\chi_\tau, \\
 \zeta_2 = &\exp(\chi_\tau (\wtilde u+\wtilde v) )-
\exp(\chi_\tau u ) - \chi_\tau(\exp  (\wtilde u+\wtilde v) - \exp u  ),
 \nonumber \\
 \zeta_3= &  \chi_\tau \wtilde \gamma (\exp((\chi_\tau -1)\wtilde v )
 -1), \nonumber
\end{align}
are some `smaller' terms. More precisely, $\zeta_j$ is controlled
by $\wtilde \beta-\beta$ thanks to \fullref{lemmatech}.
A direct consequence is
that $\sharp$ is
uniformly locally Lipschitz in the following sense:
\begin{lemma}
\label{gluingestimate0}
There exist $c, \alpha>0$, such that,
for every $N_0\geq 1$, every  $\tau$ large enough and every solution
of Seiberg--Witten
equations $(A,\Phi)$ on $M_\tau$, we have for every configuration
$(\wtilde A,\wtilde \Phi)$ on $M_\tau$
\begin{multline*}
 \|(\wtilde A,\wtilde \Phi) -
(A,\Phi)\|_{L^2_3 (g_\tau,A)}\leq \alpha \Rightarrow \\
\|  k_{(\wtilde A,\wtilde \Phi)}^{-1}\cdot (\wtilde A,
\wtilde \Phi)^\sharp-(A,\Phi)^\sharp\|_{L^2_k(g'_\tau,A^\sharp)}\leq
c \| (\wtilde A, \wtilde \Phi)-(A,\Phi)\|_{L^2_k(g_\tau,A)}\\ \qquad\qquad
\text{for all }
k=0,1,2,3.
\end{multline*}
\end{lemma}

A key estimate in order to study  the local injectivity of the
(pre) gluing is given in the next Lemma.
\begin{lemma}
\label{gluingestimate}
For every  $N_0$ large enough and  $\epsilon >0$, there exists
$\alpha>0$, such that,
for every $\tau$ large enough, and, every pair of Seiberg--Witten
solutions $(A,\Phi)$,  $(\wtilde A, \wtilde \Phi)$ on $M_\tau$ with the
conditions
$$ \delta_{1,(A,\Phi)}^*(\wtilde A-A,\wtilde \Phi-\Phi) = 0
\mbox{ and } \|(\wtilde A,\wtilde \Phi) -
(A,\Phi)\|_{L^2_3 (g_\tau,A)}\leq \alpha,$$
then
\begin{equation}\label{contrdelta1}
\left \|\delta_{1,(A,\Phi)^{\sharp}}^*\left (k_{(\wtilde A,\wtilde
\Phi)}^{-1}\cdot (\wtilde A,
\wtilde \Phi)^\sharp-(A,\Phi)^\sharp \right
)\right \|_{L^2_1(g'_\tau)}\leq \epsilon
\|( \wtilde A, \wtilde \Phi)-(A,\Phi)\|_{L^2_2(g_\tau,A)},
\end{equation}
\begin{equation}\label{contrdelta2}
\left \|\delta^{\null}_{2,(A,\Phi)^{\sharp}}\left ( k_{(\wtilde A,\wtilde
\Phi)}^{-1}\cdot (\wtilde A,
\wtilde \Phi)^\sharp-(A,\Phi)^\sharp \right )\right \|_{L^2_1(g'_\tau,
  A^\sharp)}\leq \epsilon
\|( \wtilde A, \wtilde \Phi)-(A,\Phi)\|_{L^2_2(g_\tau,A)}.
\end{equation}
\end{lemma}
\begin{proof}
We will only prove the estimate~\eqref{contrdelta2}.
The operator $\delta_1^* $  is	dealt with in completely similar way.

In order to prove the estimate~\eqref{contrdelta2}, we inspect the
contributions  of each
term of identity~\eqref{eqdiffrough}.
Let $\epsilon>0$.
We estimate the contribution of the terms $\zeta_j$ defined
at~\eqref{eq:zeta}
 first.
According to \fullref{boundedderivatives}, there exists a constant
$c_1>0$ independent of $\tau$ and the solution $(A,\Phi)$ of
Seiberg--Witten equations on $M_\tau$, such that
$$ \|\delta_{2,(A,\Phi)^\sharp} \zeta_1\|_{L^2_1(g'_\tau,A^\sharp)} \leq
c_1\|\zeta_1\|_{L^ 2_2(g'_\tau, A^\sharp)}.
$$
The derivatives of $\chi_\tau$ can be made arbitrarily small for $N_0$
large enough (see~\eqref{eq:cutoff}) so that
$$ \| \zeta_1\|_{L^2_2(g'_\tau,A^\sharp)} \leq \frac \epsilon{16c_1}
 \|\wtilde v\|_{L^ 2_2(g_\tau,\{T<\sigma_\tau<\tau \})},
$$
where $L^ 2_2(g_\tau,\{T<\sigma_\tau<\tau\})$ is the $L^2_2$ norm taken
on the
open set $\{T<\sigma_\tau<\tau\}$.
Using the uniform exponential decay of $(A,\Phi)$, we
deduce that for every $\tau$ large enough, we have, say
$$  \| \beta - 1 \|_{L^2_3(g_\tau, \{T<\sigma_\tau<\tau \}
)} \leq \alpha,
$$
hence
$$  \| \wtilde \beta - 1 \|_{L^2_3(g_\tau, \{T<\sigma_\tau<\tau \}
)} \leq 2\alpha.
$$
For $\alpha$ small enough and $\tau$ large enough,
\fullref{corexp} apply hence
$$ \|u\|_{L^2_3(g_\tau, \{T<\sigma_\tau<\tau \}
)}, \|\wtilde u+\wtilde v\|_{L^2_3(g_\tau, \{T<\sigma_\tau<\tau \}
)}  \leq 4\alpha.
$$
Thank to \fullref{lemmatech} we have the estimate
$$ \|\wtilde v + \wtilde u  - u \|_{L^2_2(g_\tau, \{T<\sigma_\tau<\tau \}
)} \leq 2  \|\exp(\wtilde v + \wtilde u)	-\exp u \|_{L^2_2(g_\tau,
\{T<\sigma_\tau<\tau \}
) }.
$$
Using the fact, the fact that $A$ decays exponentially fast toward $B$
again, the RHS is a lower bound for $2\|(\wtilde A,\wtilde
\Phi)-(A,\Phi)\|_{L^2_2(g_\tau,A)}$ for every $\tau$ large enough.

Since $\wtilde v =\Im( \wtilde u +\wtilde v - u)$, we have
$$ \|\wtilde v  \|_{L^2_2(g_\tau, \{T<\sigma_\tau<\tau \}
)} \leq 4 \|(\wtilde A,\wtilde
\Phi)-(A,\Phi)\|_{L^2_2(g_\tau,A)}.
$$
In conclusion we have we have shown that
\begin{equation}
\label{estglu1}
 \|\delta_{2,(A,\Phi)^\sharp} \zeta_1\|_{L^2_1(g'_\tau,A^\sharp)} \leq
\frac \epsilon 4 \|(\wtilde A,\wtilde
\Phi)-(A,\Phi)\|_{L^2_2(g_\tau,A)}.
\end{equation}
Similar estimates    show that for  $\alpha$
small enough, we also have
\begin{multline}
\label{estglu2}
\|(\wtilde A,\wtilde \Phi) -(A,\Phi)\|_{L^2_3(g_\tau,A)}\leq \alpha
\Rightarrow \\
\|\delta_{2,(A,\Phi)^\sharp}^{\null} \zeta_j \|_{L^2_1(g'_\tau, A^\sharp)}
\leq \frac \epsilon 4  \|(\wtilde A,\wtilde \Phi)
-(A,\Phi)\|_{L^2_2(g_\tau,A)},
\end{multline}
for $j=2,3$.

We need  to estimate the contribution of the remaining linear terms
displayed
in~\eqref{eqdiffrough}:
\begin{multline*} \delta_{2,(A,\Phi)^\sharp}^{\null} \Bigl
  (\chi_\tau((\wtilde A,\wtilde \Phi)
-(A,\Phi))\Bigr ) = \chi_\tau \delta_{2,(A,\Phi)^\sharp}^{\null} \Bigl
( (\wtilde A,\wtilde \Phi)
-(A,\Phi)\Bigr ) \\
+\Bigl ( \left (d\chi_\tau
\wedge (\wtilde A-A)\right )^+, d\chi_\tau \cdot (\wtilde \Phi -
\Phi)\Bigr ).
\end{multline*}
Notice that $\delta_{2,(A,\Phi)^\sharp}^{\null} \Bigl ( (\wtilde A,\wtilde
\Phi)
-(A,\Phi)\Bigr )$ makes sense only on the domain $\{\sigma_\tau'\leq
\tau\} = \{\sigma_\tau\leq
\tau\}$ where all bundles, metrics, etc\dots are equal.

Similarly to $\zeta_1 =\wtilde v d\chi_\tau$, we can by increasing
 $N_0$ make $d\chi_\tau$
 as small as required in order to  have the estimate
\begin{equation}
\label{estglu3}
 \left \| \left (\left (d\chi_\tau
\wedge (\wtilde A-A)\right )^+,	d\chi_\tau \cdot  (\wtilde \Phi -
\Phi)
\right )\right \|_{L^2_1(g'_\tau, A^\sharp)} \leq \frac \epsilon {12}
\|(\wtilde A,\wtilde \Phi) -(A,\Phi)\|_{L^2_1(g_\tau,A)}.
\end{equation}
The operator $ \chi_\tau\Bigl
(\delta_{2,(A,\Phi)^\sharp}^{\null} -
 \delta_{2,(A,\Phi)}^{\null}\Bigr )$ is  linear operator of degree
 $0$ and
we have
$$  \chi_\tau\Bigl (\delta_{2,(A,\Phi)^\sharp}^{\null} -
 \delta_{2,(A,\Phi)}^{\null}\Bigr ) \Bigl
( (\wtilde A,\wtilde \Phi)
-(A,\Phi)\Bigr )  =  \chi_\tau \Bigl \{ (A,\Phi)^\sharp - (A,\Phi),
(\wtilde A,\wtilde
 \Phi) - (A,\Phi)\Bigr \},
$$
where $\{\cdot,\cdot\}$ is some bilinear pairing whose explicit
expression is irrelevant to us. The uniform exponential decay property of
$\chi_\tau((A,\Phi)^\sharp - (A,\Phi))$ mentioned in
\fullref{rqexpdecay} show that for $\tau$ large enough, we have for
every pair of solutions $(A,\Phi)$,  $(\wtilde A,\wtilde \Phi)$ of
Seiberg--Witten equations
\begin{multline*}
 \left \| \chi_\tau\Bigl (\delta_{2,(A,\Phi)^\sharp}^{\null} -
 \delta_{2,(A,\Phi)}^{\null}\Bigr ) \Bigl
( (\wtilde A,\wtilde \Phi)
-(A,\Phi)\Bigr ) \right \|_{L^2_1(g'_\tau,A^\sharp)} \\
\leq  \frac \epsilon {12}
 \|(\wtilde A,\wtilde \Phi) -(A,\Phi)\|_{L^2_1(g_\tau,A)}.
\end{multline*}
Eventually, using the fact that $(A,\Phi)$ and $(\wtilde A,\wtilde
\Phi)$ are both solution of Seiberg--Witten equations we have the identity
$$ \delta_{2,(A,\Phi)}^{\null} \Bigl
( (\wtilde A,\wtilde \Phi)
-(A,\Phi)\Bigr ) + \cQ	\Bigl
( (\wtilde A,\wtilde \Phi)
-(A,\Phi)\Bigr )=0,
$$
where $\cQ$ is the quadratic term
\begin{equation}\label{eq:quaddef2}
\cQ(a,\phi) = ( \{\phi\otimes \phi^*\}_0,  a\cdot\phi ).
\end{equation}
The Sobolev multiplication theorem $L^2_2\otimes L^2_2\hookrightarrow
L^2_1$ on every compact set
together Lemmas~\ref{lemmaboundedgeometry} and
\ref{boundedderivatives}
tell us that there exists a constant $c_2>0$ independent of $(A,\Phi)$
and $(\wtilde A,\wtilde \Phi)$, such that for every
$\tau$ large enough we have
$$ \|\cQ  \Bigl
( (\wtilde A,\wtilde \Phi)
-(A,\Phi)\Bigr )\|_{L^2_1(g_\tau,A)} \leq  c_2 \|(\wtilde A,\wtilde \Phi)
-(A,\Phi)\Bigr )\|^2_{L^2_2(g_\tau,A)}.
$$
So, if we choose $\alpha$ small enough, we have
\begin{equation}
\label{estglu4}
\| \chi_\tau \delta_{2,(A,\Phi)}^{\null} \Bigl
( (\wtilde A,\wtilde \Phi)
-(A,\Phi)\Bigr )\|_{L^2_1(g'_\tau, A^\sharp)} \leq \frac \epsilon{12}
\|(\wtilde A,\wtilde \Phi)
-(A,\Phi)\Bigr )\|^2_{L^2_2(g_\tau,A)}.
\end{equation}
Remark that the above estimate is true for every $N_0\geq 1$, since
$\|d\chi_\tau\|_{C^0}$ decays as $N_0$ increases.
Adding up the estimate~\eqref{estglu1}, \eqref{estglu2},
\eqref{estglu3} and \eqref{estglu4}, we conclude
$$\left \|\delta^{\null}_{2,(A,\Phi)^{\sharp}}\left ( k_{(\wtilde A,\wtilde
\Phi)}^{-1}\cdot (\wtilde A,
\wtilde \Phi)^\sharp-(A,\Phi)^\sharp \right )\right \|_{L^2_1(g'_\tau,
  A^\sharp)}\leq \epsilon
\|( \wtilde A, \wtilde \Phi)-(A,\Phi)\|_{L^2_2(g_\tau,A)}.
$$
which proves the lemma.
\end{proof}

\section{Gluing}
Our goal is to perturb
$ \sharp\co \cM(M_\tau)\longrightarrow \cC/\cG(M'_\tau)$, in order to get
a \emph{gluing map}
$$\gG \co \cM_{\varpi_\tau}(M_\tau)\longrightarrow \cM_{\varpi_\tau}
(M'_\tau).$$
Moreover we will show that $\gG$ is a diffeomorphism.

\subsection{Construction of the gluing map}
We seek a solution of Seiberg--Witten equations of the
form
\begin{equation}\label{eqglue0}
 (A', \Phi')=(A,\Phi)^\sharp +
\delta_{2,{(A,\Phi)^\sharp}}^*(b',\psi'),
\end{equation}
where $b'$ is a purely imaginary self-dual form and $\psi'$ a section of
the negative spinor bundle $W^-_{(A,\Phi)}$ on $M'_\tau$.
The Seiberg--Witten equations for $(A',\Phi')$ viewed as equations
bearing on $(b',\psi')$ have the form
\begin{equation}
\label{eqglue}
   \Delta_{2,(A,\Phi)^\sharp} (b',\psi')
+ \cQ ( \delta_{2,(A,\Phi)^\sharp}^*(b',\psi')) = - \SW(A,\Phi)^\sharp,
\end{equation}
where
$$\Delta_{2,(A,\Phi)^\sharp}:=
\delta_{2,(A,\Phi)^\sharp}\delta_{2,(A,\Phi)^\sharp}^*$$
 and $\cQ(a',\phi')$ is nonlinear part defined at~\eqref{eq:quaddef2}.

\subsubsection{The linear problem}
As usual we start by solving the linear problem.
\begin{prop}
\label{propfirstev}
Assume that the moduli space, $\cM_{\varpi}(M)$ is unobstructed.
Then, for each $k\geq 0$,  there exists a constant $c_k>0$, such that
for every $\tau$ large
enough,
every $N_0\geq 1$ (see definition of $\chi_\tau$)
and  every solution $(A,\Phi)$ of Seiberg--Witten equations
on $M_\tau$, we have  for every pairs $(b,\psi)$ on $M_\tau$ and
$(b',\psi')$ on $M'_\tau$
\begin{equation}
\label{contrlapl2}
 \|\Delta_{2,(A,\Phi)}(\psi,b)\|_{L^2_k(g_\tau,A)}\geq c_k
  \|(\psi,b)\|_{L^2_{k+2}(g_\tau,A)},
\end{equation}
where
$\Delta_{2,(A,\Phi)}:=\delta^{\null}_{2,{(A,\Phi)}}\delta_{2,{(A,\Phi)}}
^*$, and
\begin{equation}
\label{contrlapl}
 \|\Delta_{2,(A,\Phi)^\sharp}(b',\psi')\|_{L^2_k(g'_\tau,A)}\geq c_k
  \|(b',\psi')\|_{L^2_{k+2}(g'_\tau,A^\sharp)}.
\end{equation}
\end{prop}
\begin{proof}
We will just prove the second statement \eqref{contrlapl}, since the
proof of first part of the proposition is the same with less
complications.

Let $\BBox_{(A,\Phi)^\sharp}$ be the operator
$\cD^{\null}_{(A,\Phi)^\sharp} \cD_{(A,\Phi)^\sharp}^*$, where
  $\cD^{\null}_{(A,\Phi)^\sharp} = \delta^*_{1,(A,\Phi)^\sharp}\oplus
\delta^{\null}_{2,(A,\Phi)^\sharp}$ as defined in
\fullref{seclintheo}.
Suppose that for some constant $c>0$, we have for  every $\tau$ large
enough,
every $N_0\geq 1$, we have the inequality
\begin{equation}
\label{eqopcontr}
\|\BBox_{(A,\Phi)^\sharp}(u',b',\psi')\|_{L^2(g'_\tau)}\geq c
 \|(u',b',\psi')\|_{L^2_{2}(g'_\tau,A^\sharp)} \quad \text{for all }
 (u',b',\psi').
\end{equation}
Then, the estimate on the operator $\Delta_{2,(A,\Phi)^\sharp}$
  follows. To see it we write
\begin{align}
\label{eqopsplit}
\BBox_{(A,\Phi)^\sharp}(u',b',\psi') &=
\Bigl (\Delta_{1,(A,\Phi)^\sharp} u' +
\cL^*_{(A,\Phi)^\sharp}(b',\psi') \Bigr )\\
&\hspace{1in}\oplus \Bigl (\Delta_{2,(A,\Phi)^\sharp}(b',\psi') +
\cL_{(A,\Phi)^\sharp} u'\Bigl )\notag\\
 \cL_{(A,\Phi)^\sharp}&=\delta_{2,(A,\Phi)^\sharp}
\delta_{1,(A,\Phi)^\sharp} .\tag*{\hbox{where}}
\end{align}
If $(A,\Phi)^\sharp$  were an exact  solution of Seiberg--Witten
equations, we would simply have $ \cL_{(A,\Phi)^\sharp} =0$ for
\begin{equation}\label{eqdefect}
  \cL_{(A,\Phi)^\sharp} u'  = (u'\Dirac_{A^\sharp}\Phi^\sharp, 0);
\end{equation}
nevertheless, the operator $\cL_{(A,\Phi)^\sharp}$ has operator norm
very close to $0$ since
$(A,\Phi)^\sharp$ is an approximate solution in the sense of
\fullref{lemmagoodapprox}. More precisely, we have for every
$N_0\geq 1$
and  every $\tau$ large enough
$$  \|\cL_{(A,\Phi)^\sharp} u'\|_{L^2(g'_\tau)} +
\|\cL^*_{(A,\Phi)^\sharp}
(b',\psi')\|_{L^2(g'_\tau)} \leq \frac	c2
\|(u',b',\psi')\|_{L^2({g'_\tau})}.
$$
It follows from \eqref{eqopsplit} and \eqref{eqopcontr} that
\begin{equation}
\label{eqopcontr2}
 \|\Delta_{2,(A,\Phi)^\sharp}(b',\psi')\|_{L^2(g'_\tau)} +
  \|\Delta_{1,(A,\Phi)^\sharp} u' \|_{L^2(g'_\tau)} \geq \frac c2
  \|(u',b',\psi')\|_{L^2_{2}(g'_\tau,A^\sharp)}.
\end{equation}
This gives the control that we wanted by setting $u'=0$.

It remains to show that \eqref{eqopcontr} holds to finish the proof of
the proposition. First we observe that the operator $\cD^*$ gives uniform
control on  sections supported near infinity.
\begin{lemma}\label{lem:Dstarest1}
There exist $T>0$ large enough and a constant $c>0$  such that for every
$\tau>T $, every $N_0\geq 1$ and every approximate solution
$(A,\Phi)^{\sharp}$ of Seiberg--Witten equations on
$M'_\tau$,
we have  for all sections $(u',b',\psi')$ with compact support in
$\{\sigma'_\tau >T\}\subset M'_\tau$
$$ \|\cD^*_{(A,\Phi)^\sharp}(u',b',\psi')\|_{L^2(g'_\tau)}\geq c
  \|(u',b',\psi')\|_{L^2_1(g'_\tau, A^\sharp)}.
$$
\end{lemma}
\begin{proof}
We use the Weitzenbock formula derived in \cite[Proposition 3.8]{KM} which
says that
\begin{multline*}
 \|\cD^*_{(A,\Phi)^\sharp}(u',b',\psi')\|_{L^2(g'_\tau)}  =\\
  \int_{\sigma'_\tau\geq T}
 \Big (  |du'|^2 + 4|\nabla b'|^2 +|\nabla_{A^\sharp} \psi'|^2 +
|\Phi^\sharp
 |^2(|u'|^2 + |b'|^2+ |\psi'|^2) \\
+ 4 \Bigl \langle \left ( \frac s6\id -\cW^+ \right )\cdot b',b' \Bigr
\rangle + \frac s4 |\psi'|^2 + \frac
 12\ip{F_{A^\sharp}^-\cdot \psi',\psi'} \\
 + 4\ip {b'\otimes \psi', \nabla_{A^\sharp}
 \Phi^\sharp  }  + 2\ip {\psi', u' D_{A^\sharp}\Phi^\sharp} + 2\ip
 {b',\psi'\cdot \Dirac_{A^\sharp} \Phi^\sharp }
\Big ) \vol ^{g'_\tau},
\end{multline*}
where $\cW^+$ is the positive part of the Weyl curvature.
Thanks to \fullref{lemmagoodapprox},   and the fact that the
geometries for all
$M_\tau'$ are uniformly controlled by \fullref{lemmaboundedgeometry},
we know that $s$, $\cW$,
$|\Phi^\sharp|-1$,
$F_{A^\sharp}$, $\nabla_{A^\sharp}\Phi^\sharp$,
$\Dirac_{A^\sharp}\Phi^\sharp$ are
uniformly close to zero. Therefore,  the terms of first line
control all the others and the lemma is proved.
\end{proof}

Next we prove that similar estimates hold for large $\tau$
and for global sections.
\begin{lemma}\label{lem:Dstarest}
There exists  $\kappa>0$ such that for every $\tau$ large
enough, for every $N_0\geq 1$, for every  solution of Seiberg--Witten
equations $(A,\Phi)$ and
for every (global)~section
$(u',b',\psi')$ on $M'_\tau$, we
have
$$ \|\cD^*_{(A,\Phi)^\sharp}(u',b',\psi')\|_{L^2(g'_\tau)}\geq \kappa
  \|(u',b',\psi')\|_{L^2_1(g'_\tau, A^\sharp)}.
$$
\end{lemma}
\begin{proof}
Suppose this is not true. Then, we have  sequences $N_{0,j}\geq 1$
and  $\tau_j\rightarrow
+\infty$ with $\tau_j\geq T+N_{0,j}$, some solutions of Seiberg--Witten
equations
$(A_{j},\Phi_{j})$ on $M_{\tau_j}$, and $(u'_{j},b_j',\psi_j')$ such that
$$\|\cD^*_{(A_j,\Phi_j)^\sharp}(u'_j,b'_j,\psi'_j)\|_{L^2(g'_{\tau_j})}\to
0,\quad
  \|(u_j',b'_j,\psi'_j)\|_{L^2_1(g'_{\tau_j},A^\sharp_j)}= 1.
$$
After applying gauge transformations and extracting  a
subsequence, we may assume that $(A_j,\Phi_j)$ converge on every
compact set to a solution of Seiberg--Witten equations
$(A_\infty,\Phi_\infty)$ on $M$. Then, using
the $L^2_1$ bound on $(u'_j,b'_j,\psi'_j)$, the sequence converges on
every compact subset to
 a weak limit $(u,b,\psi)\in L^2_1(g,A_\infty)$ on $M$ verifying
 $\cD^*_{(A_\infty,\Phi_\infty)} (u,b,\psi)=0$.
Using the assumption that the moduli space on $M$ is unobstructed,
this implies that $(u,b,\psi)=0$. Using the compactness of the
inclusion $L^2_{1}\subset L^2$ on compact sets, we see that the sequence
$(u'_j,b'_j,\psi'_j)$ converge strongly toward $0$ in the $L^2$--sense on
every compact set after further extraction.

Let $\chi$ be a cut-off function equal to $1$ on $\{\sigma'_\tau\leq
T\}$ and to $0$ outside $\{\sigma'_\tau\leq
T+1\}$. Then
\begin{equation}
\label{eqcutme}
 \|(u'_j,b'_j,\psi'_j)\|_{L^2(g_\tau')}\leq
\|\chi (u'_j,b'_j,\psi'_j)\|_{L^2(g_\tau')} +
\|(1-\chi)(u'_j,b'_j,\psi'_j)\|_{L^2(g_\tau')};
\end{equation}
the first term in the RHS tends to $0$ since $(u'_j,b'_j,\psi'_j)$
converges to $0$ on the compact set $\sigma'_\tau \leq {T+1}$. The
second term  is supported in $\sigma'_\tau \geq T$, hence we can
apply \fullref{lem:Dstarest1} which says that it is controlled by
$$
c\|\cD^*_{(A_j,\Phi_j)^\sharp}(1-\chi)(u'_j,b'_j,\psi'_j)\|_{L^2(g_{\tau_j}')}.
$$
The derivatives of $\chi$  are compactly supported; so the
Leibniz rule for $\cD^*$ and the fact that
$\|\cD^*_{(A_j,\Phi_j)^\sharp}(u'_j,b'_j,\psi'_j)\|_{L^2(g'_{\tau_j})}\to
  0$, shows that
$$\|\cD^*_{(A_j,\Phi_j)^\sharp}(1-\chi)(u'_j,b'_j,\psi'_j)\|_{L^2(g_{\tau_j}')}
\to 0 .$$ Therefore
$\|(1-\chi)(u'_j,b'_j,\psi'_j)\|_{L^2(g_{\tau_j}')} \to 0 $,
  hence
$\|(u'_j,b'_j,\psi'_j)\|_{L^2(g_{\tau_j}')} \to 0 $ by~\eqref{eqcutme}.
This is a contradiction and the lemma is proved.
\end{proof}

We are ready to prove that the estimate~\eqref{eqopcontr} holds.
Under the assumption of \fullref{lem:Dstarest}, we have
\begin{multline*}
  \int _ {M'_\tau } \ip{\BBox_{(A,\Phi)^\sharp}(u',b',\psi') ,
  (u',b',\psi')}
\vol^{g'_\tau } = \\
\|\cD^*_{(A,\Phi)^\sharp}(u',b',\psi')\|^2_{L^2(g'_\tau)}\geq \kappa^2
\|(u',b',\psi')\|_{L^2_1(g'_\tau,A^\sharp)}^2.
\end{multline*}
By Cauchy--Schwarz inequality, the LHS is a lower bound for
$$ \Bigl \|\BBox_{(A,\Phi)^\sharp}(u',b',\psi') \Bigr \|_{L^2(g'_\tau)}
\Bigl \|(u',b',\psi')\Bigr \|_{L^2_1(g'_\tau,A^\sharp)},
$$
and we conclude that
\begin{equation}
\label{contrlapl3}
 \|\BBox_{(A,\Phi)^\sharp}(u',b',\psi')\|_{L^2(g'_\tau)}\geq
\kappa^2\|(u',b',\psi')\|_{L^2_1(g'_\tau,A^\sharp)}.
\end{equation}
The operator $\BBox_{(A,\Phi)^\sharp}$ is elliptic of order
$2$. Hence, for some constant $c_2>0$, we have on a small  ball $B_\alpha
\subset
M'_\tau$ of radius $\alpha$, the estimate
\begin{multline*}
  \|\BBox_{(A,\Phi)^\sharp}(u',b',\psi')\|_{L^2(g'_\tau,B_\alpha)} +
  \|(u',b',\psi')\|_{L^2(g'_\tau,B_{\alpha})} \geq\\
 c_2  \|(u',b',\psi')\|_{L^2_2(g'_\tau,A^\sharp, B_{\alpha/2})}.
\end{multline*}
The fact that the approximate solution of Seiberg--Witten equations
 $(A,\Phi)^\sharp$ on $M'_\tau$ is controlled thanks to
 \fullref{boundedderivatives} and \fullref{rqexpdecay},
that the
geometry of the $M'_\tau$ is uniformly AFAK in the sense of
\fullref{lemmaboundedgeometry} implies that the constant
$c'_2$ can be chosen independently of $\tau$ ( provided it is large
enough), of the approximate solution $(A,\Phi)^\sharp$, of the
 constant $N_0\geq 1$ involved in our construction,
and of the center of the ball
$B_\alpha$. It follows that there exists a constant $c'_2>0$,
 independent of all data as $c_2$,  such that
$$ \|\BBox_{(A,\Phi)^\sharp}(u',b',\psi')\|_{L^2(g'_\tau)} +
  \|(u',b',\psi')\|_{L^2(g'_\tau)} \geq
 c'_2  \|(u',b',\psi')\|_{L^2_2(g'_\tau,A^\sharp)}.
$$
Together with \eqref{contrlapl3}, this implies that we have an
estimate of the form \eqref{contrlapl2}, thus
\fullref{propfirstev} is
proved for $k=0$. For $k\geq 1$, similar estimates hold by elliptic
regularity
as in the proof of \fullref{lemmacontrolP}.
\end{proof}
As an immediate consequence, we have the following Corollary:
\begin{cor}
\label{corunobstr}
If the Seiberg--Witten equations on $M$ are unobstructed, the
Seiberg--Witten equations on $M_\tau$ (or $M'_\tau$) are unobstructed
for every $\tau$ large enough.
\end{cor}
\begin{proof}
\fullref{propfirstev} togheter with \eqref{contrlapl2} shows in particular
that for every $\tau$
large enough, $\delta_{2,(A,\Phi)}$
is surjective for every $(A,\Phi)$, which means that
the moduli spaces $\cM_{\varpi_\tau}(M_\tau)$ and
$\cM_{\varpi_\tau}(M'_\tau)$ are
unobstructed.
\end{proof}
\begin{remark*}
In fact, we already assumed in \fullref{seclintheo} that the
Seiberg--Witten
equation on $M$ are unobstructed, by choosing a generic perturbation
$\varpi$. We will always suppose that it is the case from now on,
unless stated.
\end{remark*}
Here is another immediate consequence of \fullref{propfirstev}, directly
related to our gluing problem.
\begin{cor}
\label{corinv}
There exists a constant $c>0$ such that for every $\tau$ large enough,
every $N_0\geq 1$ and every solution $(A,\Phi)$ of Seiberg--Witten
equations on $M_\tau$, the operator
$\Delta_{2,(A,\Phi)^\sharp}\co L^2_2\to L^2$ is an isomorphism, and,
moreover,
its  inverse $\Delta_{2,(A,\Phi)^\sharp}^{-1}$ verifies
\begin{equation}\label{eqcorinv}
 c \|(b',\psi' )\|_{L^2_1(g'_\tau,A^\sharp)}\geq
  \| \Delta_{2,(A,\Phi)^\sharp}^{-1}(b',\psi')\|_{L^2_3(g'_\tau,A^\sharp)}
  \quad
  \text{for all } (b',\psi').
\end{equation}
\end{cor}

\subsubsection{The nonlinear problem}
Using the substitution
$V=\Delta_{2,(A,\Phi)^\sharp}(b',\psi')$, we can rewrite the equation
\eqref{eqglue} in the form
\begin{equation}
\label{eqglue2}
V = S_{(A,\Phi)^\sharp}(V) - \SW(A,\Phi)^\sharp,
\end{equation}
where
$$S_{(A,\Phi)^\sharp}(V)= - \cQ ( \delta_{2,(A,\Phi)^\sharp}^*
\Delta_{2,(A,\Phi)^\sharp}^{-1}V).$$
The key argument for solving Equation~\eqref{eqglue2} is that the
operator $S$ is a uniform contraction in the sense of the next lemma.
\begin{lemma}
\label{lemmacontract2}
There exist constants $\alpha>0$, $ \kappa\in(0,1/2)$, such that for
every  $\tau$ large
 enough, every $N_0\geq 1$ and every approximate solution of
 Seiberg--Witten equations on $M'_\tau$ we have
\begin{multline*}
\text{for all } V_1, V_2 \quad\|V_1\|_{L^2(g'_\tau)},
\|V_2\|_{L^2_1(g'_\tau,A^\sharp)}
\leq \alpha
\Rightarrow \\
  \|S_{(A,\Phi)^\sharp
}(V_2) - S_{(A,\Phi^\sharp)
}(V_2) \|_{L^2_1(g'_\tau,A^\sharp)}\leq
\kappa\|V_2-V_1\|_{L^2_1(g'_\tau,A^\sharp)}.
\end{multline*}
\end{lemma}
\begin{proof}
We choose $\tau$ large enough and a constant $c>0$ according to
\fullref{corinv}.
$$S_{(A,\Phi^\sharp)
}(V_2) - S_{(A,\Phi^\sharp)
}(V_2) = \cQ\Bigl ( \delta_{2,(A,\Phi)^\sharp}^*
\Delta_{2,(A,\Phi)^\sharp}^{-1}V_2\Bigr )-  \cQ\Bigl (
\delta_{2,(A,\Phi)^\sharp}^*
\Delta_{2,(A,\Phi)^\sharp}^{-1}V_1\Bigr ).
$$
Since $\cQ(a,\phi)$ is a quadratic polynomial in $(a,\phi)$, we deduce
using  the Sobolev multiplication Theorem $L^2_2\otimes
L^2_2\rightarrow L^2_1$, a control
\begin{multline}\label{estcs}
\left \|\cQ( \delta_{2,(A,\Phi)^\sharp}^*
\Delta_{2,(A,\Phi)^\sharp}^{-1}V_2 )-  \cQ( \delta_{2,(A,\Phi)^\sharp}^*
\Delta_{2,(A,\Phi)^\sharp}^{-1}V_1)\right \|_{L^2_1(g'_\tau,A^\sharp )}\\
\leq C\|
  \delta_{2,(A,\Phi)^\sharp}^*
\Delta_{2,(A,\Phi)^\sharp}^{-1}(V_1+V_2)\|_{L^2_2(g'_\tau,A^\sharp)} \|
 \delta_{2,(A,\Phi)^\sharp}^*
\Delta_{2,(A,\Phi)^\sharp}^{-1}(V_1-V_2) \|_{L^2_2(g'_\tau,A^\sharp)},
\end{multline}
where $C$ is a constant which  depend neither on  $\tau$ (large
enough) nor on	the
approximate solution of Seiberg--Witten equations $(A,\Phi)^\sharp$ on
$M'_\tau$, nor on
 the constant $N_0\geq 1$ involved in the construction of $\sharp$.

The bounds on the approximate solution $(A,\Phi)^\sharp$
deduced from those of $(A,\Phi)$,
show that there is a constant $C_2>0$ independent of
$\tau$ (large enough), $(A,\Phi)^\sharp$ and $N_0\geq 1$,
such that
\begin{equation}
\label{eqbound}
 \text{for all } (b',\psi'),\quad
\|\delta_{2,(A,\Phi)^\sharp}(b',\psi')\|_{L^2_2(g_\tau, A^\sharp)}\leq
C_2 \|(b',\psi')\|_{L^2_3(g_\tau, A^\sharp)}.
\end{equation}
In particular this holds for $(b',\psi')=
\Delta^{-1}_{2,(A,\Phi^\sharp)}(V_1\pm V_2)$. Together
with~\eqref{eqcorinv}
we deduce from~\eqref{estcs} the estimate
\begin{multline*}
\left \|\cQ( \delta_{2,(A,\Phi)^\sharp}^*
\Delta_{2,(A,\Phi)^\sharp}^{-1}V_2 )-  \cQ( \delta_{2,(A,\Phi)^\sharp}^*
\Delta_{2,(A,\Phi)^\sharp}^{-1}V_1)\right \|_{L^2_1(g'_\tau,A^\sharp)}\\
\leq  c C
C_2 \|
 V_2+V_1\|_{L^2_1(g'_\tau,A^\sharp)} \|
V_2-V_1\|_{L^2_1(g'_\tau,A^\sharp)} ,
\end{multline*}
and  the lemma holds for $\alpha = \frac \kappa {2cCC_2}$.
\end{proof}

In Since $S$ is contractant in a suitable sense and
$\SW(A,\Phi)^\sharp$ converges uniformly to $0$ in the sense of
\fullref{lemmagoodapprox}, we can solve
equation~\eqref{eqglue2} thanks to
\fullref{propcontractmapping} for $V\in L^2_1$. Then we obtain
a solution
of~\eqref{eqglue0} given by $(b',\psi' ) =
\Delta_{2,(A,\Phi)^\sharp}^{-1}V$, hence $(b',\psi')\in L^2_3$.
More precisely, we have the following theorem:

\begin{theo}[Definition of the gluing map]
\label{theogludef}
There exist  constants $\alpha,c >0$ such that for every $\tau$ large
 enough,
 every solution $(A,\Phi)$ of Seiberg--Witten equations on
$M_\tau$ and every constant $N_0\geq 1$ (see definition of $\sharp$),
there is a unique section $(b',\psi')$ on $M'_\tau$ such that
$$ \gG(A,\Phi) = (A,\Phi)^\sharp +
\delta^*_{2,(A,\Phi)^\sharp}(b',\psi')
$$
is a solution of Seiberg--Witten equations on $M'_\tau$, with
$\|(b',\psi')\|_{L^2_3(g,A^\sharp)}\leq \alpha$.

The map $\gG$ is smooth  and gauge equivariant and induces
 map
$$\gG\co \cM_{\varpi_\tau} (M_\tau) \to \cM_{\varpi_\tau}(M'_\tau),$$
furthermore
\begin{equation}\label{eqbound2}
\|(b',\psi')\|_{L^2_3(g_\tau,A^\sharp)}, \|\gG(A,\Phi) -
(A,\Phi)^\sharp\|_{L^2_2(g'_\tau,A^\sharp)}\leq
c\|\SW(A,\Phi)^\sharp\|_{L^2_1(g'_\tau,A^\sharp)}.
\end{equation}
\end{theo}
\begin{proof}
We have already proved the existence of a solution $(b',\psi')\in
L^2_3$. Therefore $\gG(A,\Phi)\in \cC_2(M'_\tau)$.
 If $(A,\Phi)\in \cC_{l}(M_\tau)$, we
 have $\SW(A,\Phi)^\sharp \in L^2_{l-1}$ and it follows by elliptic
 regularity  of Seiberg--Witten equations that $\gG(A,\Phi)\in
 \cC_{l}(M'_\tau)$. If we choose $l\geq 4$ (which is required to
 define the moduli space), we have then a
well defined induced map $\gG\co \cM_{l,\varpi_\tau} (M_\tau) \to
\cM_{l,\varpi_\tau}(M'_\tau)$.

The smoothness of  $\gG$
follows from the smoothness of $\sharp$  and from the fact that
$S_{(A,\Phi)^\sharp}$ in~\eqref{eqglue2} depends smoothly on the
parameter $(A,\Phi)^\sharp$.

The only part of the theorem left to be proved is the
estimate~\eqref{eqbound2}.
Recall that $V=\Delta_{2,(A,\Phi)}(b',\psi')$ where  $V$
is a solution of~\eqref{eqglue2} provided by
\fullref{propcontractmapping}. Therefore, $V$ verifies
$$ \| V \|_{L^2_1(g'_\tau,A^\sharp)} \leq 2
\|\SW(A,\Phi)^\sharp\|_{L^2_1(g_\tau, A^\sharp)}.
$$
Let $c_1>0$  be a constant obtained by
\fullref{propfirstev} such that
$$c_1 \|(b',\psi')\|_{L^2_3(g'_\tau,A^\sharp)}\leq   \|
\Delta_{2,(A,\Phi)^\sharp} (b',\psi')
\|_{L^2_1(g'_\tau)}= \| V \|_{L^2_1(g'_\tau,A^\sharp)}
$$
Similarly to estimate~\eqref{eqbound}, there is a constant $c_2>0$
 independent of $N_0\geq 1$, $\tau$ (large enough) and
 $(A,\Phi)^\sharp$, such that
\begin{equation}
\label{eqbound3}
 \text{for all } (b',\psi'),\quad
\|\delta^*_{2,(A,\Phi)^\sharp}(b',\psi')\|_{L^2_2(g_\tau, A^\sharp)}\leq
c_2 \|(b',\psi')\|_{L^2_3(g_\tau, A^\sharp)}.
\end{equation}
Eventually~\eqref{eqbound2} is verified for $c=\max(2c_2/c_1,2/c_1)$.
\end{proof}

\subsection{From local to global study	of the gluing map}
We will  prove that the gluing is locally injective in the following
sense:
\begin{prop}
\label{proplocinj}
There exist $N_0\geq 1$  (for the construction the
pregluing map) large enough, and a  constants $\alpha_2>0$,
such that for every $\tau$ large enough and every solution
$(A,\Phi)$ of Seiberg--Witten equations on $M_\tau$, the gluing map
$ \gG\co \cM_{\varpi_\tau} (M_\tau)\to\cM_{\varpi_\tau} (M_\tau') $
restricted to the open set
\begin{multline*}
 B({[A,\Phi]},\alpha_2) =\{ [\wtilde A,\wtilde \Phi]\in
 \cM_{\varpi_\tau}(M_\tau)
\quad | \quad \exists u\in
\cG(M_\tau),\\
\|u\cdot (\wtilde A,\wtilde \Phi)-(A,\Phi)\|_{L^2_2(g_\tau,A)} <\alpha_2\}
\end{multline*}
is an embedding.
\end{prop}
\begin{remark*}
The set $ B({[A,\Phi]},\alpha)$ is just a ball of center $[A,\Phi]$
and `$L^2_2$--radius' $\alpha_2$ in $\cM_l(M_\tau)$.
\end{remark*}
Before giving the proof of this proposition we show
how it implies
\fullref{theogluing}.
\begin{cor}
\label{cordiffeo}
For every $\tau$  and $N_0$ large enough, the gluing map
$$\gG\co \cM_{\varpi_\tau}(M_\tau)\to
\cM_{\varpi_\tau}(M'_\tau)$$
 is a
diffeomorphism.
\end{cor}
\begin{proof}
We already know that $\gG$ is  a local diffeomorphism, so we just need
to prove that it is $1:1$.
We show first that $\gG$ is globally injective for $\tau$ large
enough.  Suppose it is not true: then we have a sequence
$\tau_j\to\infty$ and solutions of Seiberg--Witten equations
$(A_j,\Phi_j)$ and  $(\wtilde A_j,\wtilde \Phi_j)$ on $M_{\tau_j}$ such that
$[A_j,\Phi_j]\neq [\wtilde A_j,\wtilde \Phi_j]$ and $\gG(A_j,\Phi_j)=
u'_j\cdot \gG
(\wtilde A_j,\wtilde \Phi_j)$ for some gauge transformations
$u'_j\in\cG_l(M'_{\tau_j})$.

After applying gauge transformations and extracting a subsequence, we
may assume that $(A_j,\Phi_j)$ and  $(\wtilde A_j,\wtilde \Phi_j)$ have
exponential decay and converge on every
compact to some solutions $(A,\Phi)$ and $(\wtilde A,\wtilde \Phi)$ of
Seiberg--Witten equations on $M$.
Notice that if $(A,\Phi)=(\wtilde A,\wtilde \Phi)$, it implies
$$ \|(A,\Phi)-(\wtilde A,\wtilde \Phi)\|_{L^2_2(g_\tau,A_j)}\to 0.
$$
We are going to see that it is indeed the case:
 $\|\gG(A_j,\Phi_j)-(A_j,\Phi_j)^\sharp\|_{L^2_2(g_\tau',
  A_j^\sharp)} \to 0$ according to \fullref{theogludef} and
  \fullref{lemmagoodapprox}. Hence $\gG (A_j,\Phi_j)$ converges on
  every compact toward
  $(A,\Phi)$ since $(A_j,\Phi_j)^\sharp$ does. Similarly $\gG
  (\wtilde A_j,\wtilde \Phi_j)$ converges to
  $(\wtilde A,\wtilde \Phi)$. The fact that $\gG
  (A_j,\Phi_j)$ and $\gG
  (\wtilde A_j,\wtilde \Phi_j)$ are gauge equivalent for each $j$ implies
  that the
  limits are also gauge equivalent.
After making further gauge transformations, we can assume that $
  (A_j,\Phi_j)$ and $
  (\wtilde A_j,\wtilde \Phi_j)$ converge toward the same limit $(A,\Phi)$
  on $M$.

Therefore we have for $j$ large
enough
$$ \|(\wtilde A_j,\wtilde \Phi_j)-
  (A_j,\Phi_j)\|_{L^2_2(g_\tau, A_j)} < \alpha,
$$
where $\alpha$ is chosen according to \fullref{proplocinj}. The fact
that
$
 [\gG (\wtilde A_j,\wtilde \Phi_j)] = [\gG (\wtilde A_j,\wtilde \Phi_j)]$
 should then imply $
 [\wtilde A_j,\wtilde \Phi_j] = [\wtilde A_j,\wtilde \Phi_j]$. This is a
 contradiction, hence $\gG$ must be injective.

We show now that $\gG$ is surjective: for $\tau$ large enough, there is
a second
gluing map
$$ \gG'\co \cM_{\varpi_\tau}(M'_\tau)\to \cM_{\varpi_\tau}(M_\tau)
$$
since $M'_\tau$ and $M_\tau$ play symmetric roles. The map $\gG'$
enjoys all the properties of $\gG$. In particular $\gG'\circ\gG$ is an
embedding of $\cM_{\varpi_\tau}(M_\tau)$ into itself. Using the fact that
$\cM_{\varpi_\tau}(M_\tau)$ is a finite dimensional compact manifold,
we conclude that
$\gG'\circ\gG$ is therefore a diffeomorphism. Therefore, $\cG$ must be
surjective, otherwise, we would have an $[A',\Phi']\not\in\Im
\gG$. On the other hand, using the fact that $\gG'\circ \gG$ is
surjective, there exists $[A,\Phi]\in\cM_{\varpi_\tau}(M\tau)$ such that
$\gG'\circ\gG[A,\Phi] = \gG'[A',\Phi']$ which contradicts
the injectivity of $\gG'$.
\end{proof}

\begin{remarks}
\label{rqisotopy}
\item An immediate consequence of \fullref{cordiffeo} is that
  $\sharp$ is an embedding for every $\tau$ and $N_0$ large enough.
\item More generally,  under the assumption of
  \fullref{theogludef} and \fullref{cordiffeo}, we can introduce the map
$$\gG_s (A,\Phi) = (A,\Phi)^\sharp + s
  \delta^*_{2,(A,\Phi)^\sharp}(b',\psi'),$$
for a parameter $s\in[0,1]$. Then, $\gG_s(\cM_{\varpi_\tau}(M_\tau))$
  realizes an
  isotopy between $\sharp(\cM_{\varpi_\tau}(M_\tau))$ and
  $\cM_{\varpi_\tau}(M'_\tau)$ in $\cC/\cG(M'_\tau)$.
\end{remarks}

We return now to the proof of \fullref{proplocinj}.
\begin{proof}
It is convenient to study locally $\gG$ in the charts provided by the
slice theorem. We begin by observing that we have charts of radius
uniformly bounded from below in \fullref{theoslice}.

Let $\alpha_1$, $\alpha_2\in (0,\alpha_1]$ and $U_{\alpha_1}$ be as in
\fullref{theoslice}, for a solution $(A,\Phi)$ of the Seiberg--Witten
equations on $M_\tau$. Put $(A',\Phi'):= \gG(A,\Phi)$ and
denote	$\alpha'_1$, $\alpha'_2$ and
$U'_{\alpha'_1}$ the analogous data on $M'_\tau$.

Let  $(\wtilde A,\wtilde \Phi)$ be a solution of Seiberg--Witten
equations such that
$$[\wtilde A,\wtilde \Phi]\in
B([A,\Phi],\alpha_2).$$
 Then, we may assume, thanks to the slice
theorem,  that up to a gauge transformation, we have $(\wtilde A,\wtilde
\Phi)\in U_{\alpha_1}$.

We are going to show that if $\alpha_1$ is chosen small enough,
independently of $(A,\Phi)$, $N_0\geq 1$ and $\tau$ large enough,
then we automatically have $[\gG(\wtilde A, \wtilde \Phi)]\in
B([\gG(A,\Phi)],\alpha'_2)$.

The next lemma is a classical application of elliptic regularity for
Seiberg--Witten equations and \fullref{prop:contrd1}. Again, the constant
involved can be
chosen uniformly thanks to the fact that it is the case for Sobolev
constants due to \fullref{lemmaboundedgeometry}, and that the moduli
spaces are uniformly bounded in the sense
of \fullref{boundedderivatives}.
\begin{lemma}
\label{lemmaclassical}
There exist  constants $C,\alpha_1>0$ and a compact set $K\subset M$
large enough such that for every every $\tau$ large
enough and every solutions $(A,\Phi)$
and $(\wtilde A,\wtilde \Phi)$ of Seiberg--Witten equations on $M_\tau$
with
$$\delta_{1,(A,\Phi)}^* \left( (\wtilde A,\wtilde \Phi)-(A,\Phi)
\right	 )= 0 \quad \mbox{ and }\quad \|(\wtilde A,\wtilde
\Phi)-(A,\Phi)\|_{L^2_2(g_\tau,A)}\leq \alpha_1,$$
 we have
 $$\|(\wtilde A,\wtilde
\Phi)-(A,\Phi)\|_{L^2_3(g_\tau,A)}\leq C\|(\wtilde A,\wtilde
\Phi)-(A,\Phi)\|_{L^2(g_\tau,K)}.
$$
\end{lemma}

The solutions provided by the gluing map have the form
\begin{align*}
(A',\Phi'):=\gG(A, \Phi) & = ( A, \Phi)^\sharp +\delta^*_{2,(A,
  , \Phi)^\sharp } (  b',\psi')\\
\gG(\wtilde A,\wtilde \Phi) & = (\wtilde A,\wtilde \Phi)^\sharp
  +\delta^*_{2,(\wtilde
  A,\wtilde \Phi)^\sharp} ( \wtilde b', \wtilde\psi')
\end{align*}
for some $(\psi, b)$ and $(\wtilde \psi',\wtilde b')$  given by
\fullref{theogludef}.  Moreover we have
\begin{align} \label{estinj0}
  \|\delta_{2,(A,\Phi)^\sharp}^*(b',\psi')\|_{L^2_2(g_\tau,A^\sharp)}
  & \leq
c\|\SW(A,\Phi)^\sharp\|_{L^2_1(g'_\tau,A^\sharp)}, \\ \label{estinj2}
\|(\wtilde
  b',\wtilde \psi')\|_{L^2_3(g_\tau, \wtilde A^\sharp)},\quad
  \|\delta_{2,(\wtilde A,\wtilde \Phi)^\sharp}^*(\wtilde
  b',\wtilde \psi')\|_{L^2_2(g_\tau, \wtilde A^\sharp)}
  & \leq
c\|\SW(\wtilde A,\wtilde \Phi)^\sharp\|_{L^2_1(g'_\tau,\wtilde A^\sharp)} ,
\end{align}
and we may assume that the above quantities are uniformly small for
$\tau$ large enough according to \fullref{lemmagoodapprox}.

Let $k_{(\wtilde A,\wtilde \Phi)\co }
W_{(A,\Phi)}\to W_{(\wtilde A,\wtilde \Phi)}$
be the isomorphism defined in
\fullref{secrough} and put
$$(\what A,\what \Phi)= k_{(\wtilde A,\wtilde \Phi)}^{-1}\cdot
(\wtilde A,\wtilde \Phi)^\sharp$$
 so that $(A,\Phi)^\sharp$ and $(\what
A,\what \Phi)$ are now defined for the same spin bundle $W_{(A,\Phi)}$.
Then
$$k_{(\wtilde A,\wtilde \Phi)}^{-1} \gG(\wtilde A,\wtilde \Phi) = (\what
  A,\what \Phi ) +\delta^*_{2,(\what A
  ,\what \Phi)} ( \what b',\what \psi'),$$
where $(\what b',\what\psi'):=k_{(\wtilde A,\wtilde
  \Phi)}^{-1}(\wtilde b',\wtilde\psi')$. Notice that with this notation,
  $k_{(\wtilde A,\wtilde
  \Phi)}$ is the identity on the bundle of self-dual forms.
By gauge invariance
\begin{align}
\label{estinj3} \|\delta_{2,(\wtilde A,\wtilde \Phi)^\sharp}^*(\wtilde
  b',\wtilde \psi')\|_{L^2_2(g_\tau, \wtilde A^\sharp)} &=
 \|\delta_{2,(\what A,\what \Phi)}^*(\what
  b',\what \psi')\|_{L^2_2(g_\tau, \what A)},\\
\|(\wtilde
  b',\wtilde \psi')\|_{L^2_3(g_\tau, \wtilde A^\sharp)} & = \|(\what
  b',\what  \psi')\|_{L^2_3(g_\tau, \what A)}.
\end{align}
 Put
$(a',  \phi')= (\what
A,\what \Phi)-(A,\Phi)^\sharp$; this is just the variation of $\sharp$
with our
rough gauge choice. If $\alpha_1$ is small enough, we get an
  $L^2_3$--estimate on $(\wtilde A,\wtilde \Phi)-(A,\Phi)$ by
  \fullref{lemmaclassical}, hence we can
  assume that \fullref{gluingestimate0} applies and we have
  eventually an
  estimate
\begin{equation}
\label{estinj1}
\|(a',	\phi')\|_{L^2_3(g'_\tau,A^\sharp)} \leq cC\|(\wtilde
  A,\wtilde \Phi)-(A,\Phi)\|_{L^2_2(g_\tau,A)} \leq cC\alpha_1.
\end{equation}
Therefore, the Sobolev multiplication theorems show that for a choice
of $\alpha_1$ small enough the $L^2_k(g'_\tau,A^\sharp)$ and
$L^2_k(g'_\tau,\what A)$ are commensurate up to $k=3$, so that, say
\begin{align} \label{estinj4}
\|\delta_{2,(\what A,\what \Phi)}^*(\what
  b',\what \psi')\|_{L^2_2(g_\tau, A^\sharp )} &\leq 2 \|\delta_{2,(\what
  A,\what \Phi)}^*(\what
  b',\what \psi')\|_{L^2_2(g_\tau, \what A)}\\ \label{estinj5}
\|(\what
  b',\what \psi')\|_{L^2_3(g_\tau, A^\sharp )} &\leq 2 \|(\what
  b',\what \psi')\|_{L^2_3(g_\tau, \what A)}.
\end{align}
Eventually, $\alpha_1$ controls the $L^2_2(g_\tau,A^\sharp)$--norm of
$( a',\phi')$ by~\eqref{estinj1}. The  \\ $L^2_2(g_\tau,A^\sharp)$ norm
of $\delta_{2,(\what A,\what \Phi)}^*(\what
  b',\what \psi')$ is controlled by
    $\|\SW(\wtilde A,\wtilde \Phi)^\sharp\|_{L^2_1(g_\tau,\wtilde A^\sharp
  )}$ thanks to~\eqref{estinj2} \eqref{estinj3} and \eqref{estinj4},
  while the one of $\delta_{2,(
    A, \Phi)^\sharp }^*(
  b', \psi')$ is controlled  via \eqref{estinj0} by
    $\|\SW(A,\Phi)^\sharp \|_{L^2_1(g_\tau, A^\sharp
  )}$.

Similarly, the $L^2_2(g_\tau, A^\sharp)$ and $L^2_2(g_\tau, A')$ are
also commensurate, and if $\tau$ is large enough, we may replace
$A^\sharp$ by $A'$ in the estimates~\eqref{estinj0} and
\eqref{estinj4}. That means that for a suitable choice of $\alpha_1$,
we will have automatically  $[\gG(\wtilde A,
  \wtilde\Phi)]\in B([\gG(A,\Phi)],\alpha'_2)$.

Therefore, using the slice theorem about $[A',\Phi']= [\gG(A,\Phi)]$,
we can recast $\gG$
into a map
$$ \gG_1\co  \cU\to U'_{\alpha'_1},
$$
where $\cU$ is a finite dimensional submanifold of $U_{\alpha_1}$
corresponding to Seiberg--Witten moduli space
$\cM_{\varpi_\tau}(M_\tau)$ in
this local
neighborhood of $[A,\Phi]$.

Everything is in order to study  the variations of $\gG$ about
$[A,\Phi]$, that is to say the variations of $\gG_1$ about  $0\in \cU$.
\begin{lemma}
\label{lemmainjest}
For every $N_0$ large enough and $\epsilon>0$, there exists $\alpha_1>0$
such that for every
$\tau$ large enough and
every pair $(A,\Phi)$ and $(\wtilde A,\wtilde \Phi)$ of solutions of
Seiberg--Witten
equations   with $(\wtilde A,\wtilde \Phi)$ in the slice about
$(A,\Phi)$ and with
$\|(\wtilde A,\wtilde \Phi)-(A,\Phi)\|_{L^2_2(g_\tau,A)}\leq \alpha_1$,
\begin{multline*}
 \|\Delta_{2,(A,\Phi)^\sharp}  \left (( \what b',\what\psi') -  ( b',\psi')
 \right
)\|_{L^2_1(g'_\tau,A^\sharp) }\\
\leq
\epsilon \left (  \|(\wtilde A,\wtilde \Phi)-(A,\Phi)\|_{L^2_2(g_\tau,A)} +
\| ( \what b',\what\psi') -  ( b',\psi')
\|_{L^2_3 (g'_\tau,A^\sharp )} \right )
\end{multline*}
where
 $(b',\what\psi')$ and $(b',\psi')$ are defined as above.
\end{lemma}
\begin{proof}
Notice first that we have
$$ \delta^*_{2,(
\what  A, \what \Phi)^\sharp} ( \what b',\what\psi') =\delta^*_{2,(
  A, \Phi)^\sharp} ( \what b',\what\psi') +
  \Bigl \{(a',\phi' ),( \what b',\what \psi' ) \Bigr \},
$$
where $\{\cdot,\cdot\}$ is a bilinear pairing with fixed coefficients,
whose particular
expression is irrelevant to us. More generally, we will denote any
bilinear pairing in this way in the proof of the lemma.

The Seiberg--Witten equations for  $k^{-1}_{(\wtilde
  A,\wtilde\Phi)}\gG(\wtilde A,\wtilde \Phi)$ give  the identity
\begin{align*}
  &\delta_{2,(A,\Phi)^\sharp} (a', \phi') + \Delta_{2,(
  A, \Phi)^\sharp} ( \what b',\what\psi' ) +
  \delta_{2,(A,\Phi)^\sharp} \left\{ (a',\phi' ),(
  \what b',\what\psi') \right\}  \\ & + \cQ\Bigl( (a', \phi')
+\delta^*_{2,(
  A, \Phi)^\sharp} ( \what b',\what\psi')+ \left\{ (a',\phi' ),(
  \what b',\what\psi') \right\}
\Bigr)	= -\SW(A,\Phi)^\sharp,
\end{align*}
where $\cQ$ is the quadratic term of the Seiberg--Witten defined
at~\eqref{eq:quaddef2}.
Taking the difference with~\eqref{eqglue}, we have
\begin{align*}
    \Delta_{2,(
      A, \Phi)^\sharp}	& \left ((  \psi', b') -  ( \what b',\what \psi')
      \right ) =
       \delta_{2,(A,\Phi)^\sharp}
(\dot
 a, \dot\phi)	+   \delta_{2,(A,\Phi)^\sharp}  \left\{ (a',\phi' ),(
  \what b',\what\psi') \right\}  \\
&\hspace{-.3in}-  \cQ\left ( \delta_{2,(A,\Phi)^\sharp}^* ( b',\psi') \right ) 
+\cQ\left ( (a', \phi')
  +\delta_{2,(A,\Phi)^\sharp}^*
 ( \what b', \what  \psi')+   \left\{ (a',\phi' ),(
  \what b',\what\psi') \right\} \right )
\end{align*}
Developing the last line, we find an expression which is formally
\begin{align}
\label{l1}    \Delta_{2,(
      A, \Phi)^\sharp}	 \left ((  \psi', b') -  (\what	b',\what \psi')
      \right ) =
      & \delta_{2,(A,\Phi)^\sharp} (
(a',\phi')    \\
\label{l2}
+&  \cQ\left ( \delta_{2,(A,\Phi)^\sharp}^* ( \what b',\what \psi') \right )
  - \cQ\left ( \delta_{2,(A,\Phi)^\sharp}^* ( b',\psi') \right ) \\
\label{l3}
+& \cQ((a', \phi' )) \\
\label{l4}
+& \cQ\left ( \left\{ (a',\phi' ),(
  \what b',\what\psi') \right\}\right ) \\
\label{l5}
+&   \delta_{2,(A,\Phi)^\sharp}  \left\{ (a',\phi' ),(
  \what b',\what\psi') \right\} \\
\label{l6}
+ & \left \{ (a',\phi'	), \left\{ (a',\phi' ),(
  \what b',\what\psi') \right\} \right \}\\
\label{l7}
+ & \left \{  \delta_{2,(A,\Phi)^\sharp}^* ( \what b',\what \psi') , \left\{
(a',\phi' ),(
  \what b',\what\psi') \right\} \right \} \\
\label{l8}
+ & \left \{  \delta_{2,(A,\Phi)^\sharp}^* ( \what b',\what \psi') ,
      (a',\phi' )\right \}
\end{align}
We explain  how to control each term of the RHS: we know that the
$L^2_3(g_\tau, A^\sharp)$--norm of $(a',\phi' )$ is controlled
by its $\alpha_1$ (see~\eqref{estinj1}). Hence we may assume that
\fullref{gluingestimate} applies. Therefore, for $N_0\geq 1$ large enough,
we have
$$\|\delta_{2,(A,\Phi)^\sharp}(a',\phi')\|_{L^2_1(g'_\tau,A^\sharp)}\leq
\frac \epsilon 7 \|(\wtilde A,\wtilde
\Phi) - (A,\Phi)\|_{L^2_2(g_\tau, A)}.$$
The term of line~\eqref{l2} is controlled using the fact that $\cQ$ is
a quadratic polynomial. Similarly to~\eqref{estcs}, we have an
estimate
\begin{multline*} \left \|\cQ\left ( \delta_{2,(A,\Phi)^\sharp}^*
( \what b',\what
\psi') \right )
  - \cQ\left ( \delta_{2,(A,\Phi)^\sharp}^* ( b',\psi') \right )
  \right \|_{L^2_1(g'_\tau,A^\sharp)}\\
\leq c \|(b',\psi')-(\what b',\what \psi')\|_{L^2_3(g'_\tau,A^\sharp)}
\|(b',\psi')+ (\what b',\what \psi')\|_{L^2_3(g'_\tau,A^\sharp)}.
\end{multline*}
Since the $L^2_3(g'_\tau,A^\sharp)$--norm of $(b',\psi')$ and $(\what
b',\what \psi')$ is arbitrarily
small for $\alpha_1$ small enough and $\tau$ large enough, we deduce
that for a suitable choice of $\alpha_1$, we will have
\begin{multline*} \left \|\cQ\left ( \delta_{2,(A,\Phi)^\sharp}^*
( \what b',\what
\psi') \right )
  - \cQ\left ( \delta_{2,(A,\Phi)^\sharp}^* ( b',\psi') \right )
  \right \|_{L^2_1(g'_\tau,A^\sharp)}\\
\leq \epsilon	\|(b',\psi')-(\what b',\what
\psi')\|_{L^2_3(g'_\tau,A^\sharp)},
\end{multline*}
for every $\tau$ large enough.

A similar technique together with \fullref{gluingestimate0}
and~\eqref{estinj1}  shows
that the term at line~\eqref{l3} is controlled
under the same circumstances by
\begin{equation*}
 \|\cQ\left ( (a',\phi')
 \right )  \|_{L^2_1(g_\tau',A^\sharp}
\leq \frac \epsilon 7 \|(\wtilde A,\wtilde
\Phi) - (A,\Phi)\|_{L^2_2(g_\tau, A)}.
\end{equation*}
For the term at line~\eqref{l4}, we use the fact that
it is this time a homogeneous polynomial expression of degree $4$ in
$(a',\psi')$
and $(\what b',\what\psi')$. Thanks to the Sobolev embedding
theorem $L^2_2 \hookrightarrow L^8$, we deduce that
$$ \|\cQ\left ( \left\{ (a',\phi' ),(
  \what b',\what\psi') \right\}\right )\|_{L^2_1(g'_\tau, A^\sharp)} \leq c
  \|(a',\phi')\|_{L^2_3(g'_\tau, A^\sharp)}^2
  \|(b',\psi')\|_{L^2_3(g'_\tau, A^\sharp)}^2 .
$$
The estimate~\eqref{estinj0} and the fact that the
$L^2_3(g'_\tau,A^\sharp)$--norm of $(\what b',\what\psi')$ is arbitrarily
small for $\alpha_1$ small enough and $\tau$ large lead to an estimate
\begin{equation*}
 \|\cQ\left (  \left\{ (a',\phi' ),(
  \what b',\what\psi') \right\}
 \right )  \|_{L^2_1(g_\tau',A^\sharp)}
\leq \frac \epsilon 7  \|(\wtilde A,\wtilde
\Phi)-(A,\Phi)\|_{L^2_2(g_\tau,A)}.
\end{equation*}
Arguing in the same manner, it is easy to show that the
$L^2_1(g'_\tau,A^\sharp)$--\-norms of~\eqref{l5}, \eqref{l6}, \eqref{l7}
and \eqref{l8} are lower bounds for $$\frac \epsilon 7	\|(\wtilde
A,\wtilde \Phi)-(A,\Phi)\|_{L^2_2(g_\tau,A)}$$
 if $\alpha_1$ is chosen
small enough and $\tau$ large enough. Summing up all the estimates, we
obtain the lemma.
\end{proof}

We return to the proof of \fullref{proplocinj}
Applying \fullref{propfirstev}, we deduce that for some
universal constant $c_1>0$, we have
$$  c_1 \|( b',\what\psi') -  ( b',\psi')\|_{L^2_3(g'_\tau,A^\sharp)} \leq
 \|\Delta_{2,(A,\Phi)^\sharp}  \left (( b',\what\psi') -  ( b',\psi')
 \right )\|_{L^2_1(g'_\tau,A^\sharp)}.
$$
Together with the \fullref{lemmainjest}  this implies
\begin{equation}\label{eqinjest6}  (c_1-\epsilon) \|( \what \psi', \what
b') -  ( b',\psi')\|_{L^2_3(g'_\tau,A^\sharp)} \leq \epsilon
\|(\wtilde A, \wtilde \Phi)-(A,\Phi)\|_{L^2_2(g_\tau,A)} .
\end{equation}
Using the estimate~\eqref{eqbound3}, that we have
$$\|\delta^*_{2,(A,\Phi)^\sharp}\Bigl ( (\what b',\what \phi')-
(b',\psi')\Bigr )\|_{L^2_2(g'_\tau, A^\sharp)}\leq
c_2 \| (\what b',\what \phi') - (b',\psi')\|_{L^2_3(g'_\tau, A^\sharp)},
$$
and it follows, once we made sure that we started with $\epsilon <
c_1$,  that
$$\|\delta^*_{2,(A,\Phi)^\sharp}\Bigl ( (\what b',\what \phi')-
(b',\psi')\Bigr )\|_{L^2_2(g'_\tau, A^\sharp)}\leq
\frac {c_2\epsilon}{c_1 - \epsilon }\|(\wtilde A, \wtilde
\Phi)-(A,\Phi)\|_{L^2_2(g_\tau,A)} .
$$
Now
\begin{multline}
k^{-1}_{(\wtilde A,\wtilde \Phi)}\gG(\wtilde A,\wtilde
\Phi)- \gG(A,\Phi) \\
= (a',\phi' ) + \delta^*_{2,(A,\Phi)^\sharp}
\left (( \what b',\what\psi') -  ( b',\psi') \right ) + \left
\{(a',\phi'),(\what b', \what \phi')\right \},
\end{multline}
and the last term is controlled as in the proof of \fullref{lemmainjest},
i.e.\ for $\tau$ large enough, we have say
$$ \|\left
\{(a',\phi'),(\what b', \what \phi')\right \} \|_{L^2_2(g'_\tau,
  A^\sharp)} \leq \epsilon \|(\wtilde A,\wtilde\Phi)-
(A,\Phi)\|_{L^2_2(g_\tau, A)}.
$$
Hence we have the estimate
\begin{multline}  \|(
a',\phi')\|_{L^2_2(g_{\tau}',A^{\sharp})} -
\left (\epsilon +\frac {c_2\epsilon}{c_1 - \epsilon }\right ) \|(\wtilde A,
\wtilde \Phi)-(A,\Phi)\|_{L^2_2(g_\tau,A)} \\
\leq\|k^{-1}_{(\wtilde A,\wtilde \Phi)}\gG(\wtilde A,\wtilde \Phi)-
\gG(A,\Phi)\|_{L^2_2(g_{\tau}',A^{\sharp})} .
\end{multline}
Then, we can apply \fullref{lemmaclassical}, and since the identity
$(a',\phi') =  (\wtilde A,\wtilde \Phi)-(A,\Phi)$ holds on every compact
set for $\tau$ large enough, we have
\begin{multline}
\label{eqimm}
 \left ( \frac 1C-\epsilon- \frac {c_2\epsilon}{c_1 - \epsilon }\right )
\|(\wtilde A, \wtilde \Phi)-(A,\Phi)\|_{L^2_2(g_\tau,A)} \\
\leq\|k^{-1}_{(\wtilde A,\wtilde \Phi)}\gG(\wtilde A,\wtilde \Phi)-
\gG(A,\Phi)\|_{L^2_2(g_{\tau}',A^{\sharp})} .
\end{multline}
If we take $\epsilon$ small enough in the first place, so that $
0<\epsilon+  \frac
{c_2\epsilon}{c_1 - \epsilon } < \frac 1{2C}$, we see immediately
that $\gG_1$ is injective.

However, a little more work is needed to see that $\gG_1$ is an
immersion at the origin.
We estimate first how far is $ k^{-1}_{(\wtilde A,\wtilde \Phi)} \gG(\wtilde
A,\wtilde \Phi)$ from
being in a Coulomb gauge w.r.t.
$ \gG(A,\Phi)$ in the next lemma.
\begin{lemma}
\label{lemmaimmest}
For every $N_0$ large enough and $\epsilon>0$, there exists $\alpha_1>0$
such that for every
$\tau$ large enough and
every pair $(A,\Phi)$ and $(\wtilde A,\wtilde \Phi)$ of solutions of
Seiberg--Witten
equations on $M_\tau$  with $(\wtilde A,\wtilde \Phi)$ in the slice about
$(A,\Phi)$ and with
$\|(\wtilde A,\wtilde \Phi)-(A,\Phi)\|_{L^2_2(g_\tau, A)}\leq \alpha_1$,
\begin{equation}
 \|\delta_{1,(A',\Phi')}^*  \left (k^{-1}_{(\wtilde A,\wtilde \Phi)}
\gG(\wtilde A,\wtilde \Phi)-(A',\Phi')
 \right
)\|_{L^2_1(g'_\tau) }
\leq
\epsilon   \|(\wtilde A,\wtilde \Phi)-(A,\Phi)\|_{L^2_2(g_\tau,A)}
\end{equation}
where
 $(A',\Phi'):=\gG(A,\Phi)$.
\end{lemma}
\begin{proof}
We have
\begin{align}
\label{eqimmest1}
\delta_{1,(A,\Phi)^\sharp}^*\left ( k^{-1}_{(\wtilde A,\wtilde
  \Phi)} \gG(\wtilde A,\wtilde
\Phi)- (A',\Phi') \right )  & = \delta_{1,(A,\Phi)^\sharp}^* (
  a', \phi ')  \\
\label{eqimmest2}
& + \delta_{1,(A,\Phi)^\sharp}^* \delta^*_{2,(A,\Phi)^\sharp}
\left (( \what b',\what\psi') -  ( b',\psi') \right ) \\
\label{eqimmest3}
& +\delta_{1,(A,\Phi)^\sharp}^* \left
\{(a',\phi'),(\what b', \what \phi')\right \}.
\end{align}
The $L^2_1(g'_\tau)$ norm of the RHS of~\eqref{eqimmest1} is controlled
via \fullref{lemmainjest}.
The operator
$\cL_{(A,\Phi)^\sharp}=\delta_{2,(A,\Phi)^{\sharp}}
\delta_{1,(A,\Phi)^{\sharp}}$
has norm very close to $0$ by \fullref{lemmagoodapprox}
and~\eqref{eqdefect}, therefore, the $L^2_1(g'_\tau)$--norm of the term
in~\eqref{eqimmest2} is controlled using~\eqref{eqinjest6}.
The term~\eqref{eqimmest3} is controlled using Sobolev multiplication
as in \fullref{lemmainjest}.
\end{proof}

An immediate consequence of \fullref{lemmaimmest} concerns the
 problem of fixing a Coulomb gauge:
\begin{cor}
\label{corcoulomb}
For every $N_0$ large enough and $\epsilon>0$, there exists $\alpha_1>0$
such that for every
$\tau$ large enough and
every pair $(A,\Phi)$ and $(\wtilde A,\wtilde \Phi)$ of solutions of
Seiberg--Witten
equations on $M_\tau$  with $(\wtilde A,\wtilde \Phi)$ in the slice about
$(A,\Phi)$ and with
$\|(\wtilde A,\wtilde \Phi)-(A,\Phi)\|_{L^2_2(g_\tau, A)}\leq \alpha_1$,
there exists a gauge transformation $u'\in \cG_2(M_\tau')$ such that
\begin{align*}
&\delta_{1,(A',\Phi')}^* \left (u'\cdot k^{-1}_{(\wtilde A,\wtilde\Phi)}
\gG(\wtilde A,\wtilde\Phi) -(A',\Phi')\right )  = 0 \\
& \|1-u'\|_{L^2_3(g'_\tau)}\leq  \epsilon \|(\wtilde A,\wtilde
\Phi)-(A,\Phi)\|_{L^2_2(g_\tau, A)} .
\end{align*}
\end{cor}
\begin{proof}
Using \fullref{lemmaimmest}, we can show that the
$L^2_{2,3}(g_\tau,A')$--norm
of $$d_0\nu^{-1}\left (   k^{-1}_{(\wtilde A,\wtilde\Phi)}
\gG(\wtilde A,\wtilde\Phi) -(A',\Phi') \right )$$
 is controlled by
$\epsilon\|(\wtilde A,\wtilde
\Phi)-(A,\Phi)\|_{L^2_2(g_\tau, A)}$ (see \fullref{refslice}). 

Then
\fullref{corF} implies that $u'=e^{v'}$, with the
$L^2_3(g'_\tau)$--norm of $v'$ controlled by $\epsilon \|(\wtilde A,\wtilde
\Phi)-(A,\Phi)\|_{L^2_2(g_\tau, A)} $. The
corollary follows using \fullref{corexp}.
\end{proof}

We can finish the proof of \fullref{proplocinj}.
Applying \fullref{corcoulomb} with
 a choice of $\epsilon$ small enough, we deduce from~\eqref{eqimm}
 that for
 some constant $c>0$, we have
$$ c\|(\wtilde A,\wtilde \Phi)-(A,\Phi)\|_{L^2_2(g_\tau,A)} \leq	\|u'\cdot
k^{-1}_{(\wtilde A,\wtilde\Phi)}
\gG(\wtilde A,\wtilde\Phi) -(A',\Phi')\|_{L^2_2(g_\tau,A)}.
$$
Since $u'\cdot k^{-1}_{(\wtilde A,\wtilde\Phi)}
\gG(\wtilde A,\wtilde\Phi) -(A',\Phi') =\gG_1\left ( (\wtilde A,\wtilde
\Phi)-(A,\Phi \right )$, the above inequality shows that $\gG_1$ is an
immersion at the origin.
\end{proof}

\subsection{Orientations}
In this section, we show that with suitable conventions, the gluing map
$\gG\co  \cM_{\varpi_\tau}(M_\tau)$ $\to \cM_{\varpi_\tau}(M'_\tau)$, is
orientation preserving.

Suppose that we are given an orientation of
$\cM_{\varpi_\tau}(M_\tau)$ that is to say an orientation of the index
bundle $\Omega(\cD)$ where $\cD_{(A,\Phi)}$ is the operator defined
at~\eqref{eqlinop}. The aim of the next section is to explain how to
deduce an orientation of $\cM_{\varpi_\tau}(M'_\tau)$ from the
orientation of $\cM_{\varpi_\tau}(M_\tau)$.

\subsubsection{Excision principle for the index}
Let $\overN$ be an almost complex manifold with boundary $\del
\overN = \del \overM =Y$. Such a manifold always exists (see
for instance \cite[Lemma 4.4]{Go}).
From this data, we can construct $N_\tau$ and $N'_\tau$ by gluing the
AFAK ends $Z$ and $Z'$ similarly to $M_\tau$ and $M'_\tau$. Moreover
$N_\tau$ and $N'_\tau$ are endowed with their canonical
\spinc--structures of almost complex type and suitable Hermitian
metrics $\wtilde g_\tau$ and $\wtilde g'_\tau$. Remark that the canonical
solutions $(B,\Psi)$ and $(B',\Psi')$ have now a natural extension
to the compact parts of
 $N_\tau$
and  $N'_\tau$. The identity $\Psi=(1,0)\in \Lambda ^{0,0}\oplus
\Lambda^{0,2}$	and the spin connection $B$ deduced from the Hermitian
metric and the Chern connection on the determinant line bundle of the
\spinc--structure  now make sense globally
on $N_\tau$. Similarly, $(B',\Psi')$ is defined globally on $N'_\tau$.

The almost complex
structure and
the Hermitian metric $\wtilde g_\tau$ on
$N_\tau$ induce a nondegenerate $2$--form $\omega_\tau$. However, this
$2$--form is not closed in general. It has  torsion (or
a Lee form)
$\theta$ defined by
$$ d\omega_\tau = \theta\wedge \omega_\tau.
$$
Notice that by construction of $\wtilde g_\tau$,  $d\omega_\tau$  is
identified with a compactly supported form
 on $N$, independent of $\tau$,  and so can $\theta$.

According to a computation of Taubes~\cite{Taub} and later
Gauduchon~\cite{G}, we have in this case
$$\Dirac_B \Phi = \Dirac^{can} \Phi + \frac 14 \theta \cdot \Phi, $$
where $\theta$ acts by Clifford product on spinors and $\Dirac
^{can}=\sqrt 2(\delb +\delb^*)$.
Therefore, $(B,\Psi)$ is solution of the perturbed Seiberg--Witten
equations
\begin{eqnarray}
\label{eqSW3}
F_A^+- \{\Phi\otimes\Phi^*\}_0 & = & F_B^+- \{\Psi\otimes\Psi^*\}_0 \\
\Dirac_A \Phi & = & \frac 14 \theta\cdot \Phi .
\end{eqnarray}
The linearization $\wtilde \cD_{(B,\Psi)}$ of the above equations differs
from the linearization	$\cD_{(B,\Psi)}$ of the standard one by a
$0$--order term. Precisely
\begin{equation}
\label{eqmodd}
 \wtilde \cD_{(B,\Psi)} (a,\phi) =  \cD_{(B,\Psi)}(a,\phi) -
(0,0,\frac 1 4 \theta\cdot \phi).
\end{equation}
According to the computation \cite[Proposition 3.8]{KM}, we see that
for every $\lambda>0$ large enough, the operator $ \wtilde
\cD_{(B,\lambda \Psi)}$ is injective (equivalently we could just scale
the metric). Replacing	$\Psi=(1,0)$ by $\lambda \Psi
=(\lambda,0)$ for such a large $\lambda$  does not affect in
anything to the theory developed
so far since the normalization for $\Psi$ was purely arbitrary.
However, for simplicity of notation, we will assume that $\lambda =1$
is a suitable value in the sequel.

 The index of  $ \wtilde \cD_{(B,\Psi)}$ equals the index
of $  \cD_{(B,\Psi)}$ for the two operator are homotopic
according to~\eqref{eqmodd}. But $\index  d:= \index \wtilde
\cD_{(B,\Psi)}$
is by definition the virtual dimension of the moduli space for a
\spinc--structure of almost complex type. This is zero, therefore $
\wtilde \cD_{(B,\Psi)}$ is an isomorphism. Hence the index bundle
$ \Omega (\wtilde \cD )$ is canonically trivial at $(B,\Psi)$ and we
deduce a compatible orientation of the index bundle.

The same remarks apply to $N'_\tau$ where $(B',\Psi')$ is solution of
some modified Seiberg--Witten equation similar to \eqref{eqSW3} whose
linearization  $\wtilde \cD_{(B',\Psi')}$ may be assumed to be an
isomorphism.

We are now going to construct a generalization of the pregluing map
$\sharp$ of \fullref{secpreg}. Let
$U\subset \cC_l(M_\tau)\times \cC_l(N'_\tau)$ be the gauge invariant  open
set that consists into pairs of configurations
$((A,\Phi), (\wtilde A',\wtilde \Phi'))$, such that  $|\beta|, |\wtilde
\beta'|>0$ (with
 $\Phi=(\beta,\gamma) $ and $\wtilde \Phi=(\wtilde \beta,\wtilde \gamma)
$ ) on the
 end $\{\sigma_\tau\geq T\}\cap \{
\wtilde \sigma'_\tau\geq T\}$.
Because of the uniform exponential decay of Seiberg--Witten solutions,
we may choose $T$ large enough so that for every  solutions of
Seiberg--Witten equations $(A,\Phi)$ on $M_\tau$, we have
$((A,\Phi),(B',\Psi'))\in U$.

Remark that the cut-off function $\chi_\tau$ defined
at~\eqref{eq:cutoff} makes sense over $M_\tau$, $M'_\tau$, $N_\tau$
and $N'_\tau$. Then we define a positive function $\mu_\tau$ on these
manifolds by the condition that
$$ \chi_\tau^2 +\mu_\tau^2 = 1.
$$
For any pair $((A,\Phi), (\wtilde A',\wtilde \Phi'))\in U$ we may assume
that,  after making a suitable gauge transformation,  $(A,\Phi)$ and
$(\wtilde A',\wtilde  \Phi' )$ are in real gauge on  the end. Then we
can define
a configuration on $M'_\tau$ by
$$ ( A', \Phi') = (  B'+\chi_\tau a+ \mu_\tau
\wtilde a',( \beta^ {\chi_\tau} (\wtilde \beta')^{\mu_\tau}, \chi_\tau
\gamma +
\mu_\tau \wtilde \gamma'))
$$
and on $N_\tau$ by
$$ (\wtilde A,\wtilde \Phi) = (B- \mu_\tau
a + \chi_\tau \wtilde a',( \beta^{- \mu_\tau} (\wtilde \beta')^
{\chi_\tau} , -\mu_\tau  \gamma+ \chi_\tau \wtilde
\gamma' ))
$$
where $A=B+a$ and $\wtilde A'=B'+\wtilde a'$ on the end. The configurations
$(\wtilde A,\wtilde
\Phi)$ and $( A',
\Phi')$ extend in a natural
way to	$M'_\tau$ and $N_\tau$.

Thus we have  defined a map
$$ f\co  \cC/\cG(M_\tau)\times \cC/\cG(N'_\tau)\supset
U/(\cG(M_\tau)\times \cG(N'_\tau))
\to \cC/\cG(M'_\tau)\times \cC/\cG(N_\tau)$$
by $f([A,\Phi],[\wtilde A',\wtilde \Phi'])= ([A',\Phi'],[\wtilde A,\wtilde
\Phi])$. In order to study the variations of $f$ we need to make a
gauge choice. The definition of a rough gauge $k$ as in
\fullref{secrough} can be extended trivially to this setting  and we
can study the variations of $k^{-1}\circ f$ around
an element of $U$ as we did for the map $\sharp$. A computation
similar to~\eqref{eqdiffrough} shows that the variation of $f$ about
$((A,\Phi), (\wtilde A',\wtilde \Phi')$ with a rough gauge fixing  is
given by the linear map
$$ F \co  T_{(A,\Phi)}\cC(M_\tau)\oplus T_{(\wtilde A',\wtilde
  \Phi')}\cC(N'_\tau) \longrightarrow T_{(A',\Phi')}\cC(M'_\tau)\oplus
  T_{(\wtilde A,\wtilde \Phi)}\cC(N_\tau)
$$
of the form
$$ F =
F_0  + d\chi_\tau F_1 + d \mu_\tau F_2,
$$
where $F_1$ and  $F_2$ are two matrices with constant coefficients whose
expression is here irrelevant to us and
$$F_0 =\left ( \begin{array}{cc}
\chi_\tau & \mu_\tau \\
-\mu_\tau & \chi_\tau
\end{array}
\right ).
$$
This latter matrix is invertible. Hence for a choice of $N_0$ large
enough, we have $d\chi_\tau$ and $d\mu_\tau$ get arbitrarily
small. Hence we can assume that
$F$ is invertible and is  homotopic to $F_0$ through isomorphisms.

  Let $P$ and $R$ be the operators
$$  P = \cD_{(A,\Phi)} \oplus \wtilde \cD_{(\wtilde A',\wtilde \Phi')},\quad
R= \cD_{(A',\Phi')} \oplus \wtilde \cD_{(\wtilde A,\wtilde \Phi)}.
$$
Since we assume that we are given an orientation of the index bundle
of $\cD_{(A,\Phi)}$ and that we have a canonical way for orienting the
index bundle of  $\cD_{(\wtilde A',\wtilde \Phi')}$, we deduce that the
index bundle of $P$ comes with an orientation.

Since $F$ is invertible, it makes sense to talk about the differential
operator $f_*P := F\circ
f\circ F^{-1}$ on $M'_\tau\times N_\tau$.
We will deduce	an orientation for the index bundle of $R$ thanks to the
next lemma.
\begin{lemma}
\label{lemmaex}
The  operator $F_0 P F^{-1}_0 - R$ is linear of order $0$ and  has
compactly supported
coefficients.
\end{lemma}
\begin{proof}
The fact that $\chi_\tau^2+\mu_\tau^ 2 =1$ implies
$$ F_0^{-1} =
\left (\begin{array}{cc}
\chi_\tau & -\mu_\tau \\
\mu_\tau & \chi_\tau
\end{array} \right ),
$$
hence
\begin{align*}
F_0 P F_0^{-1} & = \left (\begin{array}{cc}
\chi_\tau & \mu_\tau \\
-\mu_\tau & \chi_\tau
\end{array} \right )
\left (\begin{array}{cc}
\cD_{(A,\Phi)} & 0 \\
0& \wtilde \cD_{(\wtilde A',\wtilde \Phi')}
\end{array} \right )
\left (\begin{array}{cc}
\chi_\tau & -\mu_\tau \\
\mu_\tau & \chi_\tau
\end{array} \right ) \\
&= \left (\begin{array}{cc}
\chi_\tau \cD_{( A, \Phi)} \chi_\tau  + \mu_\tau \wtilde  \cD_{(\wtilde
A',\wtilde
  \Phi')} \mu_\tau
& - \chi_\tau \cD_{( A, \Phi)} \mu_\tau  + \mu_\tau \wtilde  \cD_{(\wtilde
A',\wtilde
  \Phi')} \chi_\tau \\
- \mu_\tau \cD_{( A, \Phi)} \chi_\tau  + \chi_\tau \wtilde  \cD_{(\wtilde
A',\wtilde
  \Phi')} \mu_\tau & \mu_\tau \cD_{( A, \Phi)} \mu_\tau  + \chi_\tau
  \wtilde
\cD_{(\wtilde A',\wtilde
  \Phi')} \chi_\tau
\end{array} \right )
\end{align*}
We apply the Leibniz rule together with the fact that $d\chi_\tau$ and
  $d\mu_\tau$ are supported in the annulus $\chi_\tau\mu_\tau\neq
  0$. Moreover, the operators $ \cD_{(
  A,
  \Phi)}$,  $ \cD_{( A',
  \Phi')}$,   $ \wtilde\cD_{(\wtilde A,
\wtilde	\Phi)}$ and   $\wtilde
\cD_{(\wtilde A',\wtilde
  \Phi')}$ differ only by some $0$--order terms on the annulus where
  $\chi_\tau\mu_\tau\neq 0$. We deduce that the above matrix is equal
  to
$$\left (\begin{array}{cc}
\cD_{(A',\Phi')} & 0 \\
0& \wtilde \cD_{(\wtilde A,\wtilde \Phi)}
\end{array} \right ) +	K
$$
where $K$ is a matrix with coefficients of order $0$, compactly
supported in the annulus $\chi_\tau\mu_\tau\neq 0$, and the lemma
is proved.
\end{proof}

It follows that the operators $f_*P$ and $R$ are homotopic by
\fullref{lemmaex}.
 So an orientation of $\Omega(P)$
induces an orientation of $\Omega(R)$. Remark that this definition is
consistent since $U$ is path-connected.
The index bundle of $\cD_{(\wtilde A,\wtilde \Phi)}$ is
canonically oriented on $N_\tau$, so the orientation of $\Omega(P)$
induces  finally an orientation
of $\Omega ( \cD_{(A',\Phi')})$.

\subsubsection{Gluing}
We restrict now the map $f$ to $\cM_{\varpi_\tau}(M_\tau)\times
[B,\Psi]$. We have clearly $f = \sharp \times e$, where $e$ is a map
from  $\cM_{\varpi_\tau}$ into $\cC/\cG(N_\tau)$. Because
the solutions of Seiberg--Witten equations have uniform exponential
decay, we can assume that the image of $e$ is contained in an
arbitrarily small neighborhood of $[B,\Psi]$. Hence $e$ is homotopic to
the constant map $e([A,\Phi])= [B,\Psi]$.

On the other hand, the pregluing $\sharp$ is isotopic to the gluing
map $\gG$ (see \fullref{rqisotopy}). It follows that
$f(\cM_{\varpi_\tau}(M_\tau)\times [B',\Psi'])$ is
isotopic to $\cM_{\varpi_\tau}(M'_\tau)\times [B,\Psi]$.

The first manifold is oriented by the choice of an orientation for
$\Omega (f_*P) = f_*\Omega(P)$ while the latter is oriented by
orienting $\Omega(R)$. Since those orientations are compatible and the
two manifolds are isotopic, it follows that the gluing map is
orientation preserving.
\begin{lemma}
The gluing map $\gG$ is orientation preserving.
\end{lemma}

\section{Proof of the main theorems}

\subsection[The proof of \ref{theoexcision}]
{The proof of \fullref{theoexcision}}
Let $\overZ$ be a special symplectic cobordism from$(Y,\xi)$ to
$(Y',\xi')$.
Since it is \emph{special}, there is a collar neighborhood
$[T_0,T_1)\times Y$ of $Y\subset \partial  \overZ$ and a
  contact form $\eta$ on $Y$ such that	the symplectic form $\omega$
  is given by
$\frac 12 d (t^2\eta)$. We can glue a sharp cone on the boundary $Y$
  by extending the collar neighborhood into $(0,T_1)\times Y$ together
  with its symplectic form.

The next argument now follows closely \cite[Lemma 4.1]{KM}.
Let $(0,1]\times Y'$ be a collar neighborhood of $Y'\subset \overZ$
and $f(t)$ be a smooth increasing function equals to $0$ on
$(0,1/2]$ and which tends to infinity as $t$ goes to $1$.
We perturb the symplectic form to
$$  \omega_Z = \omega + \frac 12 d(f^2 \eta'),
$$
where $\eta'$ is a contact form on $Y'$. Now, we have a noncompact
manifold $Z= ((0,T_0)\times Y) \cup (\overZ\setminus Y')$ together
with a symplectic form $\omega_Z$. We define a function $\sigma_Z$
on $Z$ by
\begin{itemize}
\item $\sigma(t,y)=t$ on $(0,T_0)\times Y$,
\item  $\sigma(t,y')=f(t)$ on $f^{-1}[T_0+1,\infty)\times Y'$,
\end{itemize}
and we extend $\sigma$ to the remaining part of $Z$ with the condition
that $T_0\leq \sigma \leq T_0+1$.

With a suitable choice of
compatible almost complex structure $J_Z$, the data
$$(Z,\omega_Z, J_Z,\sigma_Z)$$ is  an AFAK end (see \fullref{dfncob}).
Notice that the condition~\eqref{assumption} is verified as a direct
consequence of the assumption~\eqref{assumption2} for special cobordims.

The piece $(0,T_0)\times Y\subset Z$ has now a structure of almost
K\"ahler cone. We extend it to the infinite almost K\"ahler cone
$Z'= (0,\infty)\times Y$ as defined in \fullref{subsubkahler}.

Let $\overM$ be a compact manifold whose boundary $Y$ is endowed with
an element $\gs\in \spinc(\overM,\xi)$. Put $\overM'
=\overM\cup \overZ$.
We apply the construction of \fullref{secglueafak}
$$ \overM, Z \leadsto M_\tau \quad \mbox{and}\quad  \overM,
Z' \leadsto M'_\tau.
$$
The  moduli space of Seiberg--Witten equations on
$M'_\tau$ leads to the invariant $\sw_{\overM,\xi}(\gs)$ while the
moduli space on $M_\tau$ leads to $\sw_{\overM',\xi'}(j(\gs))$
(beware of the notation).
On the other hand, the gluing
\fullref{theogluing}   says that
for $\tau$ large enough,a suitable choice of perturbation and
given compatible consistent orientations, the
Seiberg--Witten moduli spaces for $M_\tau$ and $M'_\tau$ are generic
and orientation preserving diffeomorphic
via the gluing map $\gG$. Therefore   $\sw_{\overM,\xi'}(\gs)=
\pm \sw_{\overM',\xi'}(j(\gs)) $ as does the more precise
formulation  hence \fullref{theoexcision} holds.

\subsection{Surgery and monopoles}
Let $\overM^4$, $(Y,\xi)$  and $K$ a Legendrian knot be as in
\fullref{theo:main}.
\begin{nopartheorem}
\begin{result}[\cite{W}]
Let $\overZ$ be the cobordism between $Y$ and the manifold $Y'$
obtained by a $1$--handle surgery
on $Y$, or, from a $2$--handle surgery along $K$ with framing coefficient
relative to the canonical framing is $-1$.
Then $Y'$ is endowed with a contact structure $\xi'$ and  $\overZ$
has a structure of symplectic cobordism between the contact manifolds
$Y$ and $Y'$.  Moreover
the symplectic form of $\overZ$ is equal to a symplectization of
the contact structures	in a collar neighborhood
of $Y$ and $Y'$.
\end{result}
\end{nopartheorem}
This result of Weinstein  generalizes to the symplectic
category some techniques developed by Eliashberg in the case of
Stein domains~\cite{E}.
The Weinstein surgeries are particular
special symplectic cobordisms. The only thing that has to be checked
is the property~\eqref{assumption2}, namely that $i^*\co  H^1(\overZ,
Y')\to H^1(Y)$ is the zero map: take a $1$--cycle $\cS$ in $Y$. We
can always perturb $\cS$  so that it avoids a neighborhood of the
locus along
which the surgery is performed. Then $\cS$ is homologous to a $1$--cycle
$\cS'$ in
$Y'$. Let $\Omega$ be a cohomology class in $ H^1(\overZ ,Y')$.
Then $\Omega\cdot \cS = \Omega\cdot \cS' = 0$ hence $i^*\Omega = 0$.

Before we  prove \fullref{theo:main}, notice that
together with \fullref{theoexcision},
they imply an improved version of \fullref{mycortight}.
\begin{cor}
Let $Y$ be a manifold with a contact structure $\xi$, filled by a
manifold $\overM$ with
 $\sw_{\overM,\xi}\not\equiv 0$. Then, any contact manifold $(Y',\xi')$
 constructed by
applying a sequence of Weinstein surgeries to $(Y,\xi)$ is also tight.
\end{cor}
\begin{remark*}
 It is not true in general that
the category of tight contact structures is stable under Weinstein
surgery (see for instance \cite{H}). However the theorem of
Gromov--Eliashberg asserts that the category of tight
contact structures that are
weakly symplectically semi-fillable is stable under Weinstein
surgery. The above
corollary gives a new proof of this fact, since the
Seiberg--Witten invariant of a weakly symplectically semi-fillable
manifold is nonvanishing~\cite{KM}, together with  as more general
statement.

There are examples of tight nonsymplectically semi-fillable contact
structure \cite{LS}. However, to the best of our knowledge, these
examples always have a vanishing gauge theoretic (Seiberg--Witten or
Ozsv\'ath--Szab\'{o}) invariant.
Hence, it would be most interesting to find some new tight contact
structures that are not symplectically semi-fillable although fillable
by a manifold with a nonvanishing Seiberg--Witten invariant.
\end{remark*}

\subsubsection{Invariants of connected sums}
Let $\overM_1$ and $\overM_2$ be two compact manifolds with
contact boundaries $(Y_1,\xi_1)$ and $(Y_2,\xi_2)$. The connected sum
$Y_1\sharp Y_2$ is a particular case of $1$--handle surgery and
therefore carries a contact structure $\xi_1\sharp\xi_2$ by
Weinstein's result. Using the identification
 $\spinc(\overM_1,\xi_1) \times \spinc(\overM_2,\xi_2)
\simeq
\spinc(\overM_1\sqcup \overM_2,\xi_1\sqcup \xi_2)$ and
\fullref{theoexcision}, we deduce that
$$
\sw_{\overM_1,\xi_1}\sw_{\overM_2,\xi_2} =  \sw_{\overM_1\sharp \overM_2,
\xi_1\sharp\xi_2} \circ j.
$$
Similarly to  the compact case, one can easily prove that the
Seiberg--Witten invariants must vanish if the connected sum is
performed in the interior of the manifolds $\overM_j$.
The above  identity shows that Seiberg--Witten invariants behave in
utterly differently way for connected sum at the boundary.

\subsubsection[The proof of \ref{theo:main}]
{The proof of \fullref{theo:main}}

The proof will be based on the adjunction inequality.
\begin{prop}[Adjunction inequality]
\label{adjin}
Let $\overM^4$ be a manifold with a contact boundary $(Y,\xi)$
and an
element $(\gs,h)\in \spinc(\overM,\xi)$ such that
$$\sw_{\overM, \xi}(\gs,h)\neq 0.$$
 Then every  closed surfaces $ \Sigma\subset
\overM$  with $[\Sigma]^2=0$ and genus at least  $1$
verifies
$$|c_1(\gs)\cdot
 \Sigma | \leq	- \chi( \Sigma).
$$
\end{prop}
\begin{proof}
We refer the reader to	\cite{KM2} and check that the argument	 proving
the adjunction
inequality  on a
compact manifold  applies in the same way in our case.
\end{proof}

We return to the assumptions and  notation of Theorem A.
Let $\overZ$ be a special symplectic cobordism obtained by
performing Weinstein surgery along $K$.
We have a new manifold with contact boundary $(Y',\xi')$
$$ \overM' =\overM\cup_Y \overZ.
$$
Let $(\gs,h)$ be an element of $\spinc(\overM,\xi)$ such that
$\sw_{\overM,\xi}(\gs,h)\neq 0$ and let  $(\gtt,k)=j(\gs,h)\in
\spinc(\overM',\xi')$.
Let $\Sigma'$ be the closed connected surface obtained as the
union of $\Sigma$ and the core the $2$--handle added along $K$.
 Then $\Sigma'$ verifies
$$ \chi(\Sigma')=\chi(\Sigma)+1,\quad [\Sigma']\cdot[\Sigma']
=\tb(K,\sigma)-1,\quad
\langle c_1(\gtt), [\Sigma']\rangle = r(K,\sigma,\gs,h).
$$
 If $[\Sigma']^2=0$  and $\chi(\Sigma')\ge 0$
 we are done by the adjunction
 inequality.  We can reduce to this case by some standard tricks.
We first make the assumption that $\tb(K,\sigma)\geq 1$ so that
 $[\Sigma']^2\geq 0$.
The case of equality  is equivalent to the fact  that
he normal bundle of $\Sigma'$ is trivial. To  make it trivial,
we blow up $\tb(K,\sigma) - 1$ points on the core of the $2$--handle
in the interior of $ \overZ$ and denote by $E_j$ the exceptional
divisors with selfintersection $-1$. The blow-up
$\pi\co \what Z\to\overZ$  has still a structure of special symplectic
cobordism between $(Y,\xi)$ and $(Y',\xi')$.
We denote $\what M' = \overM\cup \what Z$;
similarly to $\overM'$, it is endowed with an element $(\what \gtt,k)$ of
$\spinc(\what M',\xi')$ deduced from $(\gs,h)$.
 The adjunction formula shows that the relation between
$\gtt$,  $\what\gtt $ and the proper transform $\what
 \Sigma$ of $\Sigma'$ are
$$ [\what\Sigma]= \pi^* [\Sigma'] - \sum_{j=1}^{\tb(K)-1}
E_j,\quad
  c_1(\what \gtt) = c_1(\pi^*\gtt) +
 \sum_{j=1}^{\tb(K)-1} E_j,
$$
where $E_j$ are Poincar\'e dual to the exceptional divisors (of
self-intersection $-1$).
Then
$$ [\what \Sigma]\cdot  [\what \Sigma] = 0, \quad
c_1(\what\gs)\cdot [\what \Sigma] =
r(K,\sigma) + \tb(K) - 1, \quad \chi(\what \Sigma)=
\chi(\Sigma')=\chi(\Sigma) +1,
$$
and the normal bundle of $\what\Sigma$ is trivial.

\begin{figure}[ht!]
\begin{center}
\includegraphics{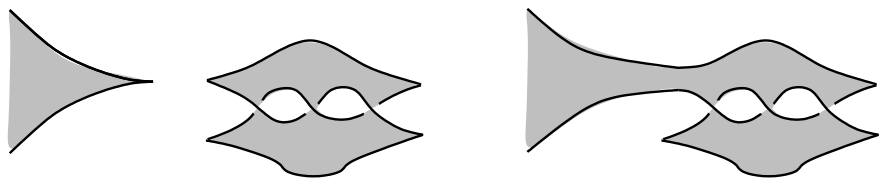}
\caption{\label{rht}Connected sum with a  right handed	Legendrian
trefoil knot}
\end{center}
\end{figure}

We have proved the theorem in the case where $\tb(K,\sigma )\geq 1$ and
$\chi(\Sigma)\leq -1$. We now show that the theorem holds in general.
For this purpose, we modify $K$ and
$\Sigma$ as follow: let $P$ be a point on $K$. Then $P$ has a
neighborhood in $Y$ contactomorphic to $\RR^3$ with its standard
contact structure.
We make a connected sum between $K$
and  a right handed Legendrian trefoil knot at $P$ as shown in
\fullref{rht}; accordingly,
add two $1$--handles to $\Sigma$ by forming the connected sum at cusps
as suggested by the gray region of \fullref{rht}. This
operation
decreases $\chi(\Sigma)$ by two,  increases
$\tb(K,\sigma)$ by $2$	but does not change the
rotation number $r(K,\sigma,\gs,h)$.

Overall the quantity $\chi(\Sigma)+\tb(K,\sigma) +|r(K,\sigma,\gs,h)|$
remains
unchanged.
 So we can always assume that $\tb(K,\sigma)\geq 1$ and $\chi(\Sigma)<0$
and Theorem A is proved.

\bibliographystyle{gtart}
\bibliography{link}

\end{document}

%% file: figs/MZ.pstex_t
\begin{picture}(0,0)%
\includegraphics{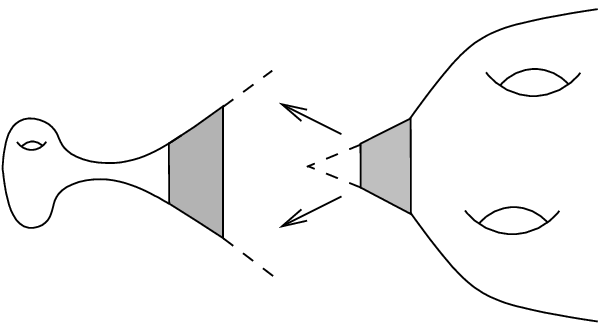}%
\end{picture}%
\setlength{\unitlength}{3947sp}%
\begingroup\makeatletter\ifx\SetFigFont\undefined%
\gdef\SetFigFont#1#2#3#4#5{%
  \reset@font\fontsize{#1}{#2pt}%
  \fontfamily{#3}\fontseries{#4}\fontshape{#5}%
  \selectfont}%
\fi\endgroup%
\begin{picture}(2880,1524)(88,-748)
\put(2431, 14){\makebox(0,0)[lb]{\smash{\SetFigFont{12}{14.4}{\rmdefault}{\mddefault}{\updefault}\tiny $Z$}}}
\put(170,-90){\makebox(0,0)[lb]{\smash{\SetFigFont{12}{14.4}{\rmdefault}{\mddefault}{\updefault}\tiny $M$}}}
\end{picture}

%% file: main_a.bbl
\begin{thebibliography}{}
\providecommand\bibmarginpar{\leavevmode\marginpar}
\def\urlstyle#1{{\tt #1}}

\bibitem{AM}
\textbf{S Akbulut}, \textbf{R Matveyev}, \emph{Exotic structures and adjunction
  inequality}, Turkish J. Math. 21 (1997) 47--53 \xox{MR}{1456158}

\bibitem{Etight}
\textbf{Y Eliashberg}, \emph{Filling by holomorphic discs and its
  applications}, from: ``Geometry of low-dimensional manifolds, 2 (Durham,
  1989)'', London Math. Soc. Lecture Note Ser. 151, Cambridge Univ. Press,
  Cambridge (1990)  45--67 \xox{MR}{1171908}

\bibitem{E}
\textbf{Y Eliashberg}, \href{http://dx.doi.org/10.1142/S0129167X90000034}
  {\emph{Topological characterization of {S}tein manifolds of dimension > 2}},
  Internat. J. Math. 1 (1990) 29--46 \xox{MR}{1044658}

\bibitem{Efill}
\textbf{Y Eliashberg}, \href{http://dx.doi.org/10.2140/gt.2004.8.277} {\emph{A
  few remarks about symplectic filling}}, Geom. Topol. 8 (2004) 277--293
  \xox{MR}{2023279}

\bibitem{EH}
\textbf{J\,B Etnyre}, \textbf{K Honda},
  \href{http://dx.doi.org/10.1007/s002220100204} {\emph{Tight contact
  structures with no symplectic fillings}}, Invent. Math. 148 (2002) 609--626
  \xox{MR}{1908061}

\bibitem{G}
\textbf{P Gauduchon}, \emph{Hermitian connections and {D}irac operators}, Boll.
  Un. Mat. Ital. B $(7)$ 11 (1997) 257--288 \xox{MR}{1456265}

\bibitem{Go}
\textbf{R\,E Gompf}, \emph{Handlebody construction of {S}tein surfaces}, Ann.
  of Math. $(2)$ 148 (1998) 619--693

\bibitem{H}
\textbf{K Honda},
  \href{http://projecteuclid.org/getRecord?id=euclid.dmj/1085598176}
  {\emph{Gluing tight contact structures}}, Duke Math. J. 115 (2002) 435--478
  \xox{MR}{1940409}

\bibitem{KM2}
\textbf{P\,B Kronheimer}, \textbf{T\,S Mrowka}, \emph{The genus of embedded
  surfaces in the projective plane}, Math. Res. Lett. 1 (1994) 797--808
  \xox{MR}{1306022}

\bibitem{KM}
\textbf{P\,B Kronheimer}, \textbf{T\,S Mrowka},
  \href{http://dx.doi.org/10.1007/s002220050183} {\emph{Monopoles and contact
  structures}}, Invent. Math. 130 (1997) 209--255 \xox{MR}{1474156}

\bibitem{LM}
\textbf{P Lisca}, \textbf{G Mati{\'c}},
  \href{http://dx.doi.org/10.1016/S0166-8641(97)00198-3} {\emph{Stein
  4--manifolds with boundary and contact structures}}, Topology Appl. 88 (1998)
  55--66 \xox{MR}{1634563}

\bibitem{LS1}
\textbf{P Lisca}, \textbf{A\,I Stipsicz},
  \href{http://dx.doi.org/10.2140/gt.2003.7.1055} {\emph{An infinite family of
  tight, not semi-fillable contact three-manifolds}}, Geom. Topol. 7 (2003)
  1055--1073 \xox{MR}{2026538}

\bibitem{LS}
\textbf{P Lisca}, \textbf{A\,I Stipsicz},
  \href{http://dx.doi.org/10.1007/s00208-003-0483-0} {\emph{Tight, not
  semi-fillable contact circle bundles}}, Math. Ann. 328 (2004) 285--298
  \xox{MR}{2030378}

\bibitem{Taub}
\textbf{C\,H Taubes}, \emph{The {S}eiberg--{W}itten invariants and symplectic
  forms}, Math. Res. Lett. 1 (1994) 809--822 \xox{MR}{1306023}

\bibitem{W}
\textbf{A Weinstein}, \emph{Contact surgery and symplectic handlebodies},
  Hokkaido Math. J. 20 (1991) 241--251 \xox{MR}{1114405}

\end{thebibliography}
